\documentclass[12pt]{article}

\usepackage[margin=1in]{geometry}
\usepackage{amsmath,amssymb,amsfonts,amsthm,mathtools}
\usepackage{enumitem}
\usepackage{microtype}
\usepackage{booktabs}
\usepackage{graphicx}
\allowdisplaybreaks[1]
\usepackage[authoryear,round]{natbib}

\usepackage[
colorlinks=true,
citecolor=blue,
linkcolor=blue,
urlcolor=blue,
pagebackref=true
]{hyperref}

\hypersetup{pdftitle={Asymptotically Bayes-Optimal Sparse Multiple Testing Under One-Factor Dependence},pdfauthor={Prasenjit Ghosh}}

\newtheorem{theorem}{Theorem}
\newtheorem{lemma}{Lemma}
\newtheorem{corollary}{Corollary}

\newtheorem{remark}{Remark}
\newtheorem{definition}{Definition}

\newcommand{\Phibar}{\overline{\Phi}}

\title{\textbf{Asymptotic Bayes Optimality of\\
	Sparse Multiple Testing Procedures\\
	under One-Factor Dependence}}
\author{
Prasenjit Ghosh\thanks{
	Affiliations when this work was completed: Department of Statistics,
	Texas A\&M University, College Station, TX 77843, USA; and Applied
	Statistics Unit, Indian Statistical Institute, Kolkata 700108, India.
	Email: \texttt{prasenjit.tamu@gmail.com}.
	}
	\hspace{1em}
	Arijit Chakrabarti\thanks{
		Applied Statistics Unit, Indian Statistical Institute,
		Kolkata 700108, India.
		Email: \texttt{arc@isical.ac.in}.
	}
}

\date{}

\begin{document}
\maketitle

\begin{abstract}
	
In this paper, we study sparse multiple testing under a Gaussian two-groups model with known one-factor covariance and centered Gaussian signals. We derive a first-order risk expansion for the observed-data Bayes Oracle, including a finite-dimensional loading average for heterogeneous loadings and an explicit limiting constant for homogeneous loadings. Under the stated factor-estimation and signal-calibration conditions, the homogeneous one-factor Oracle has strictly smaller leading risk than its independent-coordinate counterpart; observing the latent factor yields no further first-order reduction.

We establish asymptotic Bayes optimality under sparsity (ABOS) for the Maximum Residual Down (MRD) step-down procedure with Gavrilov--Benjamini--Sarkar (GBS) stagewise calibration (MRD--GBS). Uniform control of residual errors along the rejection path permits comparison with a fixed factor-known rule without requiring the rejection paths to coincide. Fixed-residual versions of GBS and Benjamini--Hochberg (BH) attain the same leading risk. We also prove ABOS for the Bayesian Step-Down (BSD) rule when supplied with the true mixture, signal variance, covariance, and loss parameters; this is an Oracle-parameterized result. When the expected signal count grows faster than the logarithm of the reciprocal signal proportion, BSD and the observed-data Oracle disagree on a number of coordinates negligible in probability relative to the expected signal count, explaining BSD's near-Oracle behavior in companion simulations. The results cover diverging and bounded positive expected signal counts, including every fixed polynomial sparsity rate up to the bounded-count regime, under the stated calibrations. A fixed-dimension simulation illustrates absolute-risk behavior but does not establish ABOS.	

\end{abstract}



\section{Introduction}
\label{sec:introduction}

Sparse signal discovery is a fundamental problem in modern statistical inference. In genomics, bioinformatics, neuroimaging, astronomy, finance, and many other fields, one seeks to identify a relatively small number of active signals among a large collection of null effects. This high-dimensional setting has made multiple testing indispensable and has generated a rich literature on procedures controlling global measures of error. The family-wise error rate and, especially, the false discovery rate (FDR) have motivated step-up, step-down, and adaptive methods; see, among many others, \citet{BH1995,BenjaminiLiu1999,BY2001,STO_2002,StoreyTaylorSiegmund2004,BKY2006,BlanchardRoquain2009,RomanoShaikhWolf2008,Sarkar2008}. These criteria provide important guarantees for simultaneous inference. Yet a procedure can satisfy an error-rate target without accurately recovering the signal configuration, and procedures with similar FDR can differ in the number of signals they miss.



The decision-theoretic optimality of multiple-testing procedures has received comparatively less attention. One line of work studies admissibility: whether a procedure can be uniformly improved in risk. Cohen and collaborators developed this perspective for multiple-endpoint testing, establishing inadmissibility results for step-up procedures and proposing covariance-aware alternatives, including the MRD step-down procedure \citep{CohenSackrowitz2005, COHEN_SACK_2007, COHEN_KOL_SACK_2007, COHEN_SACK_2008, COHEN_SACK_XU_2009}. Admissibility is a useful safeguard, but it does not quantify how close a procedure comes to the smallest attainable risk. Asymptotic Bayes optimality under sparsity (ABOS), introduced by \citet{BCFG2011}, provides such a comparison: it asks whether a procedure attains, to first order, the Bayes risk of the optimal rule under a specified sparse model and loss. It thus provides a criterion for assessing whether a computationally feasible procedure gives up anything in leading-order risk relative to the Bayes Oracle. In sparse problems, matching only the rate of the Oracle risk is not enough: the leading constant also matters. ABOS complements error-rate analyses by measuring the combined cost of false rejections and missed signals. FDR control and ABOS are distinct properties: FDR control concerns the expected fraction of false rejections under a specified sampling model, whereas ABOS compares prior-averaged classification risk with the Bayes-Oracle risk. Neither property implies the other. For example, the procedure that makes no rejections has FDR zero but misses every true signal, illustrating that FDR control alone does not guarantee useful signal recovery. In \cite{GC2016}, ABOS was established for decision rules induced by one-group continuous shrinkage priors, without an accompanying finite-sample FDR-control guarantee.

In the independent-coordinate setting, ABOS and sharp asymptotic minimaxity have been established for several thresholding and multiple-testing procedures; see, for example, \citet{Abramovich2006, BGT2008, BCFG2011, NR2012, FrommletBogdan2013, GC2016, PaulChakrabarti2025, QinGhosh2025, Ghosh2026GBSABOS}, and references therein. Related Bayes-risk analyses of the horseshoe and broader classes of one-group shrinkage rules establish Oracle-risk approximations up to multiplicative constants close to the Oracle constant, rather than the exact ratio-one conclusion defining ABOS \citep{DG2013,GTGC2015}. More recent work gives sharp risk characterizations and detection boundaries, including results for the Benjamini--Hochberg (BH) \citep{BH1995}, Gavrilov--Benjamini--Sarkar \citep{GBS2009} and $\ell$-value based procedures and for procedures induced by global--local shrinkage priors \citep{AbrahamCastilloRoquain2024, PaulGhoshChakrabarti2025, GhoshGBSMinimax2026}, only to name a few. ABOS and sharp asymptotic minimaxity are distinct: the former compares integrated risk under a specified prior with Bayes-Oracle risk, while the latter compares worst-case risk over a class of sparse signals with the minimax benchmark. Both lines of work help identify when computationally feasible procedures attain the appropriate first-order risk under sparsity.

Much less is understood when the observations are dependent. Correlation can cause multiple marginal statistics to appear large because of common variation, and marginal procedures may then misrepresent coordinate-specific evidence. But dependence can also be informative: the other coordinates can help estimate the common component and reveal signals that are difficult to distinguish from marginal observations alone. Dependence is therefore not only a complication to accommodate; when incorporated into the testing procedure, it can be a source of inferential gain. This dual role has motivated methods that control error under dependence and methods that incorporate dependence into the construction of the test statistics. Examples include empirical Bayes analyses of correlated statistics, factor-adjusted methods, dependence-aware compound decision rules, and procedures for sparse signal detection; see, for example, \citet{Efron2007, Friguet2009, LeekStorey2008, SUN_CAI_2009, HallJin2010, FanHanGu2012}. Recent shifted-BH and e-value procedures extend error-rate guarantees under their stated settings \citep{SarkarZhang2025, WangRamdas2022, LiZhang2025}. These contributions leave open the decision-theoretic question central here: can a procedure exploit the covariance structure to improve signal recovery and attain the first-order Bayes risk under the joint model?


Motivated by this, we consider in this paper the simultaneous testing of $m$ hypothesis
$$H_{0i}:\theta_i=0 \ \textrm{against} \ H_{Ai}:\theta_i\ne0,\quad i=1,\ldots,m,$$ under the sparse Gaussian two-groups model
\[
 \boldsymbol{X}\mid\boldsymbol{\theta}\sim N_m(\boldsymbol{\theta},\boldsymbol{\Sigma}_m),\qquad
 \theta_i\mid\nu_i=0\sim \delta_{\{0\}},\qquad
 \theta_i\mid\nu_i=1\sim N(0,\psi_m^2).
\]
Here, $\nu_i=0$ indicates that $H_{0i}$ is true, whereas $\nu_i=1$ indicates that $H_{0i}$ is false, and \(\delta_{\{0\}}\) denotes the Dirac delta measure degenerate at $0$. We assign independent Bernoulli$(p_m)$ priors to the indicators, with $p_m\to0$, and take the known covariance matrix $\boldsymbol{\Sigma}_m$ to be positive definite. Let \(\boldsymbol{\varphi}(\boldsymbol{X}) = (\varphi_1(\boldsymbol{X}),\ldots,\varphi_m(\boldsymbol{X}))\in\{0,1\}^m\) denote a multiple testing procedure, with \(\varphi_i=1\) indicating the rejection of the $i$-th null hypothesis. If $V_m$ and $T_m$ are the numbers of false rejections and missed signals, respectively, we measure the performance of \(\boldsymbol{\varphi}(\boldsymbol{X})\) through the additive Bayes risk
\[
 \mathcal R_m(\boldsymbol{\varphi})=\delta_{0,m}\,\mathbb E V_m+
 \delta_{A,m}\,\mathbb E T_m,
 \qquad \delta_{0,m},\delta_{A,m}>0.
\]
The expectation averages over the observations and the sparse configuration. The Bayes Oracle minimizes this risk using the entire observed vector and the specified model and loss parameters. A procedure is ABOS when the ratio of its Bayes risk to this observed-data Oracle risk tends to one. This formulation accounts jointly for both types of classification error and makes explicit which information the benchmark is allowed to use.

Under dependence, that Oracle is difficult to characterize. Although the prior indicators are independent, the likelihood couples the coordinates, so posterior inclusion probabilities are joint functions of the full observation vector. For each coordinate, evaluating its posterior inclusion probability requires summing over the $2^m$ latent configurations. As noted by \citet{XCML_2011}, exact implementation of the full Bayes Oracle therefore has computational complexity of order $O(m2^m)$. The Oracle is consequently an important ideal benchmark, but exact computation quickly becomes infeasible as $m$ grows. This difficulty motivates the Bayesian Step-Down (BSD) procedure, first proposed in our 2015 arXiv preprint \citep{GC_BSD2015}. The companion BSD study \citep{GhoshChakrabartiBSD2026} develops its sequential, posterior-guided model-pursuit interpretation: rather than summing over the full configuration space at every decision, it sequentially extends an active set using local posterior evidence while incorporating the covariance structure. In the one-factor simulations reported there, BSD often has Bayes risk and signal-recovery performance remarkably close to the numerically evaluated Bayes Oracle, even at relatively smaller dimensions such as $m=30$. The Oracle analysis and the empirical behavior of BSD thus provide a concrete setting in which it is natural to ask how well a feasible dependence-aware method can recover the risk of the full-vector decision rule.

A related empirical observation motivates the present analysis more directly. In the same one-factor comparisons, the frequentist MRD--GBS procedure often tracks BSD and the Bayes Oracle. MRD--GBS combines two components developed in the companion MRD--GBS manuscript \citep{GhoshChakrabartiMRDGBS2026}. The Maximum Residual Down (MRD) procedure of \citet{COHEN_SACK_XU_2009} uses the covariance matrix to construct standardized residual evidence on the currently active coordinates, selects the largest residual, and recomputes the residual system after each rejection in a step-down manner. The original MRD formulation requires a sequence of stagewise critical values, which were chosen numerically for particular dependence settings rather than through a unified calibration principle. The Gavrilov--Benjamini--Sarkar (GBS) constants \citep{GBS2009} supply such a systematic principle: their thresholds become less conservative as the step-down procedure advances, reflecting that hypotheses surviving earlier stages have already passed more stringent comparisons. Under independence, the MRD residuals reduce to the marginal statistics and MRD--GBS reduces exactly to GBS. The independent GBS procedure has ABOS and sharp asymptotic minimaxity results under appropriate sparse-model conditions \citep{Ghosh2026GBSABOS,GhoshGBSMinimax2026}. Those theorems do not cover dependent residuals or a data-selected sequence of active sets. The companion MRD--GBS simulation study
\citep{GhoshChakrabartiMRDGBS2026} reports encouraging empirical
performance across several dependence structures. In the settings
examined, MRD--GBS typically attains high power, empirical FDR close to the nominal level, low false non-discovery rates, and rejection counts close to the expected number of signals. Separately, the one-factor simulations in the companion BSD study
\citep{GhoshChakrabartiBSD2026} compare MRD--GBS and BSD with the
observed-data Bayes Oracle and show that MRD--GBS often tracks both.
These finite-sample findings motivate the present analysis, but do not establish ABOS. We aim to address the theoretical question of whether MRD--GBS remains Bayes-risk efficient when its covariance-adjusted residuals are updated along the data-dependent rejection path.

This calibration question also has a decision-theoretic context. Earlier work shows that several marginal procedures, including BH, can be inadmissible under natural vector-valued losses in dependent problems, whereas MRD has an admissibility foundation and the broader monotone residual-based step-down class includes the GBS-calibrated procedure \citep{CohenSackrowitz2005,GhoshChakrabartiAdmissibility2026}. The BSD rule also connects the two lines of work: at each fixed active set and coordinate, its stagewise posterior-odds score is a strictly increasing function of the squared MRD residual, placing it in the local monotone-transform class studied in the general admissibility analysis. Admissibility and ABOS nevertheless answer different questions: admissibility rules out a uniformly dominating decision rule under the specified vector loss, while ABOS identifies first-order integrated risk relative to a Bayes Oracle. Neither property entails the other. The present paper asks whether the observed near-Oracle behavior of MRD--GBS has a first-order decision-theoretic explanation. The answer requires analyzing both the benchmark and the actual adaptive procedure. We confine our attention to a known one-factor Gaussian covariance model. This provides a structured but genuinely dependent setting in which common variation can be exploited, while conditioning on the latent factor renders the coordinatewise errors independent. This structure permits a tractable first-order analysis without reducing the problem to independence or claiming results for arbitrary covariance matrices. The factor remains unobserved, so the observed-data Oracle must still account for its uncertainty; we distinguish its risk from that of the more informative factor-known Oracle. The one-factor model therefore provides a focused first setting for studying when dependence improves inference and when feasible procedures recover the joint-model Bayes risk. We also distinguish the complete adaptive MRD--GBS algorithm from fixed residualized procedures: after every rejection, MRD--GBS changes the active set and recomputes the residuals, so a result for a fixed collection of residuals does not settle the adaptive case.

There are two related difficulties in establishing ABOS for these
step-down methods under dependence. First, the step-down geometry is
stage-dependent. At each stage the rule orders current evidence, applies
a stage-specific critical value, and updates its statistics after a
rejection. For MRD--GBS, deleting a coordinate changes the active
covariance structure and the residuals. For BSD, the conditional
posterior-odds comparisons change with the active model. Consequently,
neither procedure can be analyzed as a single fixed-threshold rule.
Second, the active set at each stage is selected through the preceding
rejection history. The residuals and posterior scores are therefore
evaluated on random, history-dependent sets. A Gaussian calculation for
a fixed set cannot simply be applied after selection as though that set
had been chosen in advance. The analysis must control the statistics
over the sets that can arise along the rejection path and show that
exceptional paths do not contribute at the leading Bayes-risk scale.
These difficulties interact, but they call for separate control of the
changing step-down statistics and the selection of the active sets.

Our contributions in this paper are multi-fold:
\begin{enumerate}[label=(\roman*),leftmargin=*]
\item We derive the first-order risk of the observed-data Bayes Oracle under homogeneous and heterogeneous signed one-factor loadings. Under the stated signal calibration, a factor-known lower bound and an observed-data upper bound have the same leading constant, so observing the latent factor would not improve risk to first order.
\item We prove ABOS for the complete adaptive MRD--GBS algorithm relative to the observed-data Oracle. The proof controls residual perturbations uniformly over capped active sets visited along the data-dependent path and does not require agreement with a fixed-rule rejection path.
\item We establish ABOS for fixed leave-one-out GLS-residualized GBS and BH procedures, and for BSD with oracle-supplied model and loss inputs. Under $mp_m/\log(1/p_m)\to\infty$, its terminal rejection set also differs from the observed-data Oracle set by $o_P(mp_m)$ coordinates.
\item The results cover diverging expected signal counts and bounded expected counts bounded away from zero, including $p_m\asymp m^{-\beta}$ for each fixed $0<\beta\le1$, under the signal, level, and loss calibrations specified in the theorems. The accompanying fixed-dimension simulations report absolute Bayes risks across eleven sparsity levels and six loading patterns; they illustrate finite-sample behavior but do not establish ABOS.
\end{enumerate}

The fixed residualized-BH result is related to \citet{RoyBhandari2024}, who establish ABOS for BH-type procedures in exchangeable or equicorrelated Gaussian models using statistics centered by the sample mean, under their stated sparsity and calibration assumptions. ABOS for dependence-adjusted BH in the homogeneous equicorrelated case is therefore not new, and we do
not claim otherwise. Our fixed residualized-BH result uses leave-one-out GLS residuals and extends the analysis to heterogeneous signed one-factor loadings. The paper's main additional results concern the observed-data Oracle risk and ABOS for the complete adaptive MRD--GBS procedure and the Oracle-parameterized BSD rule. Ordinary BH applied to
raw marginal $p$-values is distinct from these centered or
GLS-residualized procedures.


The paper is organized as follows. Section~\ref{sec:framework} gives the
multivariate Gaussian model, the exact Bayes decision, the risk, and the ABOS definition. Section~\ref{sec:procedures} describes the
covariance-adaptive procedures and their one-factor forms.
Section~\ref{sec:major-results} presents the Oracle risk expansion,
followed by the MRD--GBS, BSD, and fixed-residual results.
Section~\ref{sec:simulation} reports the finite-sample comparisons, and Section~\ref{sec:discussion} discusses their interpretation, limitations, and directions for further work. Proofs are collected in the appendices.

\section{The multivariate Gaussian testing problem and the Bayes Oracle}
\label{sec:framework}

\subsection{The sampling model and sparse configurations}

We begin with the general Gaussian problem used in the companion BSD
study. Let \(\boldsymbol{X}=(X_1,\ldots,X_m)^T\) be the observed vector and suppose
\begin{equation}
\label{eq:general-sampling}
\boldsymbol{X}\mid\boldsymbol{\theta}\sim N_m(\boldsymbol{\theta},\boldsymbol{\Sigma}_m),
\qquad \boldsymbol{\theta}=(\theta_1,\ldots,\theta_m)^T\in\mathbb R^m,
\end{equation}
where \(\boldsymbol{\Sigma}_m\) is a known positive definite covariance matrix.
The objective is to identify the unknown nonzero coordinates by testing
\(H_{0i}:\theta_i=0\) against \(H_{Ai}:\theta_i\ne0\), simultaneously for
\(i=1,\ldots,m\). Conditional on a fixed mean vector, all observation
error dependence is encoded by \(\boldsymbol{\Sigma}_m\). We introduce Gaussianity
here as the sampling assumption, before imposing any particular covariance
structure. The one-factor restriction used for our risk theorems will be
introduced in Section~\ref{sec:onefactor-specialization}.

Following the two-groups formulation of \citet{BCFG2011}, let
\(\boldsymbol{\nu}=(\nu_1,\ldots,\nu_m)^T\in\mathcal M_m:=\{0,1\}^m\), with
independent \(\nu_i\sim\operatorname{Bernoulli}(p_m)\). Given \(\boldsymbol{\nu}\),
the means are independent and
\begin{equation}
\label{eq:general-prior}
\theta_i\mid\nu_i\sim(1-\nu_i)\delta_{\{0\}}+
\nu_iN(0,\psi_m^2),\qquad \psi_m^2>0.
\end{equation}
Here \(\delta_{\{0\}}\) denotes a point mass degenerate at $0$, distinct from the loss
ratio defined below. Thus the marginal prior on each mean is
\((1-p_m)\delta_{\{0\}}+p_mN(0,\psi_m^2)\). The indicator formulation
identifies a sparse discovery problem with inference on the unknown
configuration \(\boldsymbol{\nu}\). The number of active coordinates is
\(K_m=\sum_i\nu_i\sim\operatorname{Bin}(m,p_m)\); its expectation is
\(mp_m\). The two quantities should not be conflated when the expectation
is bounded.

Put \(\boldsymbol{B}_{\boldsymbol{\nu}}=\operatorname{diag}(\nu_1,\ldots,\nu_m)\) and
\(\boldsymbol{A}_{\boldsymbol{\nu}}=\boldsymbol{\Sigma}_m+\psi_m^2\boldsymbol{B}_{\boldsymbol{\nu}}\). Integrating out the mean vector gives
\begin{equation}
\label{eq:configuration-density}
\boldsymbol{X}\mid\boldsymbol{\nu}\sim N_m(\boldsymbol{0},\boldsymbol{A}_{\boldsymbol{\nu}}),\qquad
f_{\boldsymbol{\nu}}(\boldsymbol{x})=\frac{\exp(-\boldsymbol{x}^T\boldsymbol{A}_{\boldsymbol{\nu}}^{-1}\boldsymbol{x}/2)}
{(2\pi)^{m/2}|\boldsymbol{A}_{\boldsymbol{\nu}}|^{1/2}}.
\end{equation}
The prior mass of a configuration and the marginal density are
\begin{equation}
\label{eq:configuration-mixture}
\pi_m(\boldsymbol{\nu})=p_m^{|\boldsymbol{\nu}|}(1-p_m)^{m-|\boldsymbol{\nu}|},\qquad
f_m^{\boldsymbol{X}}(\boldsymbol{x})=\sum_{\boldsymbol{\nu}\in\mathcal M_m}\pi_m(\boldsymbol{\nu})f_{\boldsymbol{\nu}}(\boldsymbol{x}).
\end{equation}
Thus the model contains \(2^m\) possible active configurations. This
exponential model space is present even though each configuration has a
Gaussian likelihood and the prior indicators are independent.

\subsection{Additive loss and the exact Bayes rule}

Let \(\boldsymbol{\varphi}(\boldsymbol{X})=(\varphi_1(\boldsymbol{X}),\ldots,\varphi_m(\boldsymbol{X}))\in\{0,1\}^m\),
where \(\varphi_i=1\) means rejection of the $i$-th null hypothesis. For positive losses
\(\delta_{0,m}\) and \(\delta_{A,m}\), define
\begin{equation}
\label{eq:general-loss}
L_m(\boldsymbol{\nu},\boldsymbol{\varphi})=\sum_{i=1}^m
\{\delta_{0,m}(1-\nu_i)\varphi_i+
\delta_{A,m}\nu_i(1-\varphi_i)\}.
\end{equation}
Writing \(V_m=\sum_i(1-\nu_i)\varphi_i\) and
\(T_m=\sum_i\nu_i(1-\varphi_i)\), its Bayes risk is
\begin{equation}
\label{eq:general-risk}
\mathcal R_m(\boldsymbol{\varphi})=\mathbb E L_m(\boldsymbol{\nu},\boldsymbol{\varphi})
=\delta_{0,m}\mathbb EV_m+\delta_{A,m}\mathbb ET_m.
\end{equation}
The expectation integrates over the complete hierarchical model. Equal
unit losses give the expected total number of classification errors.
Set \(\delta_m=\delta_{0,m}/\delta_{A,m}\) and
\(f_m=(1-p_m)/p_m\); the scalar \(f_m\) is a prior odds ratio, not the
marginal density \(f_m^{\boldsymbol{X}}\).

The exact posterior inclusion probability is
\begin{equation}
\label{eq:exact-pip}
q_{i,m}(\boldsymbol{x}):=\mathbb P(\nu_i=1\mid \boldsymbol{X}=\boldsymbol{x})
=\frac{\displaystyle\sum_{\boldsymbol{\nu}\in\mathcal M_m:\nu_i=1}
\pi_m(\boldsymbol{\nu})f_{\boldsymbol{\nu}}(\boldsymbol{x})}
{\displaystyle\sum_{\boldsymbol{\nu}\in\mathcal M_m}\pi_m(\boldsymbol{\nu})f_{\boldsymbol{\nu}}(\boldsymbol{x})}.
\end{equation}
The posterior loss of rejection is \(\delta_{0,m}(1-q_{i,m})\), while
that of acceptance is \(\delta_{A,m}q_{i,m}\). Additivity therefore gives
\begin{equation}
\label{eq:exact-bayes-rule}
\varphi_{i,m}^{\mathrm{BO},\boldsymbol{X}}(\boldsymbol{x})
=I\left\{q_{i,m}(\boldsymbol{x})>\frac{\delta_m}{1+\delta_m}\right\}.
\end{equation}
Equivalently, writing \(a_{i,d}(\boldsymbol{x})=\sum_{\boldsymbol{\nu}:\nu_i=d}
\pi_m(\boldsymbol{\nu})f_{\boldsymbol{\nu}}(\boldsymbol{x})\), the rule rejects if
\(a_{i,1}(\boldsymbol{x})/a_{i,0}(\boldsymbol{x})>\delta_m\). This is the genuine observed-data
Oracle: it knows the model parameters, covariance matrix and losses, but
has access only to \(\boldsymbol{X}\). Under general dependence it is not a rule based
solely on the marginal observation \(X_i\).

Let \(\mathcal D_m^{\boldsymbol{X}}\) denote all measurable rules based on \(\boldsymbol{X}\).
The Oracle risk has the exact representations
\begin{align}
\label{eq:exact-bayes-risk}
\mathcal R_m^{\mathrm{BO},\boldsymbol{X}}
&=\inf_{\boldsymbol{\varphi}\in\mathcal D_m^{\boldsymbol{X}}}\mathcal R_m(\boldsymbol{\varphi})
=\sum_{i=1}^m\int_{\mathbb R^m}
\min\{\delta_{0,m}a_{i,0}(\boldsymbol{x}),\delta_{A,m}a_{i,1}(\boldsymbol{x})\}\,d\boldsymbol{x}\\
&=\mathbb E\sum_{i=1}^m
\min\{\delta_{0,m}(1-q_{i,m}(\boldsymbol{X})),\delta_{A,m}q_{i,m}(\boldsymbol{X})\}.
\nonumber
\end{align}
These formulas separate two tasks: evaluating posterior probabilities at
one observed vector, and integrating the optimal conditional loss over all
possible vectors. The second task is required to determine Bayes risk.

\subsection{Why the exact risk is difficult}

A direct implementation of \eqref{eq:exact-pip} sums Gaussian likelihoods
over \(2^m\) configurations, with configuration-specific determinants and
quadratic forms. Matrix identities can reduce the work per configuration,
but do not remove the exponential number of terms from this direct
representation. This is a computational obstacle, rather than a universal
complexity lower bound for every structured covariance model. Estimating
Oracle risk by simulation repeats the calculation for many independently
generated vectors, multiplying the numerical burden.

There is a separate analytical obstacle. The decision boundary in
\eqref{eq:exact-bayes-rule} is defined by a ratio of sums and generally
couples all coordinates. It need not be a rectangular region or a single
threshold on \(|X_i|\). Conditional independence in a latent representation
simplifies the likelihood, but an unobserved factor still has to be
integrated out. Even if one-dimensional quadrature makes individual
Oracle decisions numerically accessible, it does not produce a closed-form
risk expansion or justify exchanging a rare-tail approximation with the
expectation in \eqref{eq:exact-bayes-risk}. This distinction explains why
the encouraging empirical comparisons in the BSD study leave a substantive
theoretical question.

\subsection{Asymptotic Bayes optimality under sparsity}

We use the risk-ratio criterion of \citet{BCFG2011}.
\begin{definition}[Asymptotic Bayes optimality under sparsity]
\label{def:ABOS}
Along a specified parameter sequence with \(m\to\infty\) and \(p_m\to0\),
a sequence of feasible rules \(\boldsymbol{\varphi}_m\in\mathcal D_m^{\boldsymbol{X}}\) is
asymptotically Bayes optimal under sparsity if
\begin{equation}
\label{eq:ABOS-definition}
\frac{\mathcal R_m(\boldsymbol{\varphi}_m)}{\mathcal R_m^{\mathrm{BO},\boldsymbol{X}}}\to1.
\end{equation}
\end{definition}
The comparison uses the same model, losses and observed information in the
numerator and denominator. The definition allows a rule to use specified
model parameters; whether those parameters are known, estimated, or treated
as oracle inputs is part of each theorem's assumptions. ABOS is an
integrated-risk property under a specified sparse mixture. It is distinct from FDR control, a worst-case
minimax assertion, and the empirical observation of a risk ratio close to
one at finitely many dimensions. In particular, ordinary marginal BH,
factor-adjusted BH, and a procedure supplied with the true factor are
different rules and must be distinguished in any Oracle comparison.

\section{From covariance-adaptive search to one-factor testing}
\label{sec:procedures}

\subsection{BSD and MRD--GBS as covariance-adaptive searches}

The companion BSD manuscript \citep{GhoshChakrabartiBSD2026} uses a
posterior-guided model-pursuit interpretation. Rather than sum over the
whole configuration space at every decision, BSD compares the currently
retained configuration with candidate configurations containing one
additional signal. The selected coordinate is then removed from further
competition, while covariance information continues to influence subsequent
evidence. This offers a computationally tractable route through a large
model space. The BSD score and stopping rule are parameterized by the true
mixture, slab, covariance, and loss inputs. Its decisions depend on the
observed vector, but its calibration is not learned from that vector. BSD is
therefore a model-informed procedure, not the exact joint Bayes rule and not
a fully data-driven competitor. Its empirical proximity to the exact
Oracle motivates studying which components of covariance adjustment and
sequential calibration can support optimal risk.

MRD--GBS \citep{GhoshChakrabartiMRDGBS2026} takes a frequentist residual
route. For an active set \(A\), let \(\boldsymbol{\Sigma}_A\) be the principal error
covariance matrix and \(\boldsymbol{\Omega}_A=\boldsymbol{\Sigma}_A^{-1}\). The standardized MRD
statistic for \(i\in A\) is
\begin{equation}
\label{eq:mrd-general-description}
U_{i,A}=\frac{(\boldsymbol{\Omega}_A \boldsymbol{X}_A)_i}{\sqrt{(\boldsymbol{\Omega}_A)_{ii}}}
=\frac{X_i-\boldsymbol{\Sigma}_{i,A\setminus\{i\}}
\boldsymbol{\Sigma}_{A\setminus\{i\}}^{-1}\boldsymbol{X}_{A\setminus\{i\}}}
{\{\Sigma_{ii}-\boldsymbol{\Sigma}_{i,A\setminus\{i\}}
\boldsymbol{\Sigma}_{A\setminus\{i\}}^{-1}\boldsymbol{\Sigma}_{A\setminus\{i\},i}\}^{1/2}}.
\end{equation}
At stage \(k\), starting with \(A_1=\{1,\ldots,m\}\), choose
\(i_k\in\arg\max_{i\in A_k}|U_{i,A_k}|\). Reject \(H_{0i_k}\) if
\begin{equation}
\label{eq:mrd-gbs-description}
2\bar\Phi(|U_{i_k,A_k}|)\le c_k(\alpha_m),\qquad
c_k(\alpha_m)=\frac{k\alpha_m}{m+1-k(1-\alpha_m)}.
\end{equation}
After rejection set \(A_{k+1}=A_k\setminus\{i_k\}\) and recompute the
residuals; otherwise stop. At the final singleton the usual marginal
standardization applies. This combines the residual search of
\citet{COHEN_SACK_XU_2009} with the critical constants of \citet{GBS2009}.
The selected sequence is data dependent. Covariance adaptation alone does
not imply an error-rate or Bayes-risk guarantee for this complete algorithm.

\subsection{Why specialize to one-factor covariance?}
\label{sec:onefactor-specialization}

The general decision problem above supplies the benchmark and motivation.
For a tractable theoretical analysis we now impose
\begin{equation}
\label{eq:onefactor-covariance}
\boldsymbol{\Sigma}_m=\boldsymbol{D}_m+\boldsymbol{\lambda}_m\boldsymbol{\lambda}_m^T,\qquad
\boldsymbol{D}_m=\operatorname{diag}(\sigma_{1,m}^2,\ldots,\sigma_{m,m}^2),
\end{equation}
with unit marginal error variances, \(\sigma_{i,m}^2=1-\lambda_{i,m}^2\).
Equivalently,
\begin{equation}
\label{eq:framework-model}
X_i=\theta_i+\lambda_{i,m}Z+\sigma_{i,m}\varepsilon_i,
\end{equation}
where \(Z,\varepsilon_1,\ldots,\varepsilon_m\) are independent standard
normals, independent of the latent pairs \((\nu_i,\theta_i)\).
The correlation matrix $\boldsymbol{\Sigma}_m$ is known and belongs to this
one-factor class. For a loading representation we assume
\begin{equation}
\label{eq:framework-loadings}
\sup_m\max_i|\lambda_{i,m}|\le\bar\lambda<1,\qquad
\liminf_{m\to\infty}m^{-1}\sum_i\lambda_{i,m}^2\ge\kappa>0.
\end{equation}
The first bound keeps idiosyncratic variances uniformly positive; the second
ensures information for factor estimation. Fixed homogeneous loadings
\(\lambda_{i,m}=\lambda\), \(0<|\lambda|<1\), yield equicorrelation
\(\rho=\lambda^2\). Thus the original error dependence need not vanish.
Heterogeneous loadings permit different magnitudes and signs.

\begin{lemma}[Known one-factor covariance and identifiable loadings]
\label{lem:known-sigma-loadings}
In the normalized one-factor class \eqref{eq:onefactor-covariance},
if at least three loadings are nonzero, the known matrix $\boldsymbol{\Sigma}_m$
determines all loadings up to a simultaneous sign and determines every
$\sigma_{i,m}^2$ uniquely. Condition \eqref{eq:framework-loadings} ensures
that there are at least three nonzero loadings for all sufficiently large
$m$. MRD--GBS, fixed GLS residualized GBS, and residual BH are invariant
under that common sign choice. Consequently their known-loading
formulations can equivalently be implemented from known $\boldsymbol{\Sigma}_m$
within this class, without observing $Z$.
\end{lemma}

\paragraph{Implication for the information assumption.}
The loadings need not be supplied as additional information beyond the
known one-factor covariance. Their representation can be recovered from
its off-diagonal entries, as proved in
Appendix~\ref{app:covariance-sparsity}. This equivalence uses the
nondegeneracy in \eqref{eq:framework-loadings}; an unrestricted one-factor
model with only one or two nonzero loadings need not be identifiable.
A covariance matrix estimated from data is a different setting and is
not covered by this equivalence. In particular, knowing the covariance
does not reveal the realized latent factor.

Given \(Z=z\), set
\(g_{d,i}(x_i\mid z)=\phi_{\sigma_{i,m}^2+d\psi_m^2}
(x_i-\lambda_{i,m}z)\), where \(\phi_v\) is the \(N(0,v)\) density,
and \(h_i=(1-p_m)g_{0,i}+p_mg_{1,i}\). Then
\begin{equation}
\label{eq:main-oracle-quadrature}
q_{i,m}(\boldsymbol{x})=
\frac{\int_{\mathbb R}\phi(z)p_m g_{1,i}(x_i\mid z)
\prod_{j\ne i}h_j(x_j\mid z)\,dz}
{\int_{\mathbb R}\phi(z)\prod_{j=1}^m h_j(x_j\mid z)\,dz}.
\end{equation}
This exact identity sums out the configurations by conditional
factorization and then integrates the factor. It enables numerical
benchmarking without treating the realized factor as observed.

\subsection{Factor-known Oracle and feasible residualization}

If \(Z\) were observed, define
\(Y_i(Z)=(X_i-\lambda_{i,m}Z)/\sigma_{i,m}\) and
\(u_{i,m}=\psi_m^2/\sigma_{i,m}^2\). Its null and signal laws are
\(N(0,1)\) and \(N(0,1+u_{i,m})\). The factor-known Oracle rejects when
\(|Y_i(Z)|>c_{i,m}^{O}\), where, for all sufficiently large dimensions in
the calibrated regimes,
\begin{equation}
\label{eq:main-oracle-threshold}
(c_{i,m}^{O})^2=(1+u_{i,m}^{-1})
\log\{(1+u_{i,m})f_m^2\delta_m^2\}.
\end{equation}
Its risk is explicit:
\begin{equation}
\label{eq:main-factor-known-risk}
\mathcal R_m^{\mathrm{BO},Z}=
\sum_i\left[\delta_{0,m}(1-p_m)2\bar\Phi(c_{i,m}^{O})+
\delta_{A,m}p_m\{2\Phi(c_{i,m}^{O}/\sqrt{1+u_{i,m}})-1\}\right].
\end{equation}
Since \((\boldsymbol{X},Z)\) contains more information than \(\boldsymbol{X}\),
\(\mathcal R_m^{\mathrm{BO},Z}\le\mathcal R_m^{\mathrm{BO},\boldsymbol{X}}\).
The exact threshold and risk in this augmented experiment are not by
themselves formulas for the observed-data Oracle.

For a feasible rule in the nondegenerate one-factor case, use leave-one-out GLS:
\begin{equation}
\label{eq:main-gls}
D_{-i}=\sum_{j\ne i}\lambda_{j,m}^2/\sigma_{j,m}^2,\qquad
\widehat Z_{-i}=\frac{\sum_{j\ne i}\lambda_{j,m}X_j/\sigma_{j,m}^2}{D_{-i}},
\qquad
\widehat Y_i=\frac{X_i-\lambda_{i,m}\widehat Z_{-i}}{\sigma_{i,m}}.
\end{equation}
The assumptions ensure \(D_{-i}\asymp m\) uniformly, eventually.
Equivalently, the fixed residual can be written using covariance entries
as
\begin{equation}\label{eq:gls-known-sigma}
\widehat Y_i=\frac1{\sigma_i}\left\{X_i-
\frac{\sum_{j\ne i}(\boldsymbol{\Sigma}_m)_{ij}X_j/\sigma_j^2}
{\sum_{j\ne i}(1-\sigma_j^2)/\sigma_j^2}\right\}.
\end{equation}
Here the $\sigma_j^2$ are determined by $\boldsymbol{\Sigma}_m$ under
Lemma~\ref{lem:known-sigma-loadings}; this expression is independent of
the orientation chosen for the factor.
Writing \(\widehat Y_i=Y_i(Z)-\Delta_{i,m}\), the leave-one-out
construction makes \(\Delta_{i,m}\) independent of \((Y_i(Z),\nu_i)\).
It is conditionally Gaussian given the other latent indicators. This
independence allows the null-tail comparison to be made before averaging
over factor-estimation error. In the homogeneous case,
\(\widehat Z_{-i}=\{\lambda(m-1)\}^{-1}\sum_{j\ne i}X_j\).

Under exact independence, $\lambda_{i,m}=0$ for every $i$ and
$\boldsymbol\Sigma_m=\boldsymbol I_m$, so the factor-estimation formulas
above are not used; define $\widehat Y_i=X_i$. Then fixed residualized
GBS and BH reduce exactly to ordinary GBS and BH, and the MRD residuals
satisfy $U_{i,A}=X_i$ for every active set, so MRD--GBS also reduces to
GBS. Their independent-coordinate ABOS properties follow under the
respective assumptions and calibrations from \citet{BCFG2011} for BH and
\citet{Ghosh2026GBSABOS} for GBS. The dependence theorems below concern
the nondegenerate one-factor case.

Compute \(\widehat P_i=2\bar\Phi(|\widehat Y_i|)\) once and order them.
The fixed-residual GBS rule has rejection count
\begin{equation}
\label{eq:framework-procedure}
R_m=\max\{r\in\{0,\ldots,m\}:\widehat P_{(j)}\le c_j(\alpha_m)
\text{ for every }1\le j\le r\}.
\end{equation}
It rejects the first \(R_m\) hypotheses. We write rGBS for homogeneous
loadings and hrGBS for heterogeneous loadings. This feasible rule uses the
known loading structure and the nominal level, but not the unknown signal
probability or slab variance in its rejection calculation. The nominal
level is a tuning sequence.
Unlike MRD--GBS, this rule does not update residuals after deletions.

\subsection{Residualized Benjamini--Hochberg}
\label{sec:rbh-definition}
Residualized BH (rBH) uses exactly the fixed GLS residuals
$\widehat Y_i$ in \eqref{eq:main-gls} and their two-sided p-values
$\widehat P_i=2\bar\Phi(|\widehat Y_i|)$. Define
\begin{equation}\label{eq:rbh-procedure}
R_m^B=\max\bigl(\{j\in\{1,\ldots,m\}:\widehat P_{(j)}
\le j\alpha_m/m\}\cup\{0\}\bigr).
\end{equation}
Reject the first $R_m^B$ hypotheses. This is a fixed-residual step-up
procedure, not marginal BH and not an MRD procedure with BH thresholds.
The denominator in its critical line is the original dimension $m$.
Write $V_m^B,T_m^B,\mathcal R_m^B$ for its false discoveries, missed
signals and Bayes risk, respectively.

\subsection{The risk questions and calibration constants}

The theory answers two questions in order. Does the optimal observed-data
risk have the same leading term as \eqref{eq:main-factor-known-risk}?
Do the fixed-residual and adaptive rules attain that term? The first
requires a risk sandwich; the second also requires controlling the entire
step-down crossing event. Neither follows from the empirical observation
that a different adaptive procedure closely tracks the Oracle.

For homogeneous loadings write \(a_\lambda^2=1-\lambda^2\),
\(u_{m,\lambda}=\psi_m^2/a_\lambda^2\), and use
\(R_{m,\lambda}^{\mathrm{BO}}\), \(R_{m,\lambda}^{\mathrm{rGBS}}\) for
risks. In the heterogeneous setting use
\(\mathcal R_m^{\mathrm{BO},\boldsymbol{X}}\), \(\mathcal R_m^{\mathrm{hrGBS}}\).
The quantities \(R_m^{\mathrm{rGBS}}\) and \(R_m^{\mathrm{hrGBS}}\)
are rejection counts. Set \(r_{\alpha_m}=\alpha_m/(1-\alpha_m)\) and
\(z_t=\Phi^{-1}(1-t/2)\). For \(C_0>0\), define
\begin{equation}
\label{eq:framework-Qm}
Q_m=\frac1m\sum_i2\bar\Phi(\sqrt{C_0}\sigma_{i,m}).
\end{equation}
The loading bounds give fixed \(0<q_-<q_+<1\) such that
\(q_-\le Q_m\le q_+\), without requiring convergence. Homogeneous
loadings give \(C_\lambda=C_0(1-\lambda^2)\) and
\(1-Q_m=2\Phi(\sqrt{C_\lambda})-1\). The main theorems specify the
signal, level and loss calibrations under which these constants determine
risk. Both divergent and bounded positive expected counts will be treated;
vanishing expected counts are outside the stated results.

\section{Asymptotic Oracle risk and feasible optimality}
\label{sec:major-results}

We first establish the observed-data Oracle benchmark under homogeneous
and heterogeneous one-factor dependence, then prove ABOS for the complete
adaptive MRD--GBS procedure, the principal procedure of this paper. Fixed
residualized GBS and residual BH provide companion results that isolate the
effect of fixed residualization and the testing rule. Where established, we
cover both diverging and bounded positive expected signal counts. Every
risk ratio uses the same information and loss in numerator and denominator,
as required by Definition~\ref{def:ABOS}. The adaptive proof controls false
discoveries by comparison with a factor-known static rule and controls
missed signals by a separate crossing argument; the fixed-rule proofs use
direct static comparisons. The proofs and supporting lemmas appear in the
appendix.

\subsection{Observed-data Oracle under homogeneous dependence}
We first determine the risk of the genuine Bayes Oracle, which observes the full data vector but not the latent common factor. This provides the benchmark for assessing the procedures below and makes explicit when estimating the factor has no first-order cost. The following theorem establishes this benchmark under homogeneous one-factor dependence.

\begin{theorem}[Bayes Oracle risk under homogeneous one-factor dependence]
	\label{thm:one-factor-oracle-risk}
    Consider the homogeneous specialization of
    \eqref{eq:framework-model} with fixed \(0<|\lambda|<1\) and the
    sparse prior \eqref{eq:general-prior}. With the notation of
    Section~\ref{sec:procedures}, suppose
	\[
	p_m\to0,
	\qquad
	mp_m\to\infty,
	\qquad
	u_{m,\lambda}\to\infty,
	\]
	and
	\begin{equation}
		\label{eq:one-factor-calibration}
		\frac{\log\{u_{m,\lambda}f_m^2\delta_m^2\}}
		{u_{m,\lambda}}
		\to
		C_\lambda\in(0,\infty).
	\end{equation}
	Assume further that
	\begin{equation}
		\label{eq:factor-estimation-condition}
		\frac{
			u_{m,\lambda}\{1+p_m u_{m,\lambda}\}
		}{m}
		\to0.
	\end{equation}
	Then the genuine Bayes Oracle risk, based only on
	\(\boldsymbol{X}=(X_1,\ldots,X_m)\), satisfies
	\begin{equation}
		\label{eq:one-factor-oracle-risk-final}
		R_{m,\lambda}^{\mathrm{BO}}
		=
		mp_m\delta_{A,m}
		\left\{
		2\Phi(\sqrt{C_\lambda})-1
		\right\}
		\{1+o(1)\}.
	\end{equation}
\end{theorem}

\paragraph{Implication.}
The leading risk equals that of the experiment in which the common factor
is observed. Thus the latent factor need not create a first-order Oracle
risk penalty when it can be estimated sufficiently accurately. The extra
condition controls estimation on the Oracle tail scale; it is stronger
than merely requiring a consistent factor estimate. Under the common
raw-signal calibration \(\log\{\psi_m^2 f_m^2\delta_m^2\}/\psi_m^2\to C_0\),
the constant is \(C_\lambda=(1-\lambda^2)C_0\). The resulting smaller
miss-probability constant quantifies the benefit of removing shared noise,
subject to the conditions of this theorem.

\paragraph{Comparison with independence.}
Consider the independent Gaussian sequence model with the same prior signal
proportion, slab variance \(\psi_m^2\), and loss ratio, under the same
raw-signal calibration with limit \(C_0\). Its Bayes-Oracle risk has
leading term
\[
mp_m\delta_{A,m}\{2\Phi(\sqrt{C_0})-1\}.
\]
Consequently, the ratio of the dependent and independent observed-data
Oracle risks converges to
\[
\frac{2\Phi\!\left(\sqrt{(1-\lambda^2)C_0}\right)-1}
     {2\Phi(\sqrt{C_0})-1}<1,
\]
since \(0<|\lambda|<1\). Thus, under this matched calibration, the
Oracle has strictly smaller leading risk under one-factor dependence than
under independence. The common component can be estimated accurately enough
that exploiting the known covariance removes shared noise without a
first-order estimation penalty; the remaining idiosyncratic variance is
\(1-\lambda^2\). In this qualified decision-theoretic sense, dependence
is a blessing rather than a curse. This conclusion concerns the Oracle
that uses the covariance in the stated one-factor regime. It is not a claim
that arbitrary dependence helps, or that a procedure which ignores
dependence enjoys the same risk reduction.

\subsection{Observed-data Oracle with heterogeneous loadings}
We next allow the factor loadings to vary across coordinates. The resulting Oracle risk reflects the different idiosyncratic variances through the finite-sample average \(Q_m\). The theorem does not require this average to converge or the loadings to follow a limiting distribution.
\begin{theorem}[Genuine Bayes-Oracle risk under known heterogeneous one-factor dependence]
\label{thm:hetero-observed-oracle}
Assume the loading conditions \eqref{eq:framework-loadings},
\[
p_m\to0,\quad mp_m\to\infty,\quad \alpha_m\to0,\quad
\alpha_m f_m\to\infty,
\]
and
\[
\frac{2\log(f_m/r_{\alpha_m})}{\psi_m^2}\to C_0\in(0,\infty),
\qquad
\frac{\log\{1/(r_{\alpha_m}\delta_m)\}}{\psi_m^2}\to0.
\]
Then the genuine Bayes Oracle based only on
\(\boldsymbol{X}=(X_1,\ldots,X_m)\) satisfies
\begin{equation}
\label{eq:hetero-observed-oracle-risk}
\mathcal R_m^{\mathrm{BO},\boldsymbol{X}}
=
mp_m\delta_{A,m}\{1-Q_m+o(1)\}.
\end{equation}
\end{theorem}

\paragraph{Implication.}
Different coordinates contribute different first-order miss probabilities.
The expansion keeps the finite-dimensional average \(Q_m\), rather than
assuming it converges. The loading signs affect factor estimation, whereas
the leading coordinatewise risk depends on their magnitudes through
\(\sigma_{i,m}^2=1-\lambda_{i,m}^2\). The model permits a mixture of
strongly and weakly loaded coordinates under the uniform loading bounds.
The proof's leave-one-out Oracle-threshold rule is an $\boldsymbol{X}$-measurable
comparison rule, not a proposed data-driven procedure.

\paragraph{Comparison with independence under heterogeneous loadings.}
Under the assumptions of Theorem~\ref{thm:hetero-observed-oracle}, and
with the independent-coordinate benchmark matched in prior, loss, and
raw-signal calibration, that Oracle has leading missed-signal coefficient
$1-q_0$, where
$q_0=2\bar\Phi(\sqrt{C_0})$. The heterogeneous coefficient is $1-Q_m$.
For $t\in(0,\bar\lambda^2]$, define
\[
 d(t)=2\bar\Phi\bigl(\sqrt{C_0(1-t)}\bigr)-2\bar\Phi(\sqrt{C_0}),
 \qquad d(0)=0.
\]
The function $d(t)/t$ extends continuously to zero and is strictly
positive on $[0,\bar\lambda^2]$. Hence, with
$c_{C_0,\bar\lambda}=\min_{0\le t\le\bar\lambda^2}d(t)/t>0$,
the loading conditions imply
\[
 Q_m-q_0=\frac1m\sum_{i=1}^m d(\lambda_{i,m}^2)
 \ge c_{C_0,\bar\lambda}\,\frac1m\sum_{i=1}^m\lambda_{i,m}^2
 \ge c_{C_0,\bar\lambda}\,\kappa+o(1).
\]
Thus, under this matched calibration, heterogeneous one-factor dependence
also lowers the leading missed-signal coefficient by a nonvanishing
amount, without requiring $Q_m$ to converge. The comparison concerns an
Oracle that uses the known covariance; it is not a claim about arbitrary
dependence or unadjusted testing.

\subsection{Adaptive MRD--GBS: optimality of the complete algorithm}
\label{sec:adaptive-mrd-results}
Assume $\boldsymbol{\Sigma}_m$ is known and has the normalized one-factor form
\eqref{eq:onefactor-covariance}, with a representation satisfying
\eqref{eq:framework-loadings}. By
Lemma~\ref{lem:known-sigma-loadings}, this is equivalent, for the
procedures, to knowing the loadings up to a common sign.
For the actual algorithm \eqref{eq:mrd-gbs-description}, write
$R_m^M,V_m^M,T_m^M$ for rejections, false discoveries and missed signals,
and $\mathcal R_m^M=\delta_{A,m}\{\delta_m\mathbb EV_m^M+\mathbb ET_m^M\}$.
Use a fixed measurable tie-breaking convention and stop when all hypotheses
have been rejected. All losses are positive and $p_m,\alpha_m\in(0,1)$.
Throughout, $\sigma_{i,m}=\sqrt{1-\lambda_{i,m}^2}>0$.
Information ordering gives
\begin{equation}\label{mrd:eq:information-order}
\mathcal R_m^{\mathrm{BO},Z}\le\mathcal R_m^{\mathrm{BO},\boldsymbol{X}}\le\mathcal R_m^M.
\end{equation}
The following theorems concern recomputation after every rejection, with
no assumed short-path bound and no equality of adaptive and fixed
rejection paths. Appendix~\ref{app:adaptive-mrd} gives the proofs.

\begin{theorem}[MRD--GBS under general vanishing sparsity: diverging expected counts]\label{mrd:thm:diverging}
Assume the known-covariance one-factor model just specified and
\begin{gather}
 p_m=o(1),\qquad mp_m\to\infty,\qquad
 \alpha_m\to0,\qquad\alpha_mf_m\to\infty,\label{mrd:eq:div-counts}\\
 \frac{2\log(f_m/r_{\alpha_m})}{\psi_m^2}\to C_0\in(0,\infty),
 \qquad
 \frac{\log\{1/(r_{\alpha_m}\delta_m)\}}{\psi_m^2}\to0,
 \qquad\alpha_m\delta_m\to0.\label{mrd:eq:div-calibration}
\end{gather}
Use $Q_m$ and its fixed bounds $q_-,q_+$ from \eqref{eq:framework-Qm}.
Then
\begin{align}
\mathbb EV_m^M&=o(mp_m),\qquad
\delta_m\mathbb EV_m^M=o(mp_m),\label{mrd:eq:main-typeI}\\
\mathbb ET_m^M&=mp_m\{1-Q_m+o(1)\},\label{mrd:eq:main-typeII}\\
\mathcal R_m^M&=mp_m\delta_{A,m}\{1-Q_m+o(1)\},\qquad
\frac{\mathcal R_m^M}{\mathcal R_m^{\mathrm{BO},\boldsymbol{X}}}\to1.\label{mrd:eq:main-abos}
\end{align}
The sparsity assumption is $p_m=o(1)$, with no polynomial or
additional logarithmic sparsity-rate restriction. The other displayed
count, signal, level and loss assumptions remain in force. Neither a
limiting empirical loading distribution nor convergence of $Q_m$ is required.
\end{theorem}
The preceding theorem treats the regime in which the expected number of
signals diverges. We now consider the bounded-count regime, where that
expectation stays bounded above and away from zero. This covers a distinct
sparsity endpoint and shows whether the adaptive MRD--GBS rule retains its
Oracle-risk performance there.
\begin{theorem}[Adaptive MRD--GBS: bounded positive expected counts]\label{mrd:thm:bounded}
Assume \eqref{eq:framework-model} under \eqref{eq:framework-loadings}, and for constants
$0<b_-\le b_+<\infty$ suppose
\[
b_-\le mp_m\le b_+\quad\text{eventually},\qquad
\frac{\psi_m^2}{\log m}\to\frac2{C_0},\quad C_0\in(0,\infty).
\]
Suppose $\delta_m=\delta\in(0,\infty)$,
$\alpha_m\to0$, and $\log(1/\alpha_m)=o(\log m)$.
Then \eqref{mrd:eq:main-typeI}--\eqref{mrd:eq:main-abos} hold. The expectation
$mp_m$ need not converge. No additional sparsity rate is imposed.
\end{theorem}

\begin{corollary}[General vanishing sparsity under known one-factor covariance]
\label{cor:general-vanishing-known-sigma}
Assume the known covariance $\boldsymbol{\Sigma}_m$ has form
\eqref{eq:onefactor-covariance} under \eqref{eq:framework-loadings}.
Suppose
\begin{gather}
p_m=o(1),\qquad \liminf_{m\to\infty}mp_m>0,
\qquad\delta_m=\delta\in(0,\infty),\label{eq:general-sparsity-counts}\\
\frac{\psi_m^2}{2\log(1/p_m)}\to\frac1{C_0},\qquad C_0>0,
\qquad \alpha_m\to0,\qquad
\frac{\log(1/\alpha_m)}{\log(1/p_m)}\to0.
\label{eq:general-sparsity-calibration}
\end{gather}
Then MRD--GBS, fixed GLS residualized GBS, and residual BH are all ABOS
relative to the genuine observed-data Oracle. For each of their risks,
\[
\mathcal R_m=mp_m\delta_{A,m}\{1-Q_m+o(1)\},\qquad
\frac{\mathcal R_m}{\mathcal R_m^{\mathrm{BO},\boldsymbol{X}}}\to1.
\]
The same leading expansion holds for $\mathcal R_m^{\mathrm{BO},\boldsymbol{X}}$.
The expectation $mp_m$ may be bounded, divergent, or have both bounded
and divergent subsequences; no limit for it or for $Q_m$ is required.
\end{corollary}

\paragraph{Implication.}
This is a single general-sparsity statement with a fixed positive loss
ratio. It makes explicit that $p_m=o(1)$ does not need to be strengthened
by an extra rate involving $\log m$. For example both
$p_m=1/\log\log m$ and every polynomial sequence
$p_m\asymp m^{-\beta}$ with $0<\beta\le1$ are included when
\eqref{eq:general-sparsity-calibration} holds. The lower bound on expected
counts and the signal/level calibration are substantive assumptions;
$p_m=o(1)$ alone is not the whole theorem. A proof combining the two count
regimes appears in Appendix~\ref{app:covariance-sparsity}.

\begin{corollary}[All polynomial sparsity exponents]\label{mrd:cor:polynomial}
In model \eqref{eq:framework-model} under \eqref{eq:framework-loadings}, suppose, for fixed
$\beta\in(0,1]$,
\[
p_m\asymp m^{-\beta},\qquad
\frac{\psi_m^2}{\log m}\to\frac{2\beta}{C_0},\quad C_0>0,
\qquad\delta_m=\delta\in(0,\infty),
\]
and $\alpha_m\to0$, $\log(1/\alpha_m)=o(\log m)$. Then actual adaptive
MRD--GBS is ABOS relative to the observed-data Oracle, with leading risk
$mp_m\delta_{A,m}\{1-Q_m+o(1)\}$.
\end{corollary}

\begin{corollary}[Fixed equicorrelation]\label{mrd:cor:homogeneous}
Set $\lambda_{i,m}\equiv\lambda$, $0<|\lambda|<1$, and
$u_{m,\lambda}=\psi_m^2/(1-\lambda^2)$. The preceding results specialize
to fixed equicorrelation $\rho=\lambda^2$. With
$C_\lambda=(1-\lambda^2)C_0$, the leading risk becomes
\[
mp_m\delta_{A,m}\{2\Phi(\sqrt{C_\lambda})-1\}\{1+o(1)\}.
\]
For polynomial sparsity the equivalent signal scaling is
$u_{m,\lambda}/\log m\to2\beta/C_\lambda$.
\end{corollary}

\paragraph{Implication.}
These results establish ABOS for the complete adaptive procedure without
an additional sparsity-rate condition beyond the displayed count and
calibration assumptions. At stage $k$, only $k-1$ deletions are possible,
and their effect is scaled by $z_{c_k}$. Matching these two stage-dependent
quantities avoids a stronger restriction from applying the first-stage
threshold to every set within the cap. This control yields
the automatic limit \eqref{mrd:eq:automatic-rate}.
For example, $p_m=1/\log\log m$ is covered with fixed positive loss
ratio, $\alpha_m=\exp[-\sqrt{\log(1/p_m)}]$, and
$\psi_m^2\sim2\log(1/p_m)/C_0$.
More generally these choices satisfy the divergent-count calibrations
for every $p_m\to0$ with $mp_m\to\infty$. The displayed model, loading,
and level/loss assumptions remain part of the result.

\subsection{Oracle-parameterized BSD under one-factor dependence}
\label{sec:bsd-results}

BSD is included here as a model-informed Bayesian step-down rule. Its
stagewise odds score uses the true $p_m$, $\psi_m^2$, known covariance, and
loss ratio $\delta_m$; these inputs are not estimated from the observed
vector. Thus BSD is $\boldsymbol{X}$-driven in its sequential choices but is not a fully
data-driven or parameter-adaptive procedure. It is also distinct from the
joint Bayes Oracle: BSD uses local sequential model comparisons, whereas the
Oracle minimizes the full posterior expected loss over all configurations.
This distinction is essential when interpreting the ABOS result below.
For consistency with the general loss notation, put
$\Delta_m:=\delta_m=\delta_{0,m}/\delta_{A,m}$.
Under the one-factor representation, write
\[
\eta_{i,m}=\lambda_{i,m}Z+\sigma_{i,m}\varepsilon_i,
\qquad \boldsymbol{\eta}_{A,m}=(\eta_{j,m}:j\in A),
\]
so that $X_i=\theta_i+\eta_{i,m}$. For an active set $A\subseteq[m]$ and $i\in A$, let
\[
v_{i,A}=\operatorname{Var}(\eta_{i,m}\mid\boldsymbol{\eta}_{A\setminus\{i\},m}),
\qquad
U_{i,A}=\frac{X_i-\mathbb E_0(X_i\mid \boldsymbol{X}_{A\setminus\{i\}})}{\sqrt{v_{i,A}}},
\]
where the conditional variance and regression are under the known Gaussian error model.
A Gaussian conditioning calculation gives
\begin{equation}
\label{bsd:eq:conditional-variance}
v_{i,A}=\sigma_{i,m}^2+
\frac{\lambda_{i,m}^2}
{1+\sum_{j\in A\setminus\{i\}}\lambda_{j,m}^2/\sigma_{j,m}^2}.
\end{equation}
BSD uses the stagewise posterior-odds score
\begin{equation}
\label{bsd:eq:score}
S_{i,A}=\frac{p_m}{1-p_m}
\sqrt{\frac{v_{i,A}}{\psi_m^2+v_{i,A}}}
\exp\!\left\{\frac{\psi_m^2}{2(\psi_m^2+v_{i,A})}U_{i,A}^2\right\}.
\end{equation}
At each stage it selects a maximizer of $S_{i,A}$, rejects it when
$S_{i,A}>\Delta_m$, removes that coordinate, recomputes the scores, and
stops when the largest remaining score is at most $\Delta_m$. The score
uses the true $p_m$, $\psi_m^2$, covariance, and loss ratio.

\begin{theorem}[Oracle-parameterized BSD under known one-factor covariance]
\label{bsd:thm:abos}
Suppose
\[
p_m\to0,\qquad \liminf_{m\to\infty}mp_m>0,\qquad
\frac{\psi_m^2}{2\log(1/p_m)}\longrightarrow \frac1{C_0}
\quad\text{for some }C_0\in(0,\infty),
\]
and \(\Delta_m=\delta_{0,m}/\delta_{A,m}\to\Delta\in(0,\infty)\).
Assume the known one-factor loadings satisfy
\(\sup_m\max_i|\lambda_{i,m}|\le\bar\lambda<1\) and
\(\liminf_m m^{-1}\sum_i\lambda_{i,m}^2\ge\kappa>0\).
Then the BSD rule defined by \eqref{bsd:eq:score} satisfies
\[
\mathbb E V_m^{\mathrm{BSD}}=o(mp_m),\qquad
\delta_m\mathbb E V_m^{\mathrm{BSD}}=o(mp_m),\qquad
\mathbb E T_m^{\mathrm{BSD}}=mp_m\{1-Q_m+o(1)\},
\]
and hence
\[
\mathcal R_m^{\mathrm{BSD}}
=mp_m \delta_{A,m}\{1-Q_m+o(1)\}.
\]
Moreover, the observed-data Bayes-Oracle expansion under the same model is
\[
\mathcal R_m^{\mathrm{BO},\boldsymbol{X}}
=mp_m \delta_{A,m}\{1-Q_m+o(1)\},
\]
and therefore BSD is ABOS relative to that Oracle. No convergence of \(Q_m\) or limiting loading distribution is required.
\end{theorem}

\begin{corollary}[Normalized agreement of the BSD and Oracle decisions]
\label{bsd:cor:decision-agreement}
Under the assumptions of Theorem~\ref{bsd:thm:abos}, suppose in addition that
\[
\frac{mp_m}{\log(1/p_m)}\longrightarrow\infty.
\]
Let $\widehat{\mathcal S}_m^{\mathrm{BSD}}$ and
$\widehat{\mathcal S}_m^{\mathrm{BO},\boldsymbol X}$ denote the terminal
rejection sets of BSD and the observed-data Bayes Oracle, respectively. Then
\[
\frac{\bigl|\widehat{\mathcal S}_m^{\mathrm{BSD}}
\triangle\widehat{\mathcal S}_m^{\mathrm{BO},\boldsymbol X}\bigr|}{mp_m}
\xrightarrow{P}0.
\]
In particular, this conclusion holds for polynomial sparsity
$p_m\asymp m^{-\beta}$ with any fixed $0<\beta<1$.
\end{corollary}

\paragraph{Interpretation.}
The corollary gives asymptotic agreement of the terminal decision sets on
the scale of the expected number of signals. It therefore provides a
theoretical, decision-level explanation for the near-Oracle behavior of BSD
reported in the one-factor simulations of the companion BSD manuscript
\citep{GhoshChakrabartiBSD2026}. The result does not assert equality of the
two sets with probability tending to one, nor agreement of their intermediate
decision paths.

\paragraph{Implication.}
Under fixed positive loss ratio and $\liminf mp_m>0$, this theorem adds an
ABOS result for the Oracle-parameterized BSD rule over general vanishing
sparsity, including all polynomial exponents $0<\beta\le1$ when the signal
variance is calibrated as stated. The proofs of the theorem and
Corollary~\ref{bsd:cor:decision-agreement} are in Appendix~\ref{app:bsd-abos}.

\subsection{Fixed residualized GBS under homogeneous dependence}
The adaptive MRD--GBS procedure recomputes its residuals along the
rejection path. We now consider fixed residualized GBS, which applies the
GBS step-down rule to the initial GLS residuals. The following result
establishes its first-order risk under homogeneous one-factor dependence.
\begin{theorem}[ABOS of residualized GBS under homogeneous one-factor dependence]
\label{thm:rgbs-abos}
Consider model \eqref{eq:framework-model} with fixed
$0<|\lambda|<1$.  Suppose
\[
p_m\to0,
\qquad
mp_m\to\infty,
\qquad
u_{m,\lambda}\to\infty,
\]
\[
\alpha_m\to0,
\qquad
\alpha_m f_m\to\infty,
\]
and
\begin{equation}
\label{eq:rgbs-abos-calibration}
\frac{2\log(f_m/r_{\alpha_m})}{u_{m,\lambda}}
\to C_\lambda\in(0,\infty),
\qquad
r_{\alpha_m}=\frac{\alpha_m}{1-\alpha_m}.
\end{equation}
Assume also the compatibility and Type-I loss conditions
\begin{equation}
\label{eq:rgbs-abos-compatibility}
\frac{\log\{1/(r_{\alpha_m}\delta_m)\}}{u_{m,\lambda}}
\to0,
\qquad
\alpha_m\delta_m\to0.
\end{equation}
Then the residualized GBS procedure defined in
\eqref{eq:framework-procedure} is asymptotically Bayes
optimal under sparsity:
\begin{equation}
\label{eq:rgbs-abos-ratio}
\frac{R_{m,\lambda}^{\mathrm{rGBS}}}
{R_{m,\lambda}^{\mathrm{BO}}}
\to1.
\end{equation}
More precisely,
\begin{equation}
\label{eq:rgbs-risk-expansion}
R_{m,\lambda}^{\mathrm{rGBS}}
=
mp_m\delta_{A,m}
\{2\Phi(\sqrt{C_\lambda})-1\}
\{1+o(1)\}.
\end{equation}
\end{theorem}

\paragraph{Implication.}
A computationally feasible procedure attains the first-order minimum risk
without observing the factor or using the sparsity and slab variance in
its rejection rule. The loss-weighted Type-I contribution is negligible,
while the missed-signal contribution attains the Oracle constant.

\subsection{Fixed residualized GBS under heterogeneous loadings}
We extend the fixed residualized GBS result to known heterogeneous
loadings. Each coordinate is adjusted using its own loading, and the
leading risk is determined by the resulting coordinate-specific
contributions. No limiting loading distribution is imposed.
\begin{theorem}[ABOS of residualized GBS under known heterogeneous one-factor dependence]
\label{thm:hrgbs-abos}
Consider \eqref{eq:framework-model}.  Suppose
\eqref{eq:framework-loadings}
holds and
\[
p_m\to0,\qquad mp_m\to\infty,
\qquad
\alpha_m\to0,\qquad \alpha_mf_m\to\infty.
\]
Assume
\[
\frac{2\log(f_m/r_{\alpha_m})}{\psi_m^2}
\to C_0\in(0,\infty),
\]
and
\[
\frac{\log\{1/(r_{\alpha_m}\delta_m)\}}{\psi_m^2}\to0,
\qquad
\alpha_m\delta_m\to0.
\]
Then the heterogeneous residualized GBS procedure
\eqref{eq:framework-procedure} is asymptotically Bayes optimal under
sparsity relative to the genuine \(\boldsymbol{X}\)-measurable Bayes Oracle:
\begin{equation}
\label{eq:hrgbs-abos-ratio}
\frac{\mathcal R_m^{\mathrm{hrGBS}}}{\mathcal R_m^{\mathrm{BO},\boldsymbol{X}}}\to1.
\end{equation}
More precisely,
\begin{equation}
\label{eq:hrgbs-risk-expansion}
\mathcal R_m^{\mathrm{hrGBS}}
=
p_m\delta_{A,m}
\sum_{i=1}^m
\left[
2\Phi\{\sqrt{C_0(1-\lambda_{i,m}^2)}\}-1
\right]
+
o(mp_m\delta_{A,m}).
\end{equation}
\end{theorem}

\paragraph{Implication.}
The same step-down constants remain sufficient after coordinate-specific
GLS residualization. Unequal loadings change the leading risk through
\(Q_m\), but do not require a limiting loading distribution. Known
one-factor covariance supplies the loading representation through
Lemma~\ref{lem:known-sigma-loadings}. In particular, heterogeneous loadings do
not reduce this result to an exchangeable covariance model.

\subsection{Optimality with bounded positive expected counts}
The preceding result concerns fixed residualized GBS when the expected
number of signals diverges. We now consider bounded positive expected
counts, where concentration around a diverging signal count is unavailable.
The theorem gives matching risk expansions for the procedure and the
observed-data Oracle in this regime.
\begin{theorem}[Heterogeneous rGBS with bounded expected signal counts]
\label{thm:hrgbs-bounded}
Consider model \eqref{eq:framework-model} under
\eqref{eq:framework-loadings}.
Suppose that, for fixed constants \(0<b_-\le b_+<\infty\),
\[
b_-\le mp_m\le b_+\quad\text{eventually},\qquad
\frac{\psi_m^2}{\log m}\to\frac{2}{C_0},\quad C_0\in(0,\infty),
\]
and that
\[
\delta_m=\delta\in(0,\infty),\qquad
\alpha_m\to0,\qquad \log(1/\alpha_m)=o(\log m).
\]
Then, with \(Q_m\) defined by \eqref{eq:framework-Qm},
\begin{align}
\mathcal R_m^{\mathrm{BO},\boldsymbol{X}}
&=mp_m\delta_{A,m}\{1-Q_m+o(1)\},
\label{eq:bounded-hetero-oracle}\\
\mathcal R_m^{\mathrm{hrGBS}}
&=mp_m\delta_{A,m}\{1-Q_m+o(1)\}.
\label{eq:bounded-hetero-rgbs}
\end{align}
In particular, \(\mathcal R_m^{\mathrm{hrGBS}}/
\mathcal R_m^{\mathrm{BO},\boldsymbol{X}}\to1\). No convergence of \(mp_m\) or of
\(Q_m\) is required.
\end{theorem}

\paragraph{Implication.}
ABOS persists when only a bounded positive number of signals is expected.
The count may fluctuate substantially, and its expectation need not
converge. The proof does not use concentration around a diverging mean:
a signal crossing the first critical value is necessarily rejected.
The lower bound on expected counts is material; the result does not extend
by its statement to \(mp_m\to0\). Here the loss ratio is fixed.

\subsection{Homogeneous bounded-count specialization}
\begin{corollary}[Homogeneous rGBS with bounded expected signal counts]
\label{cor:rgbs-bounded}
In the homogeneous model with fixed \(0<|\lambda|<1\), suppose
\[
0<b_-\le mp_m\le b_+<\infty\quad\text{eventually},\qquad
\frac{u_{m,\lambda}}{\log m}\to\frac{2}{C_\lambda},
\quad C_\lambda\in(0,\infty).
\]
If \(\delta_m=\delta\in(0,\infty)\), \(\alpha_m\to0\), and
\(\log(1/\alpha_m)=o(\log m)\), then
\[
R_{m,\lambda}^{\mathrm{BO}}
=mp_m\delta_{A,m}\{2\Phi(\sqrt{C_\lambda})-1\}\{1+o(1)\},
\qquad
\frac{R_{m,\lambda}^{\mathrm{rGBS}}}{R_{m,\lambda}^{\mathrm{BO}}}\to1.
\]
\end{corollary}

\paragraph{Implication.}
This specializes the bounded-count theorem to fixed equicorrelation and
recovers the same homogeneous Oracle constant as the diverging-count
analysis. The two regimes require different proofs but have the same
first-order efficiency interpretation.

\subsection{Polynomial sparsity, including its bounded-count endpoint}
\begin{corollary}[ABOS throughout polynomial sparsity under one-factor dependence]
\label{cor:onefactor-polynomial}
Fix \(\beta\in(0,1]\), and suppose
\[
p_m\asymp m^{-\beta},\qquad \delta_m=\delta\in(0,\infty),\qquad
\alpha_m\to0,\qquad \log(1/\alpha_m)=o(\log m).
\]
\begin{enumerate}[label=(\roman*)]
\item In the heterogeneous model under
\eqref{eq:framework-loadings},
if \(\psi_m^2/\log m\to2\beta/C_0\), with \(C_0\in(0,\infty)\),
then hrGBS is ABOS and
\[
\mathcal R_m^{\mathrm{hrGBS}}
=mp_m\delta_{A,m}\{1-Q_m+o(1)\}.
\]
\item In the homogeneous model with fixed \(0<|\lambda|<1\), if
\(u_{m,\lambda}/\log m\to2\beta/C_\lambda\), with
\(C_\lambda\in(0,\infty)\), then rGBS is ABOS and
\[
R_{m,\lambda}^{\mathrm{rGBS}}
=mp_m\delta_{A,m}\{2\Phi(\sqrt{C_\lambda})-1\}\{1+o(1)\}.
\]
\end{enumerate}
\end{corollary}

\paragraph{Implication.}
Every fixed exponent \(0<\beta\le1\) is covered under the indicated
signal scaling. For \(\beta<1\), the expected count diverges; at
\(\beta=1\), it remains bounded above and away from zero. A level such as
\(\alpha_m=1/\log m\), for sufficiently large \(m\), satisfies the level
conditions. Signal variance is of order \(\log m\); this statement is not
an assertion of optimality at arbitrary signal strengths. The homogeneous
and heterogeneous calibrations agree after rescaling by the
idiosyncratic variance.

\subsection{Residual BH: ABOS without an adaptive sparsity rate}
We finally consider the BH step-up rule applied to the fixed GLS
residuals. This gives a simpler alternative to residualized GBS while
using the same covariance adjustment. The theorem establishes its
first-order risk in both the diverging-count and bounded positive
count regimes.
\begin{theorem}[Residual BH with divergent or bounded positive expected counts]
\label{thm:rbh-abos}
Consider \eqref{eq:framework-model} under \eqref{eq:framework-loadings},
and apply the fixed-residual BH rule \eqref{eq:rbh-procedure}.
Either of the following sets of assumptions is sufficient:
\begin{enumerate}[label=(\roman*)]
\item All the displayed assumptions of Theorem~\ref{thm:hrgbs-abos}; or
\item All the assumptions of Theorem~\ref{thm:hrgbs-bounded}.
\end{enumerate}
In case (i) there is no additional rate condition on $p_m$.
In both cases,
\begin{align}
\mathbb EV_m^B&=o(mp_m),\qquad
\delta_m\mathbb EV_m^B=o(mp_m),\label{eq:rbh-typeI}\\
\mathbb ET_m^B&=mp_m\{1-Q_m+o(1)\},\label{eq:rbh-typeII}\\
\mathcal R_m^B&=mp_m\delta_{A,m}\{1-Q_m+o(1)\},\qquad
\frac{\mathcal R_m^B}{\mathcal R_m^{\mathrm{BO},\boldsymbol{X}}}\to1.
\label{eq:rbh-abos}
\end{align}
In particular, under the conditions of
Corollary~\ref{cor:onefactor-polynomial}, residual BH is ABOS for every
fixed $0<\beta\le1$. No limiting loading distribution is required.
\end{theorem}

\paragraph{Implication.}
The same leading Bayes risk can be attained using the BH step-up line after
fixed GLS adjustment. This supplies a theoretical explanation for the
competitive residual-BH risks in the numerical study. It concerns the
explicit residuals in \eqref{eq:main-gls}, not ordinary marginal BH and
not BH applied to residuals recomputed after deletions. Together with the
adaptive MRD--GBS theorem, this identifies a common
first-order risk despite different residual updates and rejection mechanisms.

\subsection{Connection with maximum residual down statistics}
At each step, MRD uses the covariance restricted to the current active
set. Under homogeneous one-factor dependence, its standardized statistic
can also be written as an active-set leave-one-out residual. The
following corollary establishes this approximation uniformly over large
active sets; by itself, the approximation does not imply that MRD and a
fixed-residual procedure have the same rejection path.

Write \(\rho=\lambda^2\), \(\boldsymbol{\Sigma}_A=(1-\rho)\boldsymbol{I}_{s}+\rho \boldsymbol{J}_{s}\),
\(\boldsymbol{\Omega}_A=\boldsymbol{\Sigma}_A^{-1}\), and \(s=|A|\ge2\). Define
\[
U_{i,A}^{\mathrm{MRD}}=\frac{(\boldsymbol{\Omega}_A \boldsymbol{X}_A)_i}{\sqrt{(\boldsymbol{\Omega}_A)_{ii}}},
\qquad
L_{i,A}=\frac{X_i-(s-1)^{-1}\sum_{j\in A\setminus\{i\}}X_j}
{\sqrt{1-\rho}}.
\]
These are respectively the standardized MRD statistic and the active-set
leave-one-out residual; \(L_{i,[m]}=\widehat Y_i\).
\begin{corollary}[Uniform equivalence over all large active sets]
\label{cor:mrd-loo-uniform}
Under the assumptions of Theorem~\ref{thm:rgbs-abos},
\[
z_{c_1(\alpha_m)}
\sup_{A:\,|A|\ge \lfloor m/2\rfloor}
\max_{i\in A}
|U_{i,A}^{\mathrm{MRD}}-L_{i,A}|
\xrightarrow{\mathbb P}0.
\]
In particular, at the initial active set \(A=\{1,\ldots,m\}\),
\[
z_{c_1(\alpha_m)}
\max_{1\le i\le m}
|U_{i,[m]}^{\mathrm{MRD}}-\widehat Y_i|
\xrightarrow{\mathbb P}0.
\]
\end{corollary}

\paragraph{Implication.}
The approximation holds uniformly over all large active sets, including
sets selected from the data. It compares MRD to an active-set leave-one-out
residual. Fixed rGBS instead retains its initial full-sample residuals.
Consequently this correspondence alone does not identify rejection sets.
Theorems~\ref{mrd:thm:diverging}--\ref{mrd:thm:bounded} establish adaptive
ABOS by the separate capped-path argument.

\section{Simulation study under one-factor dependence}
\label{sec:simulation}

\subsection{Objectives and design}

The asymptotic results are complemented by a finite-sample comparison of
absolute Bayes risks. This experiment fixes the dimension and varies the
signal proportion; it is designed to show how the procedures behave across
sparsity levels and loading patterns, not to demonstrate convergence as
$m\to\infty$ or to establish ABOS numerically.

We set $m=200$, use symmetric unit losses
$\delta_{0,m}=\delta_{A,m}=1$, take $\alpha=0.10$, and fix the slab variance
at $\psi^2=2\log(m)$. The signal proportion ranges over
\[
 p\in\{0.01,0.02,0.03,0.05,0.07,0.10,0.15,0.20,0.25,0.30,0.35\}.
\]
For each value of $p$ and each loading pattern, we generate $G=5000$
independent data sets from the two-groups one-factor model in
equation~\eqref{eq:framework-model}. The loading patterns, matching those used in the
companion BSD study, are:
\begin{enumerate}[leftmargin=2em,itemsep=1pt]
  \item independence, $\lambda_i=0$;
  \item positive factor loadings, $\lambda_i=\sqrt{0.70}$;
  \item alternating signs, $\lambda_i=(-1)^{i+1}\sqrt{0.70}$;
  \item linear loadings from $0.45$ to $0.85$;
  \item sparse loadings, with $\lambda_i=0.85$ for the first quarter of the
        coordinates and $0.20$ otherwise; and
  \item two-block loadings, with $\lambda_i=0.85$ for the first half and
        $0.35$ otherwise.
\end{enumerate}
In every case, $\sigma_i^2=1-\lambda_i^2$, so the marginal variance of
$X_i$ under the null is one.

We compare the observed-data Bayes Oracle, BSD, fixed residualized GBS
(rGBS), adaptive MRD--GBS, residual BH, and marginal BH. The observed-data
Oracle rejects coordinate $i$ when
$q_i(\boldsymbol{X})=\Pr(\nu_i=1\mid \boldsymbol{X})>1/2$. BSD is included as an
Oracle-parameterized benchmark: its sequential posterior-odds rule uses the
true $p$, $\psi^2$, covariance matrix, and loss ratio. It is therefore not
a fully data-driven competitor. The other procedures use the same nominal
level $\alpha=0.10$. Marginal BH uses the raw two-sided marginal Gaussian
$p$-values, while residual BH applies BH to the fixed residualized
$p$-values used by rGBS. MRD--GBS updates its residual system after each
rejection.

\subsection{Oracle computation and risk estimation}

For each observation vector, the posterior inclusion probabilities are
computed by one-dimensional quadrature over the latent factor. With
$\varphi_i(\boldsymbol{X})\in\{0,1\}$ denoting a rule's decision, we estimate its Bayes
risk per hypothesis by the Rao--Blackwellized estimator
\[
\widehat{\mathcal R}_m(\boldsymbol{\varphi})
=\frac{1}{Gm}\sum_{g=1}^G\sum_{i=1}^m
\left[(1-q_i(\boldsymbol{X}^{(g)}))\varphi_i(\boldsymbol{X}^{(g)})
+q_i(\boldsymbol{X}^{(g)})(1-\varphi_i(\boldsymbol{X}^{(g)}))\right].
\]
For the observed-data Oracle this uses
$\varphi_i^{\mathrm{O}}(\boldsymbol{X})=\mathbf 1\{q_i(\boldsymbol{X})>1/2\}$.
The figures show these absolute risks, rather than ratios to the Oracle;
shaded bands are pointwise $95\%$ Monte Carlo intervals based on the
replicate-wise conditional losses. All procedures are evaluated on the same
simulated data sets.

The posterior calculation uses 64-node Gauss--Hermite quadrature centered
at a numerically determined posterior mode. Across the simulated settings,
64- and 128-node calculations on the first 20 replicates differed by at
most $5.5\times10^{-14}$. On a selected strong-correlation example,
independent adaptive quadrature differed by at most $8.5\times10^{-15}$;
in a seven-coordinate check, posterior probabilities agreed with exact
configuration enumeration within $1.7\times10^{-16}$. Direct checks of the
one-factor residual formulas against covariance-matrix calculations gave
maximum absolute discrepancies below $10^{-15}$. These checks assess
numerical implementation and are not uniform error bounds over all data.

\subsection{Results}

\begin{figure}[t]
	\centering
	\includegraphics[width=\textwidth]{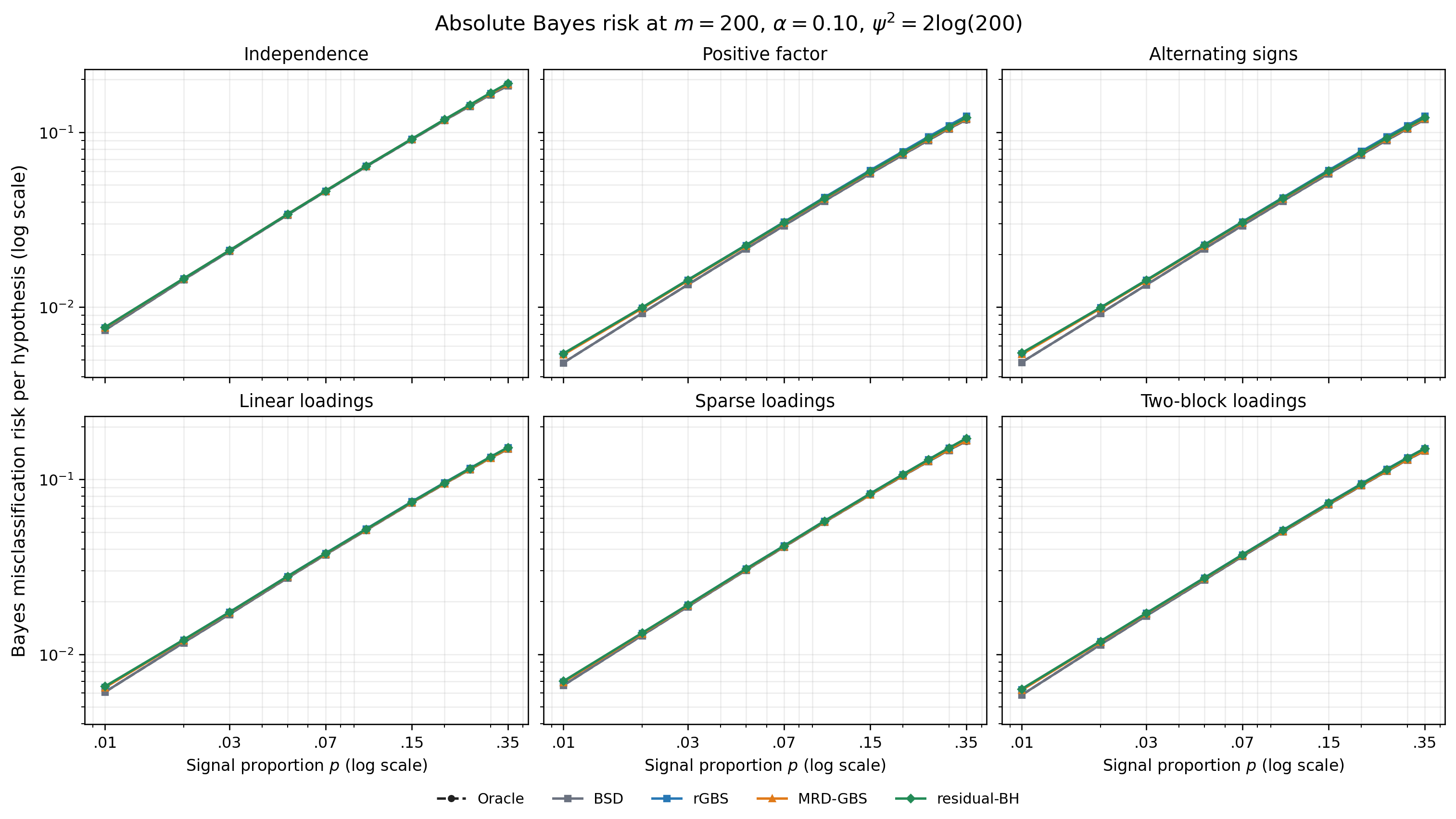}
	\caption{Absolute Bayes misclassification risk per hypothesis for the
		observed-data Oracle, Oracle-parameterized BSD, and dependence-adjusted
		procedures at $m=200$, $\alpha=0.10$, and $\psi^2=2\log(200)$. The six
		panels correspond to the one-factor loading patterns described in the
		text. The vertical axis is logarithmic; shaded bands are pointwise $95\%$
		Monte Carlo intervals.}
	\label{fig:simulation-risk-absolute}
\end{figure}

\begin{figure}[t]
	\centering
	\includegraphics[width=\textwidth]{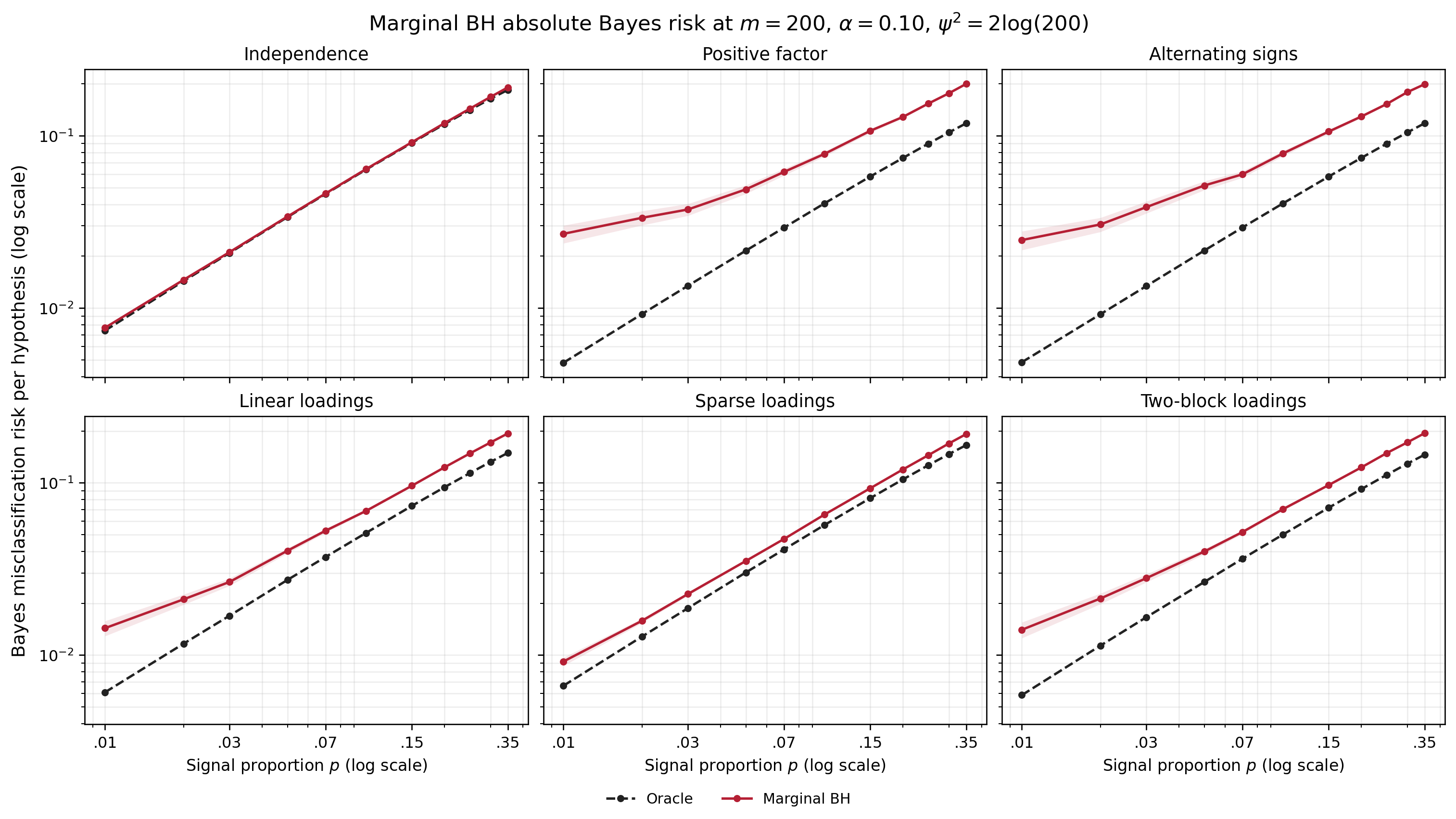}
	\caption{Absolute Bayes misclassification risk per hypothesis for marginal
		BH and the observed-data Oracle under the same simulation design. The
		vertical axis is logarithmic; shaded bands are pointwise $95\%$ Monte
		Carlo intervals.}
	\label{fig:simulation-risk-marginal-bh}
\end{figure}

Figure~\ref{fig:simulation-risk-absolute} displays absolute Bayes risks for
the Oracle, BSD, and dependence-adjusted procedures. Across all six loading
patterns, the BSD and Oracle curves are nearly indistinguishable. MRD--GBS,
rGBS, and residual BH also remain close to the Oracle, with their largest
finite-sample gaps occurring at the sparsest values of $p$ under the more
strongly dependent loading patterns. Under independence, MRD--GBS and rGBS
coincide, as expected from the absence of a factor adjustment.

For example, at $p=0.01$ under positive factor loadings, the estimated
per-hypothesis risks are approximately $0.00482$ for the observed-data
Oracle, $0.00482$ for BSD, $0.00536$ for MRD--GBS, $0.00540$ for rGBS, and
$0.00544$ for residual BH. The corresponding marginal-BH risk is about
$0.02697$. Figure~\ref{fig:simulation-risk-marginal-bh} shows that this
large separation between marginal BH and the Oracle is specific to
nontrivial dependence; under independence, marginal BH and residual BH
coincide.

The absolute risk generally rises with $p$ because the proportion of
non-null coordinates increases while $m$ and the signal variance remain
fixed. This increase should be read together with the procedure-to-Oracle
gaps, not as evidence that relative efficiency worsens. Since $m$ is fixed
throughout this experiment, these curves provide a finite-sample sparsity
profile and do not by themselves verify the asymptotic ABOS property.

\section{Discussion}
\label{sec:discussion}

This paper studies sparse multiple testing under known one-factor Gaussian
dependence through the Bayes risk of the observed-data experiment. We derive
the first-order risk of the Bayes Oracle based on the complete observation
vector and show that the factor-known experiment has the same leading risk
under the stated calibration. Relative to this observed-data benchmark, the
complete adaptive MRD--GBS procedure is ABOS. Fixed residualized GBS and BH
attain the same leading risk under their respective calibrations, and BSD
does so when its true mixture, slab, covariance, and loss parameters are
supplied.

The scope is integrated Bayes risk under the centered Gaussian two-groups
model, known one-factor covariance, and the deterministic calibrations in
the theorems. FDR control and ABOS are different targets: a result on the
integrated cost of false rejections and missed signals does not itself
provide a frequentist FDR guarantee, just as FDR control alone does not
establish first-order Bayes-risk efficiency. These results do not assert
finite-sample FDR control, cover arbitrary or estimated covariance, or
provide data-driven calibration.
The BSD theorem is additionally Oracle-parameterized, as specified in
Section~\ref{sec:bsd-results}.

The matching factor-known lower bound clarifies how dependence contributes
to this result. In the homogeneous model, adjustment for a nonzero common
factor reduces the relevant noise variance, which can lower the Oracle's
leading risk relative to the independent-coordinate benchmark. For
heterogeneous signed loadings, the coordinate-specific contributions remain
in the average $Q_m$, without requiring a limiting distribution of the
loadings. Thus, under the theorem's assumptions, uncertainty about the
latent factor does not increase the leading Bayes risk. This is a
first-order statement; it does not imply equality of the observed-data and
factor-known risks at finite dimension.

The results also connect the empirical and methodological themes of the
BSD and MRD--GBS manuscripts. The BSD study documents close agreement
between BSD and the Bayes Oracle in one-factor simulations. The MRD--GBS
study combines covariance-adaptive residualization with a systematic
stagewise calibration and reports favorable performance across several
dependence structures. The present ABOS theorem gives a first-order
explanation for the near-Oracle behavior observed in those one-factor
comparisons. It does not imply that BSD and MRD--GBS make identical
selections, have the same finite-sample risks, or are equivalent beyond the
conditions studied here.

Residualization and calibration play distinct roles. MRD constructs
residual evidence from the active covariance structure and recomputes it
after each rejection; the GBS constants determine how that evolving evidence
is translated into decisions. The adaptive proof accounts for the
data-selected path by controlling residual perturbations over capped active
sets. It does not assume that the adaptive procedure follows the same
rejection path as a fixed residual rule. The fixed residualized procedures,
in turn, help isolate the effect of applying different testing rules after
leave-one-out GLS adjustment. Their common leading risk with adaptive
MRD--GBS is an asymptotic risk statement, not a claim of identical
rejection sets or finite-sample behavior.

The fixed residualized-BH result overlaps with the centered-BH ABOS work of
\citet{RoyBhandari2024} in the homogeneous equicorrelated
setting. That overlap should be acknowledged: this paper does not claim the
first ABOS result for BH under equicorrelation. The present analysis uses
leave-one-out GLS residuals, treats heterogeneous signed loadings, and
establishes ABOS for the complete adaptive MRD--GBS path as well. Ordinary
BH applied to raw marginal $p$-values is a different procedure from
centered or GLS-residualized BH.

The simulation fixes $m=200$ and varies the signal proportion across six
loading patterns. Its absolute-risk curves describe finite-sample behavior
and connect the theory to the companion studies, but they cannot establish
convergence of an ABOS risk ratio.

These limitations also suggest a natural order for extending the theory.
The first question is whether ABOS can be established beyond one-factor
models while retaining the assumption of known covariance. Although MRD
is defined for general known positive definite covariance matrices, the
present proofs use the one-factor representation to derive the observed-data
Oracle risk and control residual errors along the adaptive path. A broader
theory would identify covariance structures for which the Oracle risk has
a tractable first-order expansion and the complete MRD--GBS procedure
attains it. Multiple-factor, block, and sparse-precision models offer
distinct settings for this investigation. Each would require connecting
the covariance geometry to the observed-data Oracle risk and controlling
the active sets visited by the procedure.

Once the ABOS theory has been extended to broader known covariance
structures, the next question is how covariance estimation changes the
benchmark and the procedure's risk. Relevant tools include covariance
shrinkage, regularized covariance estimation, adaptive thresholding, and
sparse precision-matrix estimation; see, for example,
\citet{LedoitWolf2004,BickelLevina2008,CaiLiu2011,
FriedmanHastieTibshirani2008}. The testing problem is to determine whether
estimation errors are negligible for both the Oracle risk and the residual
statistics at sparse-testing thresholds, uniformly over the active sets
visited by MRD--GBS. Estimation from an independent sample provides a
natural initial setting. Estimating covariance from the testing data
would also require accounting for its interaction with sequential
selection.

A further direction is to develop plug-in, empirical Bayes, or hierarchical
versions that retain first-order risk when model, loss, and calibration
inputs are estimated. A further question is whether a separately calibrated MRD variant can
pair Bayes-risk efficiency with a formal frequentist guarantee such as FDR
control.

Finally, ABOS under a specified prior does not imply worst-case optimality.
It would be useful to study sharp minimax risk for covariance-adaptive
procedures beyond independence and clarify its relationship with
admissibility under vector loss. The present results cover bounded
positive expected signal counts, but not the regime $mp_m\to0$, where
the Oracle normalization and the contribution of false discoveries may
change. Extensions to non-Gaussian slabs and noise distributions would
also require new tail and risk arguments.

These extensions form a progression: first broaden the dependence models
for which ABOS is known, then account for estimating covariance, and
thereafter study data-driven calibration and other optimality criteria.

\appendix
\section*{Appendix}

This appendix contains the proofs of the theoretical results. Appendix~A
derives the homogeneous observed-data Oracle risk. Appendices~B and~C establish the residual
correspondence and the adaptive MRD--GBS results. Appendix~D proves the
Oracle-parameterized BSD theorem. The later appendices treat fixed
residualized GBS, heterogeneous loadings, and bounded-count and
polynomial-sparsity results. Appendix~G proves the residual-BH result,
Appendix~H addresses known-covariance identifiability and general
vanishing sparsity.
\section{Homogeneous model and observed-data Oracle}
\label{app:homogeneous}

We consider the homogeneous one-factor sparse Gaussian model
\begin{equation}
	\label{eq:one_factor_model}
	X_i
	=
	\theta_i+\lambda Z+\sqrt{1-\lambda^2}\,\varepsilon_i,
	\qquad i=1,\ldots,m,
\end{equation}
where \(|\lambda|<1\), \(Z\sim N(0,1)\), and
\(\varepsilon_1,\ldots,\varepsilon_m\stackrel{\mathrm{i.i.d.}}{\sim}N(0,1)\),
independently of \(Z\). The latent indicators
\(\nu_1,\ldots,\nu_m\) are independent Bernoulli random variables with
\[
\mathbb P(\nu_i=1)=p_m,
\qquad
\mathbb P(\nu_i=0)=1-p_m,
\]
and
\[
\theta_i\mid \nu_i=0=0,
\qquad
\theta_i\mid \nu_i=1\sim N(0,\psi_m^2).
\]
The pairs \((\nu_i,\theta_i)\) are independent across coordinates and
independent of \(Z,\varepsilon_1,\ldots,\varepsilon_m\).
Put
\[
u_m=\psi_m^2,
\qquad
a_\lambda^2=1-\lambda^2,
\qquad
u_{m,\lambda}
=
\frac{u_m}{a_\lambda^2}
=
\frac{\psi_m^2}{1-\lambda^2}.
\]

Conditional on \(Z=z\), the null and non-null distributions of \(X_i\)
are
\[
X_i\mid Z=z,\nu_i=0
\sim
N(\lambda z,a_\lambda^2),
\]
and
\[
X_i\mid Z=z,\nu_i=1
\sim
N(\lambda z,a_\lambda^2+\psi_m^2).
\]
Equivalently, defining the residualized statistic
\begin{equation}
	\label{eq:residualized_statistic}
	Y_i(z)
	=
	\frac{X_i-\lambda z}{a_\lambda},
\end{equation}
we have
\[
Y_i(z)\mid Z=z,\nu_i=0
\sim
N(0,1),
\]
and
\[
Y_i(z)\mid Z=z,\nu_i=1
\sim
N(0,1+u_{m,\lambda}).
\]
Thus, conditional on the common factor, the problem reduces exactly to
the independent sparse Gaussian sequence model with effective signal
strength \(u_{m,\lambda}\).

Let
\[
f_m=\frac{1-p_m}{p_m},
\qquad
\delta_m=\frac{\delta_{0,m}}{\delta_{A,m}}.
\]
The conditional Bayes Oracle rejects \(H_{0i}\) whenever
\[
\frac{
	\phi\{Y_i(z);0,1+u_{m,\lambda}\}
}{
	\phi\{Y_i(z);0,1\}
}
>
\delta_m f_m.
\]
A straightforward calculation gives the equivalent rejection rule
\begin{equation}
	\label{eq:oracle_rejection_factor_known}
	Y_i(z)^2>c_{m,\lambda}^2,
\end{equation}
where
\begin{equation}
	\label{eq:oracle_threshold_factor_known}
	c_{m,\lambda}^2
	=
	\left(1+\frac1{u_{m,\lambda}}\right)
	\left[
	\log\{u_{m,\lambda}f_m^2\delta_m^2\}
	+
	\log\left(1+\frac1{u_{m,\lambda}}\right)
	\right].
\end{equation}

Consequently, the conditional Type-I and Type-II error probabilities of
the Bayes Oracle are
\begin{equation}
	\label{eq:type1_factor_known}
	t_{1,m}^{\lambda}
	=
	2\bar\Phi(c_{m,\lambda}),
\end{equation}
and
\begin{equation}
	\label{eq:type2_factor_known}
	t_{2,m}^{\lambda}
	=
	2\Phi\left(
	\frac{c_{m,\lambda}}{\sqrt{1+u_{m,\lambda}}}
	\right)-1.
\end{equation}
Therefore the conditional Bayes Oracle risk is
\begin{equation}
	\label{eq:conditional_oracle_risk}
	R_{m,\lambda}^{\mathrm{BO,cond}}
	=
	m\left[
	\delta_{0,m}(1-p_m)t_{1,m}^{\lambda}
	+
	\delta_{A,m}p_m t_{2,m}^{\lambda}
	\right].
\end{equation}
The right-hand side does not depend on \(z\). Hence the factor-known
Oracle risk is the same after integrating over the distribution of \(Z\).

Assume now that
\begin{equation}
	\label{eq:lambda_calibration}
	\frac{\log\{u_{m,\lambda}f_m^2\delta_m^2\}}
	{u_{m,\lambda}}
	\to
	C_\lambda\in(0,\infty).
\end{equation}
Then, as in the standard sparse Gaussian sequence calculation,
\[
\frac{c_{m,\lambda}^2}{1+u_{m,\lambda}}
\to
C_\lambda.
\]
Consequently,
\begin{equation}
	\label{eq:type2_factor_asymptotic}
	t_{2,m}^{\lambda}
	=
	2\Phi(\sqrt{C_\lambda})-1+o(1).
\end{equation}
Moreover, the Type-I contribution is asymptotically negligible relative
to the leading Type-II contribution under the usual sparse Bayes-risk
calibration. Hence
\begin{equation}
	\label{eq:factor_known_oracle_risk_asymptotic}
	R_{m,\lambda}^{\mathrm{BO,cond}}
	=
	mp_m\delta_{A,m}
	\left\{
	2\Phi(\sqrt{C_\lambda})-1
	\right\}
	\{1+o(1)\}.
\end{equation}

If the independence calibration satisfies
\begin{equation}
	\label{eq:independence_calibration}
	\frac{\log\{u_mf_m^2\delta_m^2\}}
	{u_m}
	\to C\in(0,\infty),
\end{equation}
then, for fixed \(|\lambda|<1\),
\[
C_\lambda
=
(1-\lambda^2)C.
\]
Therefore
\begin{equation}
	\label{eq:factor_known_oracle_risk_final}
	R_{m,\lambda}^{\mathrm{BO,cond}}
	=
	mp_m\delta_{A,m}
	\left[
	2\Phi\left\{
	\sqrt{(1-\lambda^2)C}
	\right\}
	-1
	\right]
	\{1+o(1)\}.
\end{equation}

We now turn to the genuine Bayes Oracle, which observes only
\(\boldsymbol{X}=(X_1,\ldots,X_m)\), and not the latent common factor \(Z\). The true
Bayes Oracle rejects \(H_{0i}\) whenever
\[
\mathbb P(\nu_i=1\mid \boldsymbol{X})
>
\frac{\delta_m}{1+\delta_m}.
\]
Equivalently, it rejects whenever the posterior odds exceed \(\delta_m\).

\subsection{Proof of Theorem~\ref{thm:one-factor-oracle-risk}}

\begin{proof}
	First consider the factor-known experiment, in which the common factor
	\(Z\) is observed. Conditional on \(Z=z\), define
	\begin{equation}
		\label{eq:proof-residual-statistic}
		Y_i(z)
		=
		\frac{X_i-\lambda z}{a_\lambda}.
	\end{equation}
	Then
	\begin{equation}
		\label{eq:proof-residual-distributions}
		Y_i(z)\mid \nu_i=0\sim N(0,1),
		\qquad
		Y_i(z)\mid \nu_i=1\sim N(0,1+u_{m,\lambda}).
	\end{equation}
	The factor-known Bayes Oracle rejects \(H_{0i}\) whenever
	\begin{equation}
		\label{eq:proof-likelihood-ratio-rule}
		\frac{
			\phi\{Y_i(z);0,1+u_{m,\lambda}\}
		}{
			\phi\{Y_i(z);0,1\}
		}
		>
		\delta_m f_m.
	\end{equation}
	Writing \(y=Y_i(z)\), \eqref{eq:proof-likelihood-ratio-rule} is equivalent to
	\begin{equation}
		\label{eq:proof-likelihood-ratio-expanded}
		\frac{1}{\sqrt{1+u_{m,\lambda}}}
		\exp\left[
		\frac{y^2}{2}
		\left\{
		1-\frac{1}{1+u_{m,\lambda}}
		\right\}
		\right]
		>
		\delta_m f_m.
	\end{equation}
	Hence the rejection region is
	\begin{equation}
		\label{eq:proof-oracle-rejection-region}
		y^2>c_{m,\lambda}^2,
	\end{equation}
	where
	\begin{equation}
		\label{eq:proof-one-factor-oracle-threshold}
		c_{m,\lambda}^2
		=
		\left(1+\frac1{u_{m,\lambda}}\right)
		\left[
		\log\{u_{m,\lambda}f_m^2\delta_m^2\}
		+
		\log\left(1+\frac1{u_{m,\lambda}}\right)
		\right].
	\end{equation}
	
	Therefore the factor-known Type-I and Type-II error probabilities are
	\begin{equation}
		\label{eq:proof-factor-known-type1}
		t_{1,m}^{\lambda}
		=
		2\bar\Phi(c_{m,\lambda}),
	\end{equation}
	and
	\begin{equation}
		\label{eq:proof-factor-known-type2}
		t_{2,m}^{\lambda}
		=
		2\Phi\left(
		\frac{c_{m,\lambda}}{\sqrt{1+u_{m,\lambda}}}
		\right)-1.
	\end{equation}
	Hence
	\begin{equation}
		\label{eq:proof-factor-known-risk}
		R_{m,\lambda}^{\mathrm{BO,cond}}
		=
		m\left[
		\delta_{0,m}(1-p_m)t_{1,m}^{\lambda}
		+
		\delta_{A,m}p_m t_{2,m}^{\lambda}
		\right].
	\end{equation}
	
	By \eqref{eq:one-factor-calibration} and
	\eqref{eq:proof-one-factor-oracle-threshold},
	\begin{equation}
		\label{eq:proof-threshold-ratio-limit}
		\frac{c_{m,\lambda}^2}{1+u_{m,\lambda}}
		\to
		C_\lambda.
	\end{equation}
	Consequently,
	\begin{equation}
		\label{eq:proof-factor-known-type2-asymptotic}
		t_{2,m}^{\lambda}
		=
		2\Phi(\sqrt{C_\lambda})-1+o(1).
	\end{equation}
	Moreover, by the standard Mills-ratio calculation under
	\eqref{eq:one-factor-calibration},
	\begin{equation}
		\label{eq:proof-factor-known-type1-negligible}
		m\delta_{0,m}(1-p_m)t_{1,m}^{\lambda}
		=
		o(mp_m\delta_{A,m}).
	\end{equation}
	Combining \eqref{eq:proof-factor-known-risk},
	\eqref{eq:proof-factor-known-type2-asymptotic}, and
	\eqref{eq:proof-factor-known-type1-negligible}, we obtain
	\begin{equation}
		\label{eq:proof-factor-known-risk-asymptotic}
		R_{m,\lambda}^{\mathrm{BO,cond}}
		=
		mp_m\delta_{A,m}
		\left\{
		2\Phi(\sqrt{C_\lambda})-1
		\right\}
		\{1+o(1)\}.
	\end{equation}
Since the factor-known Oracle is permitted to use the larger information
set consisting of \((\boldsymbol{X},Z)\), whereas the genuine Bayes Oracle is
restricted to observing only \(\boldsymbol{X}\), the former minimizes the Bayes risk
over a larger class of decision rules. Consequently,
\begin{equation}
	\label{eq:proof-oracle-lower-bound}
	R_{m,\lambda}^{\mathrm{BO,cond}}
	\le
	R_{m,\lambda}^{\mathrm{BO}}.
\end{equation}

Motivated by the factor-known Bayes Oracle, we construct an
\(\boldsymbol{X}\)-measurable plug-in decision rule by replacing the unobservable
residual \(Y_i(Z)\) with its leave-one-out estimator \(\widehat Y_i\),
while retaining the same Oracle threshold \(c_{m,\lambda}\). This rule is
used only to prove the observed-data Oracle upper bound; its threshold
depends on the model and loss parameters and is not a proposed data-driven
testing rule. Define
\begin{equation}
	\label{eq:proof-leave-one-out-factor-estimator}
	\widehat Z_{-i}
	=
	\frac{1}{\lambda(m-1)}
	\sum_{j\ne i}X_j
\end{equation}
and
\begin{equation}
	\label{eq:proof-plugin-residual}
	\widehat Y_i
	=
	\frac{X_i-\lambda\widehat Z_{-i}}{a_\lambda}.
\end{equation}
	Since
	\begin{equation}
		\label{eq:proof-factor-estimator-expansion}
		\widehat Z_{-i}
		=
		Z
		+
		\frac{1}{\lambda(m-1)}
		\sum_{j\ne i}\theta_j
		+
		\frac{a_\lambda}{\lambda(m-1)}
		\sum_{j\ne i}\varepsilon_j,
	\end{equation}
	we have
	\begin{equation}
		\label{eq:proof-plugin-residual-expansion}
		\widehat Y_i
		=
		Y_i(Z)-\Delta_{-i},
	\end{equation}
	where
	\begin{equation}
		\label{eq:proof-delta-definition}
		\Delta_{-i}
		=
		\frac{1}{m-1}\sum_{j\ne i}\varepsilon_j
		+
		\frac{1}{a_\lambda(m-1)}
		\sum_{j\ne i}\theta_j.
	\end{equation}
	By construction, \(\Delta_{-i}\) is independent of \((Y_i(Z),\nu_i)\).
	

Accordingly, we consider the plug-in decision rule that rejects
\(H_{0i}\) whenever
\begin{equation}
	\label{eq:proof-plugin-rejection-rule}
	\widehat Y_i^2>c_{m,\lambda}^2,
\end{equation}
where \(c_{m,\lambda}\) is the factor-known Oracle threshold defined in
\eqref{eq:proof-one-factor-oracle-threshold}. Let
\(R_{m,\lambda}^{\mathrm{plug}}\) denote the corresponding Bayes risk. Since the genuine
	Bayes Oracle minimizes Bayes risk among all \(\boldsymbol{X}\)-measurable rules,
	\begin{equation}
		\label{eq:proof-oracle-upper-bound}
		R_{m,\lambda}^{\mathrm{BO}}
		\le
		R_{m,\lambda}^{\mathrm{plug}}.
	\end{equation}
	
	Since
	\begin{equation}
		\label{eq:proof-mu-second-moment}
		\mathbb E(\theta_j^2)=p_m\psi_m^2
	\end{equation}
	and the two sums in \(\Delta_{-i}\) are independent and centered,
	\begin{equation}
		\label{eq:proof-delta-second-moment}
		\mathbb E(\Delta_{-i}^2)
		=
		\frac{1}{m-1}
		+
		\frac{p_m u_{m,\lambda}}{m-1}.
	\end{equation}
	Also, \eqref{eq:one-factor-calibration} and
	\eqref{eq:proof-one-factor-oracle-threshold} imply
	\begin{equation}
		\label{eq:proof-threshold-order}
		c_{m,\lambda}^2=O(u_{m,\lambda}).
	\end{equation}
	Hence, by \eqref{eq:factor-estimation-condition},
	\begin{align}
		c_{m,\lambda}^2\mathbb E(\Delta_{-i}^2)
		&=
		O\left(
		\frac{
			u_{m,\lambda}\{1+p_m u_{m,\lambda}\}
		}{m}
		\right)
		\notag\\
		&=
		o(1).
		\label{eq:proof-factor-error-threshold-second-moment}
	\end{align}
	Therefore
	\begin{equation}
		\label{eq:proof-factor-error-threshold-scale}
		c_{m,\lambda}\Delta_{-i}
		=
		o_{\mathbb P}(1).
	\end{equation}
	Similarly,
	\begin{equation}
		\label{eq:proof-factor-error-signal-second-moment}
		\frac{\mathbb E(\Delta_{-i}^2)}
		{1+u_{m,\lambda}}
		\le
		\frac{1+p_m u_{m,\lambda}}
		{(m-1)(1+u_{m,\lambda})}
		=
		o(1),
	\end{equation}
	and hence
	\begin{equation}
		\label{eq:proof-factor-error-signal-scale}
		\frac{\Delta_{-i}}{\sqrt{1+u_{m,\lambda}}}
		=
		o_{\mathbb P}(1).
	\end{equation}
	
	Under the alternative,
	\(Y_i(Z)\sim N(0,1+u_{m,\lambda})\), independently of
	\(\Delta_{-i}\). Hence
	\begin{equation}
		\label{eq:proof-plugin-type2-definition}
		t_{2,m}^{\mathrm{plug}}
		=
		\mathbb P\left(
		|Y_i(Z)-\Delta_{-i}|
		\le
		c_{m,\lambda}
		\mid \nu_i=1
		\right).
	\end{equation}
	Using \eqref{eq:proof-factor-error-signal-scale} and bounded convergence,
	\begin{equation}
		\label{eq:proof-plugin-type2-asymptotic-first}
		t_{2,m}^{\mathrm{plug}}
		=
		2\Phi\left(
		\frac{c_{m,\lambda}}{\sqrt{1+u_{m,\lambda}}}
		\right)-1
		+
		o(1).
	\end{equation}
	Thus, by \eqref{eq:proof-factor-known-type2-asymptotic},
	\begin{equation}
		\label{eq:proof-plugin-type2-asymptotic}
		t_{2,m}^{\mathrm{plug}}
		=
		2\Phi(\sqrt{C_\lambda})-1+o(1).
	\end{equation}
	
	We now control the Type-I error of the plug-in rule. Under the null,
	\(Y_i(Z)\sim N(0,1)\), independently of \(\Delta_{-i}\). Let
	\[
	K_{-i}
	=
	\sum_{j\ne i}\nu_j.
	\]
	Conditional on \(K_{-i}=k\), the random variable \(\Delta_{-i}\) is
	\(N(0,\tau_{m,k}^2)\), where
	\begin{equation}
		\label{eq:proof-delta-conditional-variance}
		\tau_{m,k}^2
		=
		\frac{1}{m-1}
		+
		\frac{k u_{m,\lambda}}{(m-1)^2}.
	\end{equation}
	Therefore,
	\begin{equation}
		\label{eq:proof-plugin-null-conditional-law}
		Y_i(Z)-\Delta_{-i}
		\mid K_{-i}=k,\nu_i=0
		\sim
		N(0,1+\tau_{m,k}^2).
	\end{equation}
	Hence the plug-in Type-I probability is
	\begin{equation}
		\label{eq:proof-plugin-type1-exact}
		t_{1,m}^{\mathrm{plug}}
		=
		2\mathbb E
		\left[
		\bar\Phi\left(
		\frac{c_{m,\lambda}}
		{\sqrt{1+\tau_{m,K_{-i}}^2}}
		\right)
		\right].
	\end{equation}
	Since \(\tau_{m,K_{-i}}^2\ge0\),
	\begin{equation}
		\label{eq:proof-plugin-type1-lower}
		t_{1,m}^{\mathrm{plug}}
		\ge
		t_{1,m}^{\lambda}.
	\end{equation}
	For the reverse comparison, Mills' ratio implies, uniformly for
	\(0\le s=o(1)\),
	\begin{equation}
		\label{eq:proof-normal-tail-variance-ratio}
		\frac{
			\bar\Phi\{c_{m,\lambda}/\sqrt{1+s}\}
		}{
			\bar\Phi(c_{m,\lambda})
		}
		\le
		\{1+o(1)\}
		\sqrt{1+s}
		\exp\left\{
		\frac{c_{m,\lambda}^2s}{2(1+s)}
		\right\}.
	\end{equation}
	Since
	\[
	\tau_{m,K_{-i}}^2
	\le
	\frac{1}{m-1}
	+
	\frac{u_{m,\lambda}}{m-1}
	=
	o(1),
	\]
	the preceding bound applies uniformly with
	\(s=\tau_{m,K_{-i}}^2\). Also,
	\[
	\sqrt{1+s}\le \exp(s/2),
	\qquad
	\frac{s}{1+s}\le s.
	\]
	Thus
	\begin{equation}
		\label{eq:proof-tail-ratio-exponential-bound}
		\frac{
			\bar\Phi\{c_{m,\lambda}/\sqrt{1+\tau_{m,K_{-i}}^2}\}
		}{
			\bar\Phi(c_{m,\lambda})
		}
		\le
		\{1+o(1)\}
		\exp\left\{
		\frac{1+c_{m,\lambda}^2}{2}
		\tau_{m,K_{-i}}^2
		\right\}.
	\end{equation}
	We now show that the expectation of the exponential factor is \(1+o(1)\).
	Since \(K_{-i}\sim\mathrm{Bin}(m-1,p_m)\),
	\begin{align}
		&
		\mathbb E
		\exp\left\{
		\frac{1+c_{m,\lambda}^2}{2}
		\tau_{m,K_{-i}}^2
		\right\}
		\notag\\
		&=
		\exp\left\{
		\frac{1+c_{m,\lambda}^2}{2(m-1)}
		\right\}
		\left[
		1-p_m
		+
		p_m
		\exp\left\{
		\frac{(1+c_{m,\lambda}^2)u_{m,\lambda}}
		{2(m-1)^2}
		\right\}
		\right]^{m-1}.
		\label{eq:proof-plugin-type1-mgf}
	\end{align}
	By \eqref{eq:proof-threshold-order} and
	\eqref{eq:factor-estimation-condition},
	\begin{equation}
		\label{eq:proof-mgf-first-term}
		\frac{1+c_{m,\lambda}^2}{m}
		=
		O\left(\frac{u_{m,\lambda}}{m}\right)
		=
		o(1),
	\end{equation}
	and
	\begin{align}
		&(m-1)p_m
		\left[
		\exp\left\{
		\frac{(1+c_{m,\lambda}^2)u_{m,\lambda}}
		{2(m-1)^2}
		\right\}
		-1
		\right]
		\notag\\
		&=
		O\left(
		\frac{p_m u_{m,\lambda}^2}{m}
		\right)
		=
		o(1),
		\label{eq:proof-mgf-second-term}
	\end{align}
	where the last step again follows from
	\eqref{eq:factor-estimation-condition}. Therefore
	\begin{equation}
		\label{eq:proof-plugin-type1-mgf-one}
		\mathbb E
		\exp\left\{
		\frac{1+c_{m,\lambda}^2}{2}
		\tau_{m,K_{-i}}^2
		\right\}
		=
		1+o(1).
	\end{equation}
	Combining \eqref{eq:proof-plugin-type1-exact},
	\eqref{eq:proof-tail-ratio-exponential-bound}, and
	\eqref{eq:proof-plugin-type1-mgf-one}, we obtain
	\[
	t_{1,m}^{\mathrm{plug}}
	\le
	t_{1,m}^{\lambda}\{1+o(1)\}.
	\]
	Together with \eqref{eq:proof-plugin-type1-lower}, this gives
	\begin{equation}
		\label{eq:proof-plugin-type1-equivalence}
		t_{1,m}^{\mathrm{plug}}
		=
		t_{1,m}^{\lambda}\{1+o(1)\}.
	\end{equation}
	Therefore, by \eqref{eq:proof-factor-known-type1-negligible},
	\begin{equation}
		\label{eq:proof-plugin-type1-negligible}
		m\delta_{0,m}(1-p_m)t_{1,m}^{\mathrm{plug}}
		=
		o(mp_m\delta_{A,m}).
	\end{equation}
	
	Combining \eqref{eq:proof-plugin-type2-asymptotic} and
	\eqref{eq:proof-plugin-type1-negligible},
	\begin{align}
		R_{m,\lambda}^{\mathrm{plug}}
		&=
		m\left[
		\delta_{0,m}(1-p_m)t_{1,m}^{\mathrm{plug}}
		+
		\delta_{A,m}p_m t_{2,m}^{\mathrm{plug}}
		\right]
		\notag\\
		&=
		mp_m\delta_{A,m}
		\left\{
		2\Phi(\sqrt{C_\lambda})-1
		\right\}
		\{1+o(1)\}.
		\label{eq:proof-plugin-risk-asymptotic}
	\end{align}
	Finally, by \eqref{eq:proof-oracle-lower-bound} and
	\eqref{eq:proof-oracle-upper-bound},
	\begin{equation}
		\label{eq:proof-risk-sandwich}
		R_{m,\lambda}^{\mathrm{BO,cond}}
		\le
		R_{m,\lambda}^{\mathrm{BO}}
		\le
		R_{m,\lambda}^{\mathrm{plug}}.
	\end{equation}
	Combining \eqref{eq:proof-risk-sandwich} with
	\eqref{eq:proof-factor-known-risk-asymptotic} and
	\eqref{eq:proof-plugin-risk-asymptotic} gives
	\[
	R_{m,\lambda}^{\mathrm{BO}}
	=
	mp_m\delta_{A,m}
	\left\{
	2\Phi(\sqrt{C_\lambda})-1
	\right\}
	\{1+o(1)\}.
	\]
	This completes the proof.
\end{proof}

\section{MRD--GBS: residual correspondence and proof}

Write \(\rho=\lambda^2\in(0,1)\), so that
\[
\boldsymbol{\Sigma}_m=(1-\rho)\boldsymbol{I}_{m}+\rho \boldsymbol{J}_{m}.
\]
For an active set \(A\subseteq\{1,\ldots,m\}\) of size \(s\ge2\), define the
active-set leave-one-out residual
\[
L_{i,A}
=
\frac{1}{\sqrt{1-\rho}}
\left(
X_i-\frac{1}{s-1}\sum_{j\in A,\,j\ne i}X_j
\right),
\qquad i\in A.
\]
When \(A=\{1,\ldots,m\}\), this is exactly the residual statistic
\(\widehat Y_i\) used by rGBS.

\begin{lemma}[Exact MRD--leave-one-out decomposition]
\label{lem:mrd-loo-decomposition}
For an active set \(A\) of size \(s\ge2\), let \(U_{i,A}^{\mathrm{MRD}}\)
denote the standardized MRD residual for coordinate \(i\in A\). Put
\[
\beta_s=\frac{\rho}{1+(s-2)\rho},
\qquad
h_s=\left\{\frac{1+(s-1)\rho}{1+(s-2)\rho}\right\}^{1/2},
\]
and
\[
d_s=\frac{1}{s-1}-\beta_s
=\frac{1-\rho}{(s-1)\{1+(s-2)\rho\}}.
\]
Then, exactly,
\begin{equation}
\label{eq:mrd-loo-exact}
U_{i,A}^{\mathrm{MRD}}
=
\frac{1}{h_s}
\left[
L_{i,A}
+
\frac{d_s}{\sqrt{1-\rho}}
\sum_{j\in A,\,j\ne i}X_j
\right].
\end{equation}
Consequently, for fixed \(\rho\in(0,1)\),
\begin{equation}
\label{eq:mrd-loo-bound}
|U_{i,A}^{\mathrm{MRD}}-L_{i,A}|
\le
\frac{C_\rho}{s}|L_{i,A}|
+
\frac{C_\rho}{s^2}
\left|\sum_{j\in A,\,j\ne i}X_j\right|
\end{equation}
for a constant \(C_\rho<\infty\) independent of \(A,i,s\).
\end{lemma}

\begin{proof}
For the equicorrelation matrix on \(A\),
\[
\boldsymbol{\Sigma}_A^{-1}
=
\frac{1}{1-\rho}
\left[
\boldsymbol{I}_{s}-\frac{\rho}{1+(s-1)\rho}\boldsymbol{J}_{s}
\right].
\]
Hence the MRD conditional-residual numerator is proportional to
\[
X_i-\beta_s\sum_{j\in A,\,j\ne i}X_j,
\qquad
\beta_s=\frac{\rho}{1+(s-2)\rho},
\]
while its Gaussian error conditional variance, holding the mean vector fixed, is
\[
\operatorname{Var}(X_i\mid \boldsymbol{X}_{A\setminus\{i\}},\boldsymbol{\theta}_A)
=(1-\rho)\frac{1+(s-1)\rho}{1+(s-2)\rho}
=(1-\rho)h_s^2.
\]
Since
\[
\beta_s=\frac{1}{s-1}-d_s,
\]
substitution gives \eqref{eq:mrd-loo-exact}.  Finally,
\[
|h_s^{-1}-1|=O(s^{-1}),
\qquad d_s=O(s^{-2}),
\]
uniformly for fixed \(\rho\), which yields \eqref{eq:mrd-loo-bound}.
\end{proof}

\subsection{Proof of Corollary~\ref{cor:mrd-loo-uniform}}

\begin{proof}
For every \(A\) with \(|A|=s\ge m/2\),
\[
|L_{i,A}|
\le \frac{2}{\sqrt{1-\rho}}\max_{1\le j\le m}|X_j|
\]
and
\[
\left|\sum_{j\in A,\,j\ne i}X_j\right|
\le s\max_{1\le j\le m}|X_j|.
\]
Thus Lemma~\ref{lem:mrd-loo-decomposition} implies
\[
\sup_{A:\,|A|\ge m/2}\max_{i\in A}
|U_{i,A}^{\mathrm{MRD}}-L_{i,A}|
\le \frac{C_\rho}{m}\max_{1\le j\le m}|X_j|.
\]
Under the homogeneous one-factor model,
\[
\max_j|X_j|
=O_{\mathbb P}\left(1+\sqrt{(1+u_{m,\lambda})\log m}\right).
\]
The ABOS assumptions give \(u_{m,\lambda}=O(\log m)\) and
\(z_{c_1(\alpha_m)}=O(\sqrt{\log m})\). Therefore the preceding display,
after multiplication by \(z_{c_1(\alpha_m)}\), is
\[
O_{\mathbb P}\left(\frac{\log^{3/2}m}{m}\right)=o_{\mathbb P}(1).
\]
\end{proof}

\begin{remark}[Residual correspondence and adaptive optimality]
The preceding corollary compares MRD residuals with active-set leave-one-out
residuals and is not itself a transfer of complete rejection paths.
Appendix~\ref{app:adaptive-mrd} proves adaptive ABOS by a separate argument,
using a deterministic cap, a rejection-set embedding and risk-weighted
exceptional-event bounds. The approximation is matched to the stage
and imposes no additional sparsity rate or assumed short selected path.
\end{remark}


\section{Adaptive MRD--GBS: proofs without an additional sparsity rate}
\label{app:adaptive-mrd}
The notation is that of Section~\ref{sec:adaptive-mrd-results}.
Choose any loading representation recovered from the known covariance
by Lemma~\ref{lem:known-sigma-loadings}. A common sign change is absorbed
by $Z\mapsto-Z$, preserving its distribution and independence. Thus the
factor representation used below is a proof device compatible with the
same observed-data model and the same implemented algorithm.
Throughout, $p_m=o(1)$ means exactly $p_m\to0$. No additional logarithmic
sparsity rate is assumed: \eqref{mrd:eq:automatic-rate} below is a
consequence of this limit.
For $A\subseteq[m]$, put $\boldsymbol{D}_A=\operatorname{diag}(\sigma_i^2:i\in A)$
and $\boldsymbol{\lambda}_A=(\lambda_i:i\in A)^T$, so
$\boldsymbol{\Sigma}_A=\boldsymbol{D}_A+\boldsymbol{\lambda}_A\boldsymbol{\lambda}_A^T$.
The cap is an analytical device; the algorithm is not truncated and does
not use $p_m$ in its residual computations. The perturbation argument matches the number of possible deletions to
the threshold at each stage, instead of applying the first threshold
to every set within the cap.

\paragraph{Proof roadmap and assumption dependencies.}
The common structural assumptions are the known one-factor model
\eqref{eq:framework-model}--\eqref{eq:framework-loadings}; they give
uniformly positive idiosyncratic variances and independent coordinates
after conditioning on the factor. The two risk regimes add different
calibrations: Theorem~\ref{mrd:thm:diverging} uses
\eqref{mrd:eq:div-counts}--\eqref{mrd:eq:div-calibration}, whereas
Theorem~\ref{mrd:thm:bounded} assumes $mp_m$ is bounded above and away
from zero, $\psi_m^2/\log m\to2/C_0$, fixed positive loss ratio,
$\alpha_m\to0$, and $\log(1/\alpha_m)=o(\log m)$. The proof uses these assumptions in
four steps.
\begin{enumerate}[leftmargin=2em]
\item Lemma~\ref{mrd:lem:cap} uses only $p_m\to0$ and a positive lower
  bound on $mp_m$ to choose a deterministic cap larger than the expected
  number of signals but still $o(m)$. Lemma~\ref{mrd:lem:uniform} adds
  the loading bounds, $\alpha_m\to0$, $\alpha_m/p_m\to\infty$, and
  $\psi_m^2=O\{\log(e/p_m)\}$ to control residual perturbations
  uniformly over every active set within that cap. The condition
  $\alpha_m\delta_m\to0$ is used when converting its exceptional-event
  bound to a loss-weighted risk bound.
\item Lemmas~\ref{mrd:lem:transfers} and~\ref{mrd:lem:inclusion} turn
  that uniform event into deterministic threshold transfers and embed
  adaptive rejections in a static factor-known GBS rule.
  Lemma~\ref{mrd:lem:cap-tail} controls the event that this static rule reaches
  the cap. Together with the independent-null GBS bound in
  Lemma~\ref{mrd:lem:static-V}, these results yield the Type-I bounds in
  Lemma~\ref{mrd:lem:typeI}; no equality of the two rejection paths is
  used.
\item For diverging counts, Lemma~\ref{mrd:lem:crossing} uses the first
  signal calibration in \eqref{mrd:eq:div-calibration}, $mp_m\to\infty$,
  and the loading bounds to show that enough factor-known signal
  $p$-values cross the GBS thresholds. Lemma~\ref{mrd:lem:adaptive-crossing}
  transfers this crossing to the adaptive path using the same uniform
  event. For bounded counts, the proof instead counts signals crossing
  the first-stage threshold and bounds the exceptional contribution by
  Cauchy--Schwarz.
\item Lemma~\ref{mrd:lem:oracle} uses the signal, level, and loss
  calibrations in either theorem to obtain the factor-known Oracle risk.
  The information ordering
  \eqref{mrd:eq:information-order}, this risk expansion, and the Type-I
  and Type-II bounds from the preceding steps give the observed-data
  Oracle sandwich and the asserted ABOS limits.
\end{enumerate}

\subsection{A uniform approximation over a capped collection of active sets}

The cap below is deterministic; it is not an assumed bound on the random
MRD path. The probability of reaching it will be proved small separately.

\begin{lemma}[A cap under unrestricted vanishing sparsity]\label{mrd:lem:cap}
Suppose $p_m\to0$ and $mp_m\ge b_*>0$ eventually. Define
\begin{equation}\label{mrd:eq:cap}
L_m=\lfloor m\sqrt{p_m}\rfloor.
\end{equation}
Then $L_m\to\infty$, $L_m=o(m)$, $L_m/(mp_m)\to\infty$ and
$L_m\ge c\sqrt m$ eventually. For $b_m=\log(e/p_m)$,
\begin{equation}\label{mrd:eq:automatic-rate}
\sqrt{p_m}\,b_m^{3/2}\to0.
\end{equation}
\end{lemma}
\begin{proof}
The unrounded cap is at least $\sqrt{b_*m}$; hence integer rounding is
negligible. Its ratio to $m$ is $\sqrt{p_m}\to0$, and its ratio to
$mp_m$ is $p_m^{-1/2}\to\infty$. Finally, if $v_m=\log(1/p_m)\to\infty$,
then $\sqrt{p_m}b_m^{3/2}=e^{-v_m/2}(1+v_m)^{3/2}\to0$.
This last limit is automatic, not an additional assumption.
\end{proof}

Define the factor-known independent residuals and weighted residual terms
\[
Y_i=\theta_i/\sigma_i+\varepsilon_i,\qquad
B_i=\frac{\lambda_i}{\sigma_i}Y_i,
\qquad \widetilde P_i=2\Phibar(|Y_i|).
\]
The $Y_i$ (and the $B_i$) are independent across coordinates, and
$Y_i\mid\nu_i=0\sim N(0,1)$,
$Y_i\mid\nu_i=1\sim N(0,1+\psi_m^2/\sigma_i^2)$.
Let $z_t=\Phi^{-1}(1-t/2)$.

\begin{lemma}[Exact active-set MRD formula]\label{mrd:lem:exact}
For $i\in A$, put
\[
T_{A,-i}=\sum_{j\in A\setminus\{i\}}\lambda_j^2/\sigma_j^2,
\quad N_{A,-i}=\sum_{j\in A\setminus\{i\}}B_j,
\quad h_{i,A}=\left(1+\frac{\lambda_i^2/\sigma_i^2}{1+T_{A,-i}}\right)^{1/2}.
\]
Then, exactly,
\begin{equation}\label{mrd:eq:exact-active}
U_{i,A}=h_{i,A}^{-1}
\left[Y_i+\frac{\lambda_i}{\sigma_i(1+T_{A,-i})}
\{Z-N_{A,-i}\}\right].
\end{equation}
Under \eqref{eq:framework-loadings}, uniformly for $m-|A|\le L_m=o(m)$ and
$i\in A$, $T_{A,-i}\ge cm$ and
\begin{equation}\label{mrd:eq:deterministic-error}
\sup_{A:\,m-|A|\le L_m}\max_{i\in A}|U_{i,A}-Y_i|
\le\frac C m\left[|Z|+\left|\sum_{j=1}^m B_j\right|
+(L_m+2)\max_j|Y_j|\right].
\end{equation}
Constants may depend on the loading bounds, but not on $A,i,m$.
\end{lemma}
\begin{proof}
Write $T_A=\sum_{j\in A}\lambda_j^2/\sigma_j^2$. Sherman--Morrison gives
\[
\boldsymbol{\Omega}_A=\boldsymbol{D}_A^{-1}-
\frac{\boldsymbol{D}_A^{-1}\boldsymbol{\lambda}_A\boldsymbol{\lambda}_A^T\boldsymbol{D}_A^{-1}}{1+T_A}.
\]
It follows, by expanding the diagonal and numerator, that
\[
U_{i,A}=
\frac{X_i-\lambda_i(1+T_{A,-i})^{-1}
\sum_{j\in A\setminus\{i\}}\lambda_j X_j/\sigma_j^2}
{\sigma_i h_{i,A}}.
\]
Now $\sum_{j\in A\setminus\{i\}}\lambda_jX_j/\sigma_j^2
=T_{A,-i}Z+N_{A,-i}$ and $X_i=\sigma_iY_i+\lambda_iZ$.
Substitution proves \eqref{mrd:eq:exact-active}.

By \eqref{eq:framework-loadings}, the full sum
$\sum_j\lambda_j^2/\sigma_j^2\ge cm$. Each summand is bounded above,
so deletion of at most $L_m+1=o(m)$ terms leaves a sum at least $c'm$.
Consequently $|h_{i,A}^{-1}-1|\le C/m$ and
$|\lambda_i|/\{\sigma_i(1+T_{A,-i})\}\le C/m$ uniformly.
Furthermore,
\[
|N_{A,-i}|\le\left|\sum_j B_j\right|+(L_m+1)\max_j|B_j|,
\qquad \max_j|B_j|\le C\max_j|Y_j|.
\]
The triangle inequality in \eqref{mrd:eq:exact-active} now yields
\eqref{mrd:eq:deterministic-error}. No conditioning on the selected active
set has been used; the inequality holds simultaneously for every such set.
\end{proof}

\begin{lemma}[Quantitative stage-matched uniform perturbation]
\label{mrd:lem:uniform}
Assume \eqref{eq:framework-loadings}, $mp_m\ge b_*>0$, $p_m\to0$,
$\alpha_m\to0$, $\alpha_m/p_m\to\infty$ and
$\psi_m^2=O\{\log(e/p_m)\}$. For fixed $a,B>0$, define
\begin{equation}\label{mrd:eq:stage-event}
\mathcal E_m(a)=\left\{
\max_{1\le k\le L_m} z_{c_k(\alpha_m)}
\sup_{A:\,m-|A|\le k-1}\max_{i\in A}|U_{i,A}-Y_i|\le a\right\}.
\end{equation}
There are constants $c_a,C_{a,B}>0$ such that, eventually,
\begin{equation}\label{mrd:eq:uniform-tail}
\Pr\{\mathcal E_m(a)^c\}\le C_{a,B}m^{-B}
+2\exp\{-c_am/(\log m)^2\}.
\end{equation}
If also $\alpha_m\delta_m\to0$, then
\begin{equation}\label{mrd:eq:uniform-risk}
(1+\delta_m)m\Pr\{\mathcal E_m(a)^c\}=o(mp_m).
\end{equation}
No sparsity rate beyond $p_m\to0$ and the displayed count assumption
is needed for this lemma.
\end{lemma}
\begin{proof}
Put $b_m=\log(e/p_m)$ and $A_m=\log(1/\alpha_m)$.
The assumptions imply $p_m\ge b_*/m$ and eventually $\alpha_m\ge p_m$,
so $b_m=O(\log m)$, $0\le A_m\le b_m$ and
$z_{c_1(\alpha_m)}^2=O(\log m)$.
The independent centered variables $Y_i$ and $B_i$ have sub-Gaussian
variance proxies bounded by $Cb_m$. Indeed, for every real $t$,
\[
\mathbb E e^{tY_i}=(1-p_m)e^{t^2/2}
+p_m e^{(1+\psi_m^2/\sigma_i^2)t^2/2}
\le e^{Cb_mt^2/2},
\]
and the loading bounds give the same conclusion for
$B_i=(\lambda_i/\sigma_i)Y_i$.

For $1\le k\le L_m$, let
$H_{m,k}=\sup_{|J|\le k}\sum_{j\in J}|B_j|$ and
$t_{m,k}=\log(em/k)$. By extending a smaller set to size $k$ and
maximizing over signs,
\[
H_{m,k}=\max_{|J|=k}\max_{s\in\{-1,1\}^{J}}\sum_{j\in J}s_jB_j.
\]
Each signed sum has variance proxy at most $Ckb_m$. Since
$2^k\binom mk\le\exp\{Ck\log(em/k)\}$, its exponential tail and a
union bound give, for a sufficiently large $C_B$,
\[
\Pr\left\{H_{m,k}>
C_B\sqrt{b_m}\left[k\sqrt{t_{m,k}}+\sqrt{k\log m}\right]\right\}
\le m^{-(B+2)}.
\]
A further union bound over at most $m$ possible $k$ proves that,
except on probability at most $m^{-(B+1)}$, simultaneously for all
$k\le L_m$,
\begin{equation}\label{mrd:eq:deleted-stage}
H_{m,k}\le C_B\sqrt{b_m}
\{k\sqrt{t_{m,k}}+\sqrt{k\log m}\}.
\end{equation}
This bound integrates over the centered Gaussian slab; it is not a
uniform assertion over arbitrary deterministic signal means.

If $m-|A|\le k-1$, the set excluded from $N_{A,-i}$ has at most $k$
elements, including $i$. The exact residual formula therefore gives
\begin{equation}\label{mrd:eq:stage-error}
\sup_{A:\,m-|A|\le k-1}\max_{i\in A}|U_{i,A}-Y_i|
\le\frac C m\left\{|Z|+\left|\sum_iB_i\right|
+\max_i|Y_i|+H_{m,k}\right\}.
\end{equation}
The denominators are uniformly of order $m$ because $L_m=o(m)$;
no conditioning on a selected set is involved.
Also $c_k(\alpha_m)\ge k\alpha_m/(m+1)$, so the Gaussian tail bound
$z_t^2\le2\log(2/t)$ gives
\begin{equation}\label{mrd:eq:stage-quantile}
z_{c_k(\alpha_m)}^2\le C(t_{m,k}+A_m).
\end{equation}
Write $x=k/m\le L_m/m\le\sqrt{p_m}$. The deterministic part of
\eqref{mrd:eq:deleted-stage}, after multiplication by
$z_{c_k}/m$, is at most
\[
C_Bx\sqrt{b_m\log(e/x)\{\log(e/x)+A_m\}}
\le C_B\sqrt{b_m}\,x\{\log(e/x)+A_m\}.
\]
The function $x\mapsto x\log(e/x)$ is increasing for $0<x\le1$.
Consequently, uniformly for these $k$,
\begin{equation}\label{mrd:eq:stage-deletion-limit}
\sqrt{b_m}\,x\{\log(e/x)+A_m\}
\le\sqrt{b_m p_m}\{\log(e/\sqrt{p_m})+A_m\}
\le C\sqrt{p_m}\,b_m^{3/2}\to0.
\end{equation}
The concentration remainder, similarly scaled, is at most
\[
C_B\sqrt{\frac{b_m\log m}{m}\,
 x\{\log(e/x)+A_m\}}
\le C_Bb_m\sqrt{\frac{\log m}{m}}\to0,
\]
where $x\log(e/x)\le1$, $xA_m\le b_m$, and $b_m\ge1$ were used.

Gaussian mixture tails and a union bound further give
$\max_i|Y_i|\le C_B\sqrt{b_m\log m}$ and
$|Z|\le C_B\sqrt{\log m}$ outside probability $C_Bm^{-B}$.
Their contributions in \eqref{mrd:eq:stage-error}, even scaled by the
largest quantile $z_{c_1}$, are at most
$C_B\sqrt{b_m}\log m/m=o(1)$.
For each fixed $a>0$, all terms other than the full sum are therefore
at most $a/2$ eventually, simultaneously over the stages and their sets.
Independence and the sub-Gaussian proxies give
\[
\Pr\left\{\left|\sum_iB_i\right|>cam/z_{c_1}\right\}
\le2\exp\{-c_am/(b_m\log m)\}
\le2\exp\{-c_am/(\log m)^2\}.
\]
Since $z_{c_k}\le z_{c_1}$, this handles the remaining term uniformly.
This proves \eqref{mrd:eq:uniform-tail} for every fixed $B>0$.
Finally, if $\alpha_m\delta_m\to0$, then eventually
$\delta_m\le1/\alpha_m$ and
$(1+\delta_m)/p_m\le C_{b_*}(m+m^2)$. Choose $B>3$ to obtain
\eqref{mrd:eq:uniform-risk}.
\end{proof}

\subsection{Embedding adaptive rejections without matching paths}

\begin{lemma}[Two deterministic transfers]\label{mrd:lem:transfers}
Under the hypotheses of Lemma~\ref{mrd:lem:uniform}, there is a fixed $a_0>0$
such that on $\mathcal E_m(a_0)$, eventually and simultaneously for
$1\le k\le L_m$, $m-|A|\le k-1$ and $i\in A$,
\begin{align}
2\Phibar(|U_{i,A}|)\le c_k(\alpha_m)
&\ \Longrightarrow\ \widetilde P_i\le c_k(2\alpha_m),\label{mrd:eq:transfer-upper}\\
\widetilde P_i\le\tfrac12 c_k(\alpha_m)
&\ \Longrightarrow\ 2\Phibar(|U_{i,A}|)\le c_k(\alpha_m).
\label{mrd:eq:transfer-lower}
\end{align}
Here $2\alpha_m<1$ eventually.
\end{lemma}
\begin{proof}
Put $A_{m,k}=m+1-k$. Directly,
\[
\frac{c_k(2\alpha_m)}{c_k(\alpha_m)}
=\frac{2(A_{m,k}+k\alpha_m)}{A_{m,k}+2k\alpha_m}
\ge\frac32
\]
eventually, uniformly for $k\le L_m=o(m)$. Also
$c_{L_m}(\alpha_m)\to0$. Mills' ratio gives, for fixed $a>0$,
\[
\frac{\Phibar(z_t-a/z_t)}{\Phibar(z_t)}\to e^a,
\qquad
\frac{\Phibar(z_t+a/z_t)}{\Phibar(z_t)}\to e^{-a}
\quad(t\downarrow0).
\]
Choose $a_0>0$ such that $e^{2a_0}<3/2$ and $e^{-2a_0}>1/2$.
These limits hold uniformly for $t=c_k(\alpha_m)$, $k\le L_m$, since
the largest such $t$ tends to zero. They imply
\[
z_{c_k(2\alpha_m)}\le z_{c_k(\alpha_m)}-
\frac{2a_0}{z_{c_k(\alpha_m)}},\qquad
z_{c_k(\alpha_m)/2}\ge z_{c_k(\alpha_m)}+
\frac{2a_0}{z_{c_k(\alpha_m)}}.
\]
On $\mathcal E_m(a_0)$ the perturbation at stage $k$, for every
set with at most $k-1$ deletions, is at most $a_0/z_{c_k}$. The two triangle inequalities for absolute
values give \eqref{mrd:eq:transfer-upper} and \eqref{mrd:eq:transfer-lower}.
\end{proof}

Let $\widetilde R_m',\widetilde V_m'$ be the rejection and false-rejection
counts of ordinary GBS applied to the independent factor-known p-values
$\widetilde P_i$ at level $2\alpha_m$.

\begin{lemma}[Pathwise inclusion up to the deterministic cap]\label{mrd:lem:inclusion}
On $\mathcal E_m(a_0)$, if $R_m^M<L_m$, every adaptive MRD--GBS
rejection belongs to the rejection set of factor-known GBS at $2\alpha_m$.
If $R_m^M\ge L_m$, then $\widetilde R_m'\ge L_m$.
Consequently, without any assumed bound on $R_m^M$,
\begin{equation}\label{mrd:eq:pathwise-V}
V_m^M\le\widetilde V_m'
+mI\{\mathcal E_m(a_0)^c\}+mI\{\widetilde R_m'\ge L_m\}.
\end{equation}
\end{lemma}
\begin{proof}
Suppose the adaptive rule makes at least $r\le L_m$ rejections in the
order $i_1,\ldots,i_r$. At stage $k\le r$ exactly $k-1$ indices have been deleted, so \eqref{mrd:eq:transfer-upper} gives
$\widetilde P_{i_k}\le c_k(2\alpha_m)$. For every $j\le r$, the first
$j$ selected indices therefore all have factor-known p-values at most
$c_j(2\alpha_m)$. Thus
$\widetilde P_{(j)}\le c_j(2\alpha_m)$ for all $j\le r$, and
$\widetilde R_m'\ge r$. This reasoning does not require matching the
selected order with the factor-known order.

If the adaptive rule stops at $r<L_m$, $r=0$ is trivial. For $r>0$, every
selected index has $\widetilde P_i\le c_r(2\alpha_m)
\le c_{\widetilde R_m'}(2\alpha_m)$. Such an index must be rejected by
the static step-down rule. Otherwise, if $\widetilde R_m'<m$, an excluded
p-value at or below $c_{\widetilde R_m'}$ would force
$\widetilde P_{(\widetilde R_m'+1)}\le c_{\widetilde R_m'}
<c_{\widetilde R_m'+1}$, contradicting the stopping definition.
For $\widetilde R_m'=m$ inclusion is automatic. If the adaptive rule
reaches $L_m$, use the first part with $r=L_m$.
The three-event decomposition gives \eqref{mrd:eq:pathwise-V}.
\end{proof}

\begin{lemma}[The cap is reached with sufficiently small probability]\label{mrd:lem:cap-tail}
Under the preceding assumptions, for some $c>0$ and eventually,
\begin{equation}\label{mrd:eq:cap-tail}
\mathbb P(\widetilde R_m'\ge L_m)\le2e^{-cL_m}.
\end{equation}
Consequently
$(1+\delta_m)m\mathbb P(\widetilde R_m'\ge L_m)=o(mp_m)$ when
$\alpha_m\delta_m\to0$.
\end{lemma}
\begin{proof}
Let $t_m=c_{L_m}(2\alpha_m)$. Since $L_m=o(m)$,
$t_m\le4\alpha_m L_m/m$ eventually. If $\widetilde R_m'\ge L_m$,
at least $L_m$ p-values are at most $t_m$. Unless $K_m\ge L_m/2$,
at least $L_m/2$ of these are null p-values. Write
\[
N_{0,m}(t_m)=\sum_i I\{\nu_i=0,\ \widetilde P_i\le t_m\}.
\]
These indicators are independent Bernoulli variables of common parameter
$(1-p_m)t_m$. Hence $\mathbb EN_{0,m}(t_m)\le4\alpha_m L_m$, and
\[
\mathbb P(\widetilde R_m'\ge L_m)
\le\mathbb P(K_m\ge L_m/2)+\mathbb P(N_{0,m}(t_m)\ge L_m/2).
\]
The binomial Chernoff bound $\mathbb P(W\ge x)\le(e\mathbb EW/x)^x$
for $x\ge\mathbb EW$ yields
\[
\mathbb P(\widetilde R_m'\ge L_m)
\le\left(\frac{2e mp_m}{L_m}\right)^{L_m/2}
+(8e\alpha_m)^{L_m/2}\le2e^{-cL_m}.
\]
Both bases tend to zero. The cap grows at least as
$c\sqrt m$, and $(1+\delta_m)/p_m=O(m^2)$ as in
Lemma~\ref{mrd:lem:uniform}; the loss-weighted conclusion follows.
This is the step that avoids a circular assumption that the adaptive
path is already short. A Markov bound using only the expected static
rejection count would generally be too weak here.
\end{proof}

\subsection{False-discovery risk}

The following static fact is already proved in the fixed-rGBS analysis.
The full proof appears in Lemma~\ref{lem:factor-known-gbs-type-one};
its argument allows arbitrary non-null p-value distributions.
\begin{lemma}[Static factor-known false discoveries]\label{mrd:lem:static-V}
For independent factor-known p-values, with independent uniform nulls,
GBS at level $\beta\in(0,1)$ satisfies
\[
\mathbb E\{\widetilde V_m(\beta)\mid K_m=k\}
\le\frac\beta{1-\beta}(k+1),\qquad
\mathbb E\widetilde V_m(\beta)\le\frac\beta{1-\beta}(mp_m+1).
\]
Signal p-values may have different distributions across coordinates.
\end{lemma}

\begin{lemma}[Adaptive Type-I bound]\label{mrd:lem:typeI}
Under the assumptions of Lemma~\ref{mrd:lem:uniform} and
$\alpha_m\delta_m\to0$,
\[
\mathbb EV_m^M=o(mp_m),\qquad
\delta_m\mathbb EV_m^M=o(mp_m).
\]
Moreover $\mathbb P(R_m^M\ge L_m)\to0$.
\end{lemma}
\begin{proof}
Lemma~\ref{mrd:lem:static-V} gives
$\mathbb E\widetilde V_m'\le
2\alpha_m(mp_m+1)/(1-2\alpha_m)=O(\alpha_mmp_m)$, since $mp_m\ge b$.
Take expectations in \eqref{mrd:eq:pathwise-V} and use
Lemmas~\ref{mrd:lem:uniform} and \ref{mrd:lem:cap-tail}. The two exceptional
terms are $o(mp_m)$ both without loss weighting and after multiplication
by $\delta_m$. The static term is negligible because
$\alpha_m\to0$ and $\alpha_m\delta_m\to0$.
More explicitly, writing $F_m=\Pr\{\mathcal E_m(a_0)^c\}$ and
$H_m=\Pr(\widetilde R_m'\ge L_m)$,
\[
\frac{(1+\delta_m)\mathbb EV_m^M}{mp_m}
\le \frac{2\alpha_m(1+\delta_m)}{1-2\alpha_m}
\left(1+\frac1{mp_m}\right)
+\frac{1+\delta_m}{p_m}(F_m+H_m)\longrightarrow0.
\]
On $\mathcal E_m(a_0)$, reaching the adaptive cap implies reaching the
static cap, so the probability assertion follows as well.
\end{proof}

\subsection{Missed signals when the expected count diverges}

\begin{lemma}[Tightened factor-known signal crossing]\label{mrd:lem:crossing}
Under \eqref{eq:framework-loadings}, \eqref{mrd:eq:div-counts}, and the first calibration
in \eqref{mrd:eq:div-calibration}, for each fixed $0<\epsilon<q_-$, put
$r_m=\lfloor(Q_m-\epsilon)mp_m\rfloor$. Then
\begin{equation}\label{mrd:eq:tight-crossing}
\mathbb P\{\widetilde P^{1}_{(j)}\le\tfrac12c_j(\alpha_m),
\quad1\le j\le r_m\}\to1,
\end{equation}
where the superscript 1 denotes ordering only over non-null coordinates.
The event is defined to fail if $K_m<r_m$.
\end{lemma}
\begin{proof}
Let
\[
G_{i,m}(t)=2\Phibar\left(\frac{z_t}{\sqrt{1+\psi_m^2/\sigma_i^2}}\right),
\qquad \overline G_m(t)=m^{-1}\sum_iG_{i,m}(t).
\]
For any fixed $a\in(0,1)$, replacing $\alpha_m$ by $a\alpha_m$ leaves
the calibration constant $C_0$ unchanged. At
$t_m=(Q_m-\epsilon)a\alpha_mp_m$, uniformly in $i$,
\[
\frac{z_{t_m}}{\sqrt{1+\psi_m^2/\sigma_i^2}}
=\sqrt{C_0}\sigma_i+o(1),\qquad
\overline G_m(t_m)=Q_m+o(1).
\]
Here $Q_m-\epsilon$ is bounded above and away from zero; Gaussian
quantile asymptotics and the loading bounds justify uniformity.
The normal hazard $h(x)=\phi(x)/\Phibar(x)$ is increasing on $[0,\infty)$:
$h'(x)=h(x)\{h(x)-x\}>0$ by Mills' inequality. For $S>1$,
\[
\frac{d}{dz}\log\frac{\Phibar(z/S)}{\Phibar(z)}
=h(z)-S^{-1}h(z/S)>0.
\]
Thus each $G_{i,m}(t)/t$, and its average, is non-increasing in $t$.
There is a fixed $\eta>0$ such that, eventually,
\[
\overline G_m(q a\alpha_mp_m)\ge(1+\eta)q,
\qquad0<q\le Q_m-\epsilon.
\]
Indeed, with $q_m=Q_m-\epsilon$, eventually
$\overline G_m(q_ma\alpha_mp_m)\ge Q_m-\epsilon/2$; hence its ratio to
$q_m$ is at least $1+\epsilon/(2q_+)$. Choose any fixed
$0<\eta<\epsilon/(2q_+)$.

Put $d_j=c_j(\alpha_m)/2$. For $j\ge2$ eventually
$d_j\ge j\alpha_m/(2m)$. Use the preceding bound with $a=1/2$ and
$q=j/(mp_m)$ to obtain
\[
\mathbb EN_{1,m}(d_j)=mp_m\overline G_m(d_j)\ge(1+\eta)j,
\quad 2\le j\le r_m,
\]
where $N_{1,m}(t)=\sum_i I\{\nu_i=1,\widetilde P_i\le t\}$ is a sum
of independent, possibly nonidentical Bernoulli variables.

The earliest ranks require a separate bound. Write
$A_m=\log(1/\alpha_m)$, $B_m=\log(1/p_m)$, $H_m=\log(mp_m)$.
Then $\log(1/d_1)=A_m+B_m+H_m+O(1)$ and
$\psi_m^2\sim2(A_m+B_m)/C_0$.
As $\psi_m^2/\sigma_i^2\ge\psi_m^2$, set
$y_m=z_{d_1}/\sqrt{1+\psi_m^2}$ to obtain
\[
\frac{y_m^2}{H_m}\le C\left(\frac1{H_m}+\frac1{A_m+B_m}\right)\to0,
\]
using $z_t^2\le2\log(1/t)$. Mills' lower bound gives
\[
M_m:=mp_m\overline G_m(d_1)
\ge2e^{H_m}\Phibar(y_m)\to\infty.
\]
Choose integers $\ell_m\to\infty$ with
$\ell_m=o(M_m)$ and $\ell_m=o(mp_m)$. A Chernoff bound shows
$\mathbb P\{N_{1,m}(d_1)<\ell_m\}\to0$, which handles all
$j\le\ell_m$ by monotonicity of $d_j$.
For the remaining ranks the mean bound and the multiplicative Chernoff
inequality give
\[
\sum_{j=\ell_m+1}^{r_m}\mathbb P\{N_{1,m}(d_j)<j\}
\le\sum_{j=\ell_m+1}^{\infty}e^{-c_\eta j}\to0.
\]
Together these bounds yield \eqref{mrd:eq:tight-crossing}.
\end{proof}

\begin{lemma}[Crossing along the adaptive path]\label{mrd:lem:adaptive-crossing}
Under Theorem~\ref{mrd:thm:diverging}'s assumptions, for $r_m$ as above,
$\mathbb P(R_m^M\ge r_m)\to1$, and
\begin{equation}\label{mrd:eq:adaptive-T}
\mathbb ET_m^M\le mp_m\{1-Q_m+\epsilon\}+o(mp_m).
\end{equation}
\end{lemma}
\begin{proof}
By Lemma~\ref{mrd:lem:cap}, $r_m<L_m$ eventually. On the event
\eqref{mrd:eq:tight-crossing} intersected with $\mathcal E_m(a_0)$,
suppose the adaptive rule were to stop at stage $k\le r_m$.
Only $k-1$ coordinates have been deleted, so at least one of the first
$k$ factor-known signal order statistics remains active. Its factor-known
p-value is at most $c_k(\alpha_m)/2$. By \eqref{mrd:eq:transfer-lower},
its current MRD residual passes the stage-$k$ threshold. The maximum
absolute residual therefore passes as well, contradicting stopping.
Consequently the adaptive rejection count is at least $r_m$ on this
intersection, whose probability tends to one.

It is not necessary to assert that the first $r_m$ signal indices
are themselves rejected. Instead, use the exact identity
$R_m^M=S_m^M+V_m^M$, where $S_m^M$ counts true discoveries. Pathwise,
\[
T_m^M=K_m-R_m^M+V_m^M
\le K_m-r_m+V_m^M+r_m I\{R_m^M<r_m\}.
\]
Taking expectations, invoking Lemma~\ref{mrd:lem:typeI}, and using
$r_m=(Q_m-\epsilon)mp_m+O(1)$ proves \eqref{mrd:eq:adaptive-T}.
This directly addresses the possibility that some early adaptive
rejections are nulls.
\end{proof}

\subsection{Oracle normalization and proofs of the main results}

\begin{lemma}[Factor-known Oracle leading risk]\label{mrd:lem:oracle}
Under either Theorem~\ref{mrd:thm:diverging}'s or
Theorem~\ref{mrd:thm:bounded}'s assumptions,
\begin{equation}\label{mrd:eq:oracle-leading}
\mathcal R_m^{\mathrm{BO},Z}=mp_m\delta_{A,m}\{1-Q_m+o(1)\}.
\end{equation}
There are constants $0<q_-<q_+<1$ such that
$q_-\le Q_m\le q_+$ eventually.
\end{lemma}
\begin{proof}
Put $u_i=\psi_m^2/\sigma_i^2$. The factor-known Oracle rejects
$|Y_i|>c_i^O$, where
\[
(c_i^O)^2=(1+u_i^{-1})\log\{(1+u_i)f_m^2\delta_m^2\}.
\]
For diverging counts, compatibility gives
\[
\frac{\log(f_m\delta_m)}{\psi_m^2}
=\frac{\log(f_m/r_{\alpha_m})+\log(r_{\alpha_m}\delta_m)}{\psi_m^2}
\to C_0/2.
\]
For bounded counts, the same limit follows from
$\log f_m=\log m+O(1)$ and fixed positive loss ratio.
In both cases $\psi_m^2\to\infty$ and, uniformly over coordinates,
\[
\frac{(c_i^O)^2}{1+u_i}=C_0\sigma_i^2+o(1),\qquad
c_i^O\asymp\sqrt{u_i}.
\]
Writing $L_i=\log\{(1+u_i)f_m^2\delta_m^2\}$,
\[
f_m\delta_m\exp\{-(c_i^O)^2/2\}
=(1+u_i)^{-1/2}\exp\{-L_i/(2u_i)\}=O(u_i^{-1/2})
\]
uniformly. Mills' upper bound yields
$f_m\delta_m2\Phibar(c_i^O)=O(u_i^{-1})=o(1)$ uniformly.
The alternative miss probability is
$2\Phi(c_i^O/\sqrt{1+u_i})-1
=2\Phi(\sqrt{C_0}\sigma_i)-1+o(1)$ uniformly.
Summing these two error contributions proves \eqref{mrd:eq:oracle-leading}.
Finally $\sqrt{1-\bar\lambda^2}\le\sigma_i\le1$, so the averages $Q_m$
stay bounded away from both zero and one.
\end{proof}

\subsubsection*{Proof of Theorem~\ref{mrd:thm:diverging}}
\begin{proof}
The calibration, $mp_m\to\infty$, and $\alpha_m/p_m\to\infty$ imply
$\psi_m^2=O\{\log(e/p_m)\}$, as
$\log(f_m/r_{\alpha_m})=O\{\log(e/p_m)\}$ and its ratio to $\psi_m^2$ has a
positive finite limit. Thus all hypotheses of Lemma~\ref{mrd:lem:uniform}
are satisfied. Lemma~\ref{mrd:lem:typeI} proves \eqref{mrd:eq:main-typeI}.
For every fixed $0<\epsilon<q_-$,
Lemma~\ref{mrd:lem:adaptive-crossing} gives
\[
\frac{\mathcal R_m^M}{mp_m\delta_{A,m}}
\le1-Q_m+\epsilon+o(1).
\]
By \eqref{mrd:eq:information-order} and Lemma~\ref{mrd:lem:oracle},
\[
\mathcal R_m^M\ge\mathcal R_m^{\mathrm{BO},\boldsymbol{X}}\ge
mp_m\delta_{A,m}\{1-Q_m+o(1)\}.
\]
Although $Q_m$ need not converge, $1-Q_m\ge1-q_+>0$. Consequently
\[
1\le\liminf_m\frac{\mathcal R_m^M}{\mathcal R_m^{\mathrm{BO},\boldsymbol{X}}}
\le\limsup_m\frac{\mathcal R_m^M}{\mathcal R_m^{\mathrm{BO},\boldsymbol{X}}}
\le1+\frac{\epsilon}{1-q_+}.
\]
Let $\epsilon\downarrow0$. The same inequalities show that the difference
between $\mathcal R_m^M/(mp_m\delta_{A,m})$ and $1-Q_m$ tends to zero.
Subtracting the negligible weighted Type-I contribution gives
\eqref{mrd:eq:main-typeII}. This also recovers the genuine observed-data Oracle
expansion without assuming equality of the two information sets.
\end{proof}

\subsubsection*{Proof of Theorem~\ref{mrd:thm:bounded}}
\begin{proof}
Here $p_m\asymp1/m$, so $mp_m\ge b_*>0$ eventually,
$\alpha_m/p_m=m^{1-o(1)}\to\infty$, and $\psi_m^2=O\{\log(e/p_m)\}$.
Lemmas~\ref{mrd:lem:uniform}--\ref{mrd:lem:typeI} therefore apply.
Set $c_{1,m}=\alpha_m/(m+\alpha_m)$ and call a signal strong if
\[
\widetilde P_i\le c_{1,m}/2.
\]
On $\mathcal E_m(a_0)\cap\{\widetilde R_m'<L_m\}$, the adaptive rule
stops after $R_m^M<L_m$ rejections, at stage $R_m^M+1\le L_m$.
If a strong signal remained
active at its stopping stage $k\le L_m$, then
$\widetilde P_i\le c_{1,m}/2\le c_k(\alpha_m)/2$.
Lemma~\ref{mrd:lem:transfers} would force its current residual, hence the
maximum residual, to pass the threshold. This contradicts stopping.
Therefore every strong signal is rejected on this event.

Let $N_m^{\mathrm{strong}}$ count such signals and let $\mathcal B_m$
be the complement of the preceding event. Then, pathwise,
\[
T_m^M\le K_m-N_m^{\mathrm{strong}}+K_m I(\mathcal B_m).
\]
The exceptional probability tends to zero by
Lemmas~\ref{mrd:lem:uniform} and \ref{mrd:lem:cap-tail}. Since
$\mathbb EK_m^2=mp_m(1-p_m)+(mp_m)^2=O((mp_m)^2)$,
Cauchy--Schwarz gives
$\mathbb E\{K_m I(\mathcal B_m)\}=o(mp_m)$.
Further,
\[
z_{c_{1,m}/2}^2=2\log m\{1+o(1)\},\qquad
\frac{z_{c_{1,m}/2}}{\sqrt{1+\psi_m^2/\sigma_i^2}}
=\sqrt{C_0}\sigma_i+o(1)
\]
uniformly, so
$\mathbb EN_m^{\mathrm{strong}}=mp_m\{Q_m+o(1)\}$.
Thus $\mathbb ET_m^M\le mp_m\{1-Q_m+o(1)\}$.
The negligible Type-I risk and the same Oracle sandwich as above prove
all conclusions. No concentration of $K_m$ at a diverging deterministic
quantity has been used.
\end{proof}

\subsubsection*{Proofs of Corollaries~\ref{mrd:cor:polynomial} and \ref{mrd:cor:homogeneous}}
\begin{proof}
For $0<\beta<1$, $mp_m\to\infty$ and
$\log f_m=\beta\log m+O(1)$. The sub-polynomial decrease of the level
implies $\alpha_mf_m=m^{\beta-o(1)}\to\infty$,
$2\log(f_m/r_{\alpha_m})/\psi_m^2\to C_0$,
$\log\{1/(r_{\alpha_m}\delta)\}/\psi_m^2\to0$, and
$\alpha_m\delta\to0$. No additional sparsity rate needs to be checked.
For $\beta=1$, the two-sided bounds on $mp_m$ and the stated variance
scaling give Theorem~\ref{mrd:thm:bounded}. Finally constant loadings satisfy
\eqref{eq:framework-loadings}; replacing $\psi_m^2$ by
$(1-\lambda^2)u_{m,\lambda}$ gives $C_\lambda=(1-\lambda^2)C_0$ and the
homogeneous risk formula.
\end{proof}

\section{Oracle-parameterized BSD: ABOS and decision-set agreement}
\label{app:bsd-abos}

The result below concerns the BSD score in the one-factor model of
Section~\ref{sec:framework}. The known-covariance representation is
\eqref{eq:heterogeneous-factor-model}; the loading conditions are
\eqref{eq:hetero-loading-upper}--\eqref{eq:hetero-loading-information}.
Write $\delta_m=\delta_{0,m}/\delta_{A,m}$, consistently with
Section~\ref{sec:framework}. BSD uses these true model and loss inputs.
Let $b_{i,m}$ denote the factor-known posterior-odds threshold, so that
\begin{equation}
\label{bsd:eq:factor-threshold}
b_{i,m}^2=\left(1+\frac{\sigma_{i,m}^2}{\psi_m^2}\right)
\left[2\log\!\left\{\Delta_m\frac{1-p_m}{p_m}\right\}
+\log\!\left(1+\frac{\psi_m^2}{\sigma_{i,m}^2}\right)\right].
\end{equation}
For later use, define $Y_i=(X_i-\lambda_{i,m}Z)/\sigma_{i,m}$,
$B_i=(\lambda_{i,m}/\sigma_{i,m})Y_i$,
$T_{A,-i}=\sum_{j\in A\setminus\{i\}}\lambda_{j,m}^2/\sigma_{j,m}^2$,
$N_{A,-i}=\sum_{j\in A\setminus\{i\}}B_j$, and
$h_{i,A}=\{1+(\lambda_{i,m}^2/\sigma_{i,m}^2)/(1+T_{A,-i})\}^{1/2}$.
The exact active-set residual identity is
\begin{equation}
\label{bsd:eq:residual-identity}
U_{i,A}=h_{i,A}^{-1}\left[Y_i+
\frac{\lambda_{i,m}}{\sigma_{i,m}(1+T_{A,-i})}
\{Z-N_{A,-i}\}\right].
\end{equation}

\begin{proof}
The proof has four parts: identify the factor-known likelihood-ratio
threshold, control all active-set residuals before a deterministic cap,
transfer threshold crossings to the adaptive BSD path, and show that the
cap is reached with negligible probability.

\medskip
\noindent\emph{1. Factor-known threshold and its risk.}
If \(Z\) were observed, the posterior odds for coordinate \(i\) would be
\[
O_{i,m}^Z=
\frac{p_m}{1-p_m}
\sqrt{\frac{\sigma_{i,m}^2}{\psi_m^2+\sigma_{i,m}^2}}
\exp\left\{\frac{\psi_m^2}
{2(\psi_m^2+\sigma_{i,m}^2)}Y_i^2\right\}.
\]
Thus \(O_{i,m}^Z>\Delta_m\) is equivalent to
\(|Y_i|>b_{i,m}\), with \(b_{i,m}\) given by \eqref{bsd:eq:factor-threshold}. Uniformly in \(i\),
\[
\frac{b_{i,m}^2}
{1+\psi_m^2/\sigma_{i,m}^2}
=
\frac{\sigma_{i,m}^2}{\psi_m^2}
\left[
2\log\left\{\Delta_m\frac{1-p_m}{p_m}\right\}
+\log\left(1+\frac{\psi_m^2}{\sigma_{i,m}^2}\right)
\right]
=C_0\sigma_{i,m}^2+o(1).
\tag{BSD.4}
\]
The error is uniform because the loadings are bounded away from unit
magnitude and \(\Delta_m\) is bounded above and away from zero. Therefore
\[
\mathbb P_{\nu_i=1}(|Y_i|\le b_{i,m})
=2\Phi(\sqrt{C_0}\sigma_{i,m})-1+o(1)
\tag{BSD.5}
\]
uniformly in \(i\).

Under a null, Mills' inequality and \eqref{bsd:eq:factor-threshold}--(BSD.4) give, uniformly in \(i\),
\[
\mathbb P_{\nu_i=0}(|Y_i|>b_{i,m})
\le C\frac{p_m}{\sqrt{\log(1/p_m)\,\log(1/p_m)}}=
O\!\left(\frac{p_m}{\log(1/p_m)}\right).
\tag{BSD.6}
\]
(The constant may depend on \(C_0,\bar\lambda,\Delta\).) Thus the factor-known Bayes Oracle has risk
\[
\mathcal R_m^{\mathrm{BO},Z}
=mp_m \delta_{A,m}\{1-Q_m+o(1)\}.
\]
Its Type-I risk is \(o(mp_m\delta_{A,m})\), while its Type-II risk is
\(mp_m\delta_{A,m}\{1-Q_m+o(1)\}\).

\medskip
\noindent\emph{2. Uniform residual control over a capped set family.}
Let
\[
L_m=\lfloor m\sqrt{p_m}\rfloor,\qquad
q_m=\max_{i\le m}b_{i,m},\qquad
\mathcal E_m(a)=
\left\{q_m\sup_{\substack{A\subseteq[m],\,m-|A|\le L_m\\i\in A}}
|U_{i,A}-Y_i|\le a\right\}.
\]
The loading and count assumptions imply \(L_m=o(m)\), \(L_m\to\infty\),
\(L_m/(mp_m)\to\infty\), \(q_m^2\asymp\log(1/p_m)\), and
\(\log(1/p_m)=O(\log m)\). From \eqref{bsd:eq:residual-identity}, uniformly over the displayed sets,
\[
|U_{i,A}-Y_i|
\le \frac C m\left(
|Z|+\left|\sum_{j=1}^m B_j\right|
+\max_j|Y_j|+H_{m,L_m+1}\right),
\quad
H_{m,L}=\max_{|J|\le L}\sum_{j\in J}|B_j|.
\tag{BSD.7}
\]
The independent centered variables \(B_i\) and \(Y_i\) have sub-Gaussian
variance proxies bounded by \(C\log(e/p_m)\). A union bound over subsets
and signs yields, for every fixed \(B>0\), with probability at least
\(1-O(m^{-B})\),
\[
H_{m,L_m+1}\le C_B\sqrt{\log(e/p_m)}
\left[
(L_m+1)\sqrt{\log\{em/(L_m+1)\}}+\sqrt{(L_m+1)\log m}
\right],
\tag{BSD.8}
\]
and simultaneously
\[
\max_i|Y_i|\le C_B\sqrt{\log(e/p_m)\log m},\quad
\left|\sum_iB_i\right|\le C_B\sqrt{m\log(e/p_m)\log m},\quad
|Z|\le C_B\sqrt{\log m}.
\]
Multiplying (BSD.7) by \(q_m=O(\sqrt{\log(e/p_m)})\), the leading subset
term is bounded by
\[
C_B\sqrt{p_m}\{\log(e/p_m)\}^{3/2}=o(1).
\]
The remaining terms also tend to zero, using
\(p_m\ge c/m\) (from \(\liminf mp_m>0\)) and
\(\log(1/p_m)=O(\log m)\). For the sum term, a sub-Gaussian tail bound
at level \(a m/q_m\) gives an additional probability at most
\(2\exp\{-c_a m/\log^2m\}\). Consequently, for every fixed \(a>0\) and
every \(B>0\),
\[
\mathbb P\{\mathcal E_m(a)^c\}
\le C_{a,B}m^{-B}+2\exp\{-c_a m/\log^2m\}.
\tag{BSD.9}
\]
In particular, \(\mathbb P\{\mathcal E_m(a)^c\}=o(p_m)\).

Also, for \(m-|A|\le L_m\), the loading condition implies
\[
\sup_{i\in A}|v_{i,A}-\sigma_{i,m}^2|=O(m^{-1}).
\tag{BSD.10}
\]
If \(b_{i,A}\) denotes the exact BSD residual threshold obtained by
replacing \(\sigma_{i,m}^2\) by \(v_{i,A}\) in \eqref{bsd:eq:factor-threshold}, differentiation of
the right side of \eqref{bsd:eq:factor-threshold} with respect to \(v\), together with
\(\psi_m^2\asymp\log(1/p_m)\), gives
\[
\sup_{m-|A|\le L_m,\,i\in A}|b_{i,A}^2-b_{i,m}^2|=O(m^{-1}),
\qquad
q_m|b_{i,A}-b_{i,m}|=o(1).
\tag{BSD.11}
\]
Hence for each fixed \(a>0\), on \(\mathcal E_m(a)\) there is a constant
\(d_a<\infty\), independent of \(m,A,i\), such that
\[
S_{i,A}>\Delta_m\ \Longrightarrow\
|Y_i|>b_{i,m}-d_a/b_{i,m},
\tag{BSD.12}
\]
and
\[
|Y_i|>b_{i,m}+d_a/b_{i,m}\ \Longrightarrow\
S_{i,A}>\Delta_m.
\tag{BSD.13}
\]
These implications hold simultaneously for every active set within the cap.

\medskip
\noindent\emph{3. The adaptive path and the cap.}
Define
\[
K_m=\sum_i\nu_i,\qquad
W_m=\sum_i(1-\nu_i)
\mathbf 1\{|Y_i|>b_{i,m}-d_a/b_{i,m}\}.
\]
The pairs \((\nu_i,Y_i)\) are independent across \(i\). By the same Mills
bound as in (BSD.6),
\[
\mathbb E W_m=O\!\left(\frac{mp_m}{\log(1/p_m)}\right)=o(mp_m),
\qquad \mathbb E K_m=mp_m.
\tag{BSD.14}
\]
If BSD reaches \(L_m\) rejections while \(\mathcal E_m(a)\) holds, each of
the first \(L_m\) rejections is either a signal or is counted by \(W_m\).
Thus \(K_m+W_m\ge L_m\). Since \(L_m/(mp_m)\to\infty\), a Chernoff bound
for the independent bounded summands defining \(K_m+W_m\) gives
\[
\mathbb P\{K_m+W_m\ge L_m\}
\le \left(\frac{e\,\mathbb E(K_m+W_m)}{L_m}\right)^{L_m}
\le \exp\{-cL_m\log(1/p_m)\}.
\tag{BSD.15}
\]
Together with (BSD.9), this shows that the probability of reaching the cap,
plus the uniform-residual exceptional probability, is \(o(p_m)\).

On the event that BSD stops before the cap, every missed signal \(i\)
satisfies \(S_{i,A}\le\Delta_m\) at the terminal active set. By (BSD.13),
\[
T_m^{\mathrm{BSD}}
\le \sum_i\nu_i\mathbf 1\{|Y_i|\le b_{i,m}+d_a/b_{i,m}\}
+m\mathbf 1\{\mathcal E_m(a)^c\}
+m\mathbf 1\{\text{BSD reaches }L_m\}.
\tag{BSD.16}
\]
Likewise, each false rejection made before the cap has
\(|Y_i|>b_{i,m}-d_a/b_{i,m}\), so
\[
V_m^{\mathrm{BSD}}
\le W_m+
m\mathbf 1\{\mathcal E_m(a)^c\}
+m\mathbf 1\{\text{BSD reaches }L_m\}.
\tag{BSD.17}
\]
The adaptive ordering need not agree with the ordering of the \(Y_i\)'s:
(BSD.16) uses only the stopping rule, and (BSD.17) uses only the fact that a
rejected candidate has stagewise odds above threshold.

\medskip
\noindent\emph{4. Risk expansion and comparison with the Oracle.}
For each \(i\), the \(O(1/b_{i,m})\) margin in (BSD.16) vanishes after
division by the signal standard deviation, so, uniformly in \(i\),
\[
\mathbb P_{\nu_i=1}(|Y_i|\le b_{i,m}+d_a/b_{i,m})
=2\Phi(\sqrt{C_0}\sigma_{i,m})-1+o(1).
\]
The margin in (BSD.17) changes the null tail in (BSD.6) by at most a fixed
multiplicative factor, so its loss-weighted expectation remains
\(o(mp_m \delta_{A,m})\). Equations (BSD.9), (BSD.14)--(BSD.17) now imply
\[
\mathbb E T_m^{\mathrm{BSD}}=mp_m\{1-Q_m+o(1)\},\qquad
\Delta_m\mathbb E V_m^{\mathrm{BSD}}=o(mp_m).
\]
This proves the displayed BSD risk expansion. Finally, observing \(Z\) gives at least as much information as observing
\(\boldsymbol{X}\) alone, so the Bayes risks obey
\[
\mathcal R_m^{\mathrm{BO},Z}
\le \mathcal R_m^{\mathrm{BO},\boldsymbol{X}}
\le \mathcal R_m^{\mathrm{BSD}}.
\]
The first and last risks have the same expansion, as proved above and in
Steps 1--3. Hence the observed-data Oracle has that expansion as well and
\(\mathcal R_m^{\mathrm{BSD}}/\mathcal R_m^{\mathrm{BO},\boldsymbol{X}}\to1\). The ratio
is well-defined asymptotically because
\[
1-Q_m\ge 2\Phi\!\left(\sqrt{C_0(1-\bar\lambda^2)}\right)-1>0.
\]
\end{proof}

\begin{proof}[Proof of Corollary~\ref{bsd:cor:decision-agreement}]
Write $\ell_m=\log(1/p_m)$ and let
\[
\widehat{\mathcal S}_m^Z
=\{i:|Y_i(Z)|>b_{i,m}\}
\]
be the factor-known Bayes decision set from the proof of
Theorem~\ref{bsd:thm:abos}. We compare each procedure with this set.

First, on the uniform residual event and before BSD reaches the cap, the
crossing implications (BSD.12)--(BSD.13) show that
any disagreement between BSD and $\widehat{\mathcal S}_m^Z$ is confined to
\[
\left\{i:\bigl||Y_i(Z)|-b_{i,m}\bigr|\le d_a/b_{i,m}\right\}.
\]
The uniform Gaussian density bounds at the thresholds give conditional
boundary probabilities $O(p_m/\ell_m)$ under a null and
$O(1/\ell_m)$ under a signal. The exceptional-event and cap probabilities
in (BSD.9) and (BSD.15) are $o(p_m)$. Hence
\begin{equation}
\label{bsd:eq:set-distance-factor-known}
\mathbb E\bigl|\widehat{\mathcal S}_m^{\mathrm{BSD}}
\triangle\widehat{\mathcal S}_m^Z\bigr|
=O(mp_m/\ell_m)+o(mp_m)=o(mp_m).
\end{equation}

Next let
\[
\pi_{i,m}=\mathbb P(\nu_i=1\mid\boldsymbol X,Z),\qquad
q_{i,m}=\mathbb P(\nu_i=1\mid\boldsymbol X),\qquad
t_m=\frac{\delta_m}{1+\delta_m}.
\]
Conditional independence given $Z$ implies
$q_{i,m}=\mathbb E(\pi_{i,m}\mid\boldsymbol X)$, and the factor-known and
observed-data Oracle decisions are $I(\pi_{i,m}>t_m)$ and
$I(q_{i,m}>t_m)$, respectively. Put
\[
D_m=\sum_{i=1}^m\frac{\lambda_{i,m}^2}{\sigma_{i,m}^2},\qquad
\widetilde Z_m=D_m^{-1}\sum_{i=1}^m
\frac{\lambda_{i,m}X_i}{\sigma_{i,m}^2}.
\]
The loading assumptions imply $D_m\asymp m$, and the model gives
\[
\mathbb E(\widetilde Z_m-Z)^2
=D_m^{-2}\left\{D_m+p_m\psi_m^2
\sum_{i=1}^m\frac{\lambda_{i,m}^2}{\sigma_{i,m}^4}\right\}
=O(m^{-1}),
\]
since $p_m\psi_m^2=O(p_m\ell_m)=o(1)$. The posterior mean
$\overline Z_m=\mathbb E(Z\mid\boldsymbol X)$ minimizes mean squared
error, so
\begin{equation}
\label{bsd:eq:posterior-factor-variance}
\mathbb E\operatorname{Var}(Z\mid\boldsymbol X)=O(m^{-1}).
\end{equation}

The conditional inclusion probability is a logistic function of
$Y_i(z)^2$, where $Y_i(z)=(X_i-\lambda_{i,m}z)/\sigma_{i,m}$. Since
$b_{i,m}\asymp\sqrt{\ell_m}$ uniformly and the loadings and residual
standard deviations are uniformly bounded, its derivative in $z$ is
uniformly bounded by $C\sqrt{\ell_m}$. The conditional variance inequality
for Lipschitz functions and \eqref{bsd:eq:posterior-factor-variance}
therefore yield
\begin{equation}
\label{bsd:eq:posterior-inclusion-variance}
\sum_{i=1}^m\mathbb E(q_{i,m}-\pi_{i,m})^2
=\sum_{i=1}^m\mathbb E\operatorname{Var}(\pi_{i,m}\mid\boldsymbol X)
=O(\ell_m).
\end{equation}

Because $t_m$ stays bounded away from zero and one, the Gaussian density
bounds at $b_{i,m}$ also imply, uniformly for sufficiently small
$\varepsilon>0$,
\begin{equation}
\label{bsd:eq:oracle-margin-bound}
\sum_{i=1}^m\mathbb P\bigl(|\pi_{i,m}-t_m|\le\varepsilon\bigr)
\le Cmp_m\varepsilon/\ell_m.
\end{equation}
Indeed, this event restricts $|Y_i(Z)|$ to a band of width
$O(\varepsilon/\sqrt{\ell_m})$ about $b_{i,m}$; the null density there is
$O(p_m/\sqrt{\ell_m})$, and the signal density is
$O(1/\sqrt{\ell_m})$ before multiplication by its prior probability.
For any $\varepsilon>0$,
\[
I\{I(q_{i,m}>t_m)\ne I(\pi_{i,m}>t_m)\}
\le I\{|\pi_{i,m}-t_m|\le\varepsilon\}
+\varepsilon^{-2}(q_{i,m}-\pi_{i,m})^2.
\]
Combining \eqref{bsd:eq:posterior-inclusion-variance} and
\eqref{bsd:eq:oracle-margin-bound},
\[
\mathbb E\bigl|\widehat{\mathcal S}_m^{\mathrm{BO},\boldsymbol X}
\triangle\widehat{\mathcal S}_m^Z\bigr|
\le Cmp_m\varepsilon/\ell_m+C\ell_m/\varepsilon^2.
\]
Take $\varepsilon_m=(\ell_m/(mp_m))^{1/4}$. The condition
$mp_m/\ell_m\to\infty$ makes $\varepsilon_m\to0$, and the last display is
$o(mp_m)$. Together with \eqref{bsd:eq:set-distance-factor-known} and the
triangle inequality, this proves that the expected symmetric difference
between the BSD and observed-data Oracle rejection sets is $o(mp_m)$.
Markov's inequality completes the proof.
\end{proof}

\section{Fixed residualized GBS under homogeneous dependence}
\label{app:rgbs-homogeneous}

\subsection{Procedure and ABOS objective}

For the observed-data procedure, define the leave-one-out factor estimator
\[
\widehat Z_{-i}
=
\frac{1}{\lambda(m-1)}\sum_{j\ne i}X_j,
\]
and the corresponding residual statistic and two-sided $p$-value
\[
\widehat Y_i
=
\frac{X_i-\lambda\widehat Z_{-i}}{a_\lambda},
\qquad
\widehat P_i=2\bar\Phi(|\widehat Y_i|).
\]
Let
\[
\widehat P_{(1)}\le\cdots\le\widehat P_{(m)}
\]
denote the ordered residualized $p$-values.  At nominal level $\alpha_m$,
the Gavrilov--Benjamini--Sarkar critical constants are
\begin{equation}
\label{eq:gbs-critical-constants}
 c_j(\alpha_m)
 =
 \frac{j\alpha_m}{m+1-j(1-\alpha_m)},
 \qquad j=1,\ldots,m.
\end{equation}
The residualized GBS (rGBS) procedure is the step-down rule with rejection
count
\begin{equation}
\label{eq:rgbs-rejection-count-definition}
 R_m^{\mathrm{rGBS}}
 =
 \max\left\{
 r:\widehat P_{(j)}\le c_j(\alpha_m)
 \text{ for every }j=1,\ldots,r
 \right\},
\end{equation}
with the convention that the maximum is zero if the first inequality fails.
Thus the GBS calibration itself is unchanged; the distinction from ordinary
GBS is that it is applied to leave-one-out factor-residualized $p$-values,
which are generally dependent.

The following auxiliary results show that, under sparse
one-factor asymptotics and the GBS calibration, this feasible observed-data
procedure has the same first-order Bayes risk as the genuine Bayes Oracle.
For brevity, throughout the intermediate lemmas we write $c_j=c_j(\alpha_m)$.

\begin{lemma}[Residualized signal distribution approximation]
	\label{lem:residual-signal-distribution}
	Assume the conditions of Theorem~\ref{thm:one-factor-oracle-risk}.
	Let
	\[
	\widehat P_i
	=
	2\bar\Phi(|\widehat Y_i|)
	\]
	and
	\[
	G_m(t)
	=
	2\bar\Phi\left(
	\frac{z_t}{\sqrt{1+u_{m,\lambda}}}
	\right),
	\qquad
	z_t=\Phi^{-1}\left(1-\frac t2\right).
	\]
	Define
	\[
	\widehat G_m(t)
	=
	\mathbb P(\widehat P_i\le t\mid \nu_i=1).
	\]
	Then, for every fixed \(A>0\),
	\[
	\sup_{0<t\le A\alpha_m p_m}
	\left|
	\widehat G_m(t)-G_m(t)
	\right|
	\to0.
	\]
\end{lemma}

\begin{proof}
	Let
	\[
	S_m=\sqrt{1+u_{m,\lambda}},
	\qquad
	W_i=\frac{Y_i(Z)}{S_m}.
	\]
	Under \(\nu_i=1\),
	\[
	W_i\sim N(0,1),
	\]
	and \(W_i\) is independent of \(\Delta_{-i}\). Moreover,
	\[
	\widehat Y_i
	=
	Y_i(Z)-\Delta_{-i}
	=
	S_m\left(W_i-D_{m,-i}\right),
	\qquad
	D_{m,-i}
	=
	\frac{\Delta_{-i}}{S_m}.
	\]
	By the factor-estimation condition, as shown in the proof of
	Theorem~\ref{thm:one-factor-oracle-risk},
	\[
	D_{m,-i}=o_{\mathbb P}(1).
	\]
	
	For \(0<t<1\), let
	\[
	z_t=\Phi^{-1}\left(1-\frac t2\right).
	\]
	Then
	\[
	\widehat P_i\le t
	\quad\Longleftrightarrow\quad
	|\widehat Y_i|\ge z_t
	\quad\Longleftrightarrow\quad
	|W_i-D_{m,-i}|
	\ge
	\frac{z_t}{S_m}.
	\]
	Hence
	\[
	\widehat G_m(t)
	=
	\mathbb P\left(
	|W_i-D_{m,-i}|
	\ge
	\frac{z_t}{S_m}
	\right),
	\]
	whereas
	\[
	G_m(t)
	=
	\mathbb P\left(
	|W_i|
	\ge
	\frac{z_t}{S_m}
	\right).
	\]
	
	Let \(F_m\) denote the distribution function of
	\(W_i-D_{m,-i}\). Since \(D_{m,-i}=o_{\mathbb P}(1)\) and
	\(W_i\sim N(0,1)\), Slutsky's theorem gives
	\[
	W_i-D_{m,-i}
	\Rightarrow
	N(0,1).
	\]
	Since the standard normal distribution function \(\Phi\) is continuous,
	Polya's theorem yields
	\[
	\sup_{x\in\mathbb R}|F_m(x)-\Phi(x)|\to0.
	\]
	Therefore, uniformly over \(b\ge0\),
	\begin{align*}
		&\left|
		\mathbb P(|W_i-D_{m,-i}|\ge b)
		-
		\mathbb P(|W_i|\ge b)
		\right|  \\
		&\le
		\left|1-F_m(b)-\{1-\Phi(b)\}\right|
		+
		\left|F_m(-b)-\Phi(-b)\right| \\
		&\le
		2\sup_{x\in\mathbb R}|F_m(x)-\Phi(x)|
		\to0.
	\end{align*}
	Since the preceding bound is uniform over \(b\ge0\), we may take
	\[
	b=\frac{z_t}{S_m}
	\]
	and then take the supremum over \(0<t\le A\alpha_m p_m\). Thus
	\[
	\sup_{0<t\le A\alpha_m p_m}
	\left|
	\widehat G_m(t)-G_m(t)
	\right|
	\to0.
	\]
	This completes the proof.
\end{proof}
	
\begin{lemma}[Residual deterministic signal-crossing bound]
	\label{lem:residual-deterministic-signal-crossing}
	Assume the conditions of Theorem~\ref{thm:one-factor-oracle-risk}.
	Furthermore, assume that \(\alpha_m\to0\) and
	\begin{equation}
		\label{eq:residual-gbs-calibration}
		\frac{2\log(f_m/r_{\alpha_m})}{u_{m,\lambda}}
		\to C_\lambda\in(0,\infty),
		\qquad
		r_{\alpha_m}
		=
		\frac{\alpha_m}{1-\alpha_m}.
	\end{equation}
	Let
	\[
	q_* = 2\bar\Phi(\sqrt{C_\lambda}).
	\]
	Then, for every \(\varepsilon\in(0,q_*)\), there exists
	\(\eta_\varepsilon>0\) such that, for all sufficiently large \(m\),
	\[
	\widehat G_m(q\alpha_m p_m)
	\ge
	(1+\eta_\varepsilon)q
	\]
	uniformly over
	\[
	0<q\le q_*-\varepsilon.
	\]
\end{lemma}

\begin{proof}
	For \(0<t<1\), let
	\[
	z_t=\Phi^{-1}\left(1-\frac t2\right),
	\]
	and recall that
	\[
	G_m(t)
	=
	2\bar\Phi\left(
	\frac{z_t}{\sqrt{1+u_{m,\lambda}}}
	\right)
	\]
	is the signal distribution function in the corresponding factor-known
	residual experiment.
	
	We first compare \(\widehat G_m\) with \(G_m\). Under \(\nu_i=1\),
	\[
	Y_i(Z)\sim N(0,1+u_{m,\lambda}),
	\]
	independently of \(\Delta_{-i}\). Let
	\[
	K_{-i}=\sum_{j\ne i}\nu_j.
	\]
	Conditional on \(K_{-i}=k\), we have
	\[
	\Delta_{-i}\sim N(0,\tau_{m,k}^2),
	\qquad
	\tau_{m,k}^2
	=
	\frac{1}{m-1}
	+
	\frac{k u_{m,\lambda}}{(m-1)^2}.
	\]
	Therefore,
	\[
	\widehat Y_i
	=
	Y_i(Z)-\Delta_{-i}
	\mid K_{-i}=k,\nu_i=1
	\sim
	N(0,1+u_{m,\lambda}+\tau_{m,k}^2).
	\]
	Hence
	\[
	\mathbb P(\widehat P_i\le t\mid K_{-i}=k,\nu_i=1)
	=
	2\bar\Phi\left(
	\frac{z_t}
	{\sqrt{1+u_{m,\lambda}+\tau_{m,k}^2}}
	\right).
	\]
	Since \(\tau_{m,k}^2\ge0\), we have
	\[
	\frac{z_t}
	{\sqrt{1+u_{m,\lambda}+\tau_{m,k}^2}}
	\le
	\frac{z_t}
	{\sqrt{1+u_{m,\lambda}}},
	\]
	and therefore
	\[
	2\bar\Phi\left(
	\frac{z_t}
	{\sqrt{1+u_{m,\lambda}+\tau_{m,k}^2}}
	\right)
	\ge
	2\bar\Phi\left(
	\frac{z_t}
	{\sqrt{1+u_{m,\lambda}}}
	\right)
	=
	G_m(t).
	\]
	Averaging over \(K_{-i}\), we obtain the pointwise bound
	\begin{equation}
		\label{eq:residual-G-dominates-factor-known-G}
		\widehat G_m(t)\ge G_m(t),
		\qquad 0<t<1.
	\end{equation}
	
    Put \(q_0=q_*-\varepsilon\) and \(t_m=q_0\alpha_mp_m\).
    The calibration gives
    \[
    \frac{z_{t_m}^2}{u_{m,\lambda}}\to C_\lambda,
    \qquad G_m(t_m)\to q_*.
    \]
    For completeness, \(G_m(t)/t\) is non-increasing in \(t\).
    Indeed, put \(S=\sqrt{1+u_{m,\lambda}}>1\),
    \(z=z_t\), and \(h(x)=\phi(x)/\bar\Phi(x)\).
    The normal hazard \(h\) is increasing on \([0,\infty)\), since
    \(h'(x)=h(x)\{h(x)-x\}>0\) by Mills' inequality.
    Therefore
    \[
    \frac{d}{dz}\log\frac{\bar\Phi(z/S)}{\bar\Phi(z)}
    =h(z)-S^{-1}h(z/S)>0.
    \]
    As \(z_t\) decreases in \(t\), the asserted monotonicity follows.
    Consequently, for \(0<q\le q_0\),
    \[
    G_m(q\alpha_mp_m)\ge\frac{q}{q_0}G_m(t_m).
    \]
    Since \(q_*/q_0>1\), choose a fixed \(\eta_\varepsilon>0\)
    with \(1+\eta_\varepsilon<q_*/q_0\).  For all sufficiently large
    \(m\), this proves
    \[
    G_m(q\alpha_mp_m)\ge(1+\eta_\varepsilon)q,
    \qquad 0<q\le q_*-\varepsilon.
    \]
	Combining this with
	\eqref{eq:residual-G-dominates-factor-known-G}, we get
	\[
	\widehat G_m(q\alpha_m p_m)
	\ge
	G_m(q\alpha_m p_m)
	\ge
	(1+\eta_\varepsilon)q
	\]
	uniformly over
	\[
	0<q\le q_*-\varepsilon.
	\]
	This completes the proof.
\end{proof}

\begin{lemma}[Residual initial-crossing condition]
	\label{lem:residual-initial-crossing}
	Assume the conditions of Theorem~\ref{thm:one-factor-oracle-risk}.
	Furthermore, assume that \(\alpha_m\to0\) and that the residual GBS
	calibration condition \eqref{eq:residual-gbs-calibration} holds. Let
	\[
	c_1
	=
	\frac{\alpha_m}{m+\alpha_m}
	\]
	denote the first GBS critical constant. Then
	\[
	mp_m\widehat G_m(c_1)\to\infty
	\qquad\text{as }m\to\infty.
	\]
\end{lemma}

\begin{proof}
	By the conditional Gaussian representation used in the proof of
	Lemma~\ref{lem:residual-deterministic-signal-crossing}, under
	\(\nu_i=1\),
	\[
	\widehat Y_i
	\mid K_{-i}=k
	\sim
	N(0,1+u_{m,\lambda}+\tau_{m,k}^2),
	\]
	where
	\[
	\tau_{m,k}^2
	=
	\frac{1}{m-1}
	+
	\frac{k u_{m,\lambda}}{(m-1)^2}.
	\]
	Thus, for every \(0<t<1\),
	\[
	\widehat G_m(t)
	=
	2\mathbb E\left[
	\bar\Phi\left(
	\frac{z_t}{\sqrt{1+u_{m,\lambda}+\tau_{m,K_{-i}}^2}}
	\right)
	\right],
	\]
	where
	\[
	z_t=\Phi^{-1}\left(1-\frac t2\right).
	\]
	Since \(\tau_{m,K_{-i}}^2\ge0\), we have
	\[
	\widehat G_m(t)
	\ge
	2\bar\Phi\left(
	\frac{z_t}{\sqrt{1+u_{m,\lambda}}}
	\right)
	=
	G_m(t).
	\]
	Therefore,
	\[
	mp_m\widehat G_m(c_1)
	\ge
	mp_mG_m(c_1).
	\]
    Write
    \[
    A_m=\log(1/\alpha_m),\quad B_m=\log(1/p_m),\quad
    H_m=\log(mp_m).
    \]
    Then \(A_m+B_m\to\infty\), \(H_m\to\infty\), and
    \[
    \log(1/c_1)=A_m+B_m+H_m+o(1),\qquad
    u_{m,\lambda}=\frac{2}{C_\lambda}(A_m+B_m)\{1+o(1)\}.
    \]
    Set \(y_m=z_{c_1}/\sqrt{1+u_{m,\lambda}}\).
    The Gaussian bound \(z_t^2\le2\log(1/t)\) gives
    \[
    \frac{y_m^2}{H_m}
    \le C\left\{\frac1{H_m}+\frac1{A_m+B_m}\right\}\to0.
    \]
    Mills' lower bound yields
    \[
    mp_mG_m(c_1)=2e^{H_m}\bar\Phi(y_m)
    \ge\frac{c\exp(H_m-y_m^2/2)}{1+y_m}\to\infty.
    \]
	Hence
	\[
	mp_m\widehat G_m(c_1)\to\infty.
	\]
	This proves the lemma.
\end{proof}

\begin{lemma}[Uniform leave-one-out residual perturbation]
	\label{lem:uniform-residual-perturbation}
	Let
	\[
	\Delta_{-i}
	=
	\frac{1}{m-1}\sum_{j\ne i}\varepsilon_j
	+
	\frac{1}{a_\lambda(m-1)}
	\sum_{j\ne i}\theta_j .
	\]
	Suppose that
	\[
	\frac{\log(m/\alpha_m)\{1+p_m u_{m,\lambda}\}}{m}
	\to0
	\]
	and
	\[
	\frac{
		\sqrt{\log(m/\alpha_m)}
		\sqrt{(1+u_{m,\lambda})\log m}
	}{m}
	\to0.
	\]
	Then
	\[
	\sqrt{\log(m/\alpha_m)}
	\max_{1\le i\le m}|\Delta_{-i}|
	=
	o_{\mathbb P}(1).
	\]
	In particular, if \(t_m\ge c_1=\alpha_m/(m+\alpha_m)\), then
	\[
	z_{t_m}\max_{1\le i\le m}|\Delta_{-i}|
	=
	o_{\mathbb P}(1),
	\]
	where
	\[
	z_t=\Phi^{-1}\left(1-\frac t2\right).
	\]
\end{lemma}

\begin{proof}
	Put
	\[
	B_j=\varepsilon_j+\frac{\theta_j}{a_\lambda}.
	\]
	Then
	\[
	\Delta_{-i}
	=
	\frac{1}{m-1}\sum_{j\ne i}B_j
	=
	\frac{1}{m-1}
	\left(
	\sum_{j=1}^m B_j-B_i
	\right).
	\]
	Hence
	\begin{equation}
		\label{eq:uniform-delta-decomposition}
		\max_{1\le i\le m}|\Delta_{-i}|
		\le
		\frac{1}{m-1}
		\left|\sum_{j=1}^m B_j\right|
		+
		\frac{1}{m-1}
		\max_{1\le i\le m}|B_i|.
	\end{equation}
	
	The variables \(B_1,\ldots,B_m\) are independent and centered. Moreover,
	\[
	\operatorname{Var}(B_j)
	=
	\mathbb E(B_j^2)
	=
	1+p_m u_{m,\lambda}.
	\]
	Therefore, by Chebyshev's inequality,
	\[
	\sum_{j=1}^m B_j
	=
	O_{\mathbb P}\left(
	\sqrt{m\{1+p_m u_{m,\lambda}\}}
	\right).
	\]
	Consequently,
	\[
	\sqrt{\log(m/\alpha_m)}
	\frac{1}{m-1}
	\left|\sum_{j=1}^m B_j\right|
	=
	O_{\mathbb P}
	\left[
	\sqrt{
		\frac{\log(m/\alpha_m)\{1+p_m u_{m,\lambda}\}}{m}
	}
	\right]
	=
	o_{\mathbb P}(1).
	\]
	
	Next we control the maximum term. Conditional on \(\nu_i=0\),
	\(B_i\sim N(0,1)\), while conditional on \(\nu_i=1\),
	\(B_i\sim N(0,1+u_{m,\lambda})\). Hence, for every \(x>0\),
	\begin{align*}
		\mathbb P(|B_i|>x)
		&=
		(1-p_m)\mathbb P(|N(0,1)|>x)
		+
		p_m\mathbb P(|N(0,1+u_{m,\lambda})|>x) \\
		&\le
		2\exp\left\{
		-\frac{x^2}{2(1+u_{m,\lambda})}
		\right\}.
	\end{align*}
	By the union bound,
	\[
	\mathbb P\left(
	\max_{1\le i\le m}|B_i|>x
	\right)
	\le
	2m\exp\left\{
	-\frac{x^2}{2(1+u_{m,\lambda})}
	\right\}.
	\]
	Taking
	\[
	x=M\sqrt{(1+u_{m,\lambda})\log m},
	\qquad M>2,
	\]
	gives
	\[
	\max_{1\le i\le m}|B_i|
	=
	O_{\mathbb P}\left(
	\sqrt{(1+u_{m,\lambda})\log m}
	\right).
	\]
	Thus
	\[
	\sqrt{\log(m/\alpha_m)}
	\frac{1}{m-1}
	\max_{1\le i\le m}|B_i|
	=
	O_{\mathbb P}
	\left[
	\frac{
		\sqrt{\log(m/\alpha_m)}
		\sqrt{(1+u_{m,\lambda})\log m}
	}{m}
	\right]
	=
	o_{\mathbb P}(1).
	\]
	Combining this with \eqref{eq:uniform-delta-decomposition} proves
	\[
	\sqrt{\log(m/\alpha_m)}
	\max_{1\le i\le m}|\Delta_{-i}|
	=
	o_{\mathbb P}(1).
	\]
	
	Finally, if \(t_m\ge c_1\), then \(z_{t_m}\le z_{c_1}\). Since
	\[
	c_1=\frac{\alpha_m}{m+\alpha_m},
	\]
	we have
	\[
	z_{c_1}^2
	=
	O\{\log(m/\alpha_m)\}.
	\]
	Therefore
	\[
	z_{t_m}
	\max_{1\le i\le m}|\Delta_{-i}|
	\le
	z_{c_1}
	\max_{1\le i\le m}|\Delta_{-i}|
	=
	o_{\mathbb P}(1).
	\]
	This completes the proof.
\end{proof}

\begin{lemma}[Tightened factor-known empirical signal crossing]
\label{lem:tightened-factor-known-crossing}
Assume
\[
p_m\to0,\qquad mp_m\to\infty,\qquad
\alpha_m\to0,\qquad \alpha_m f_m\to\infty,
\]
and
\[
\frac{2\log(f_m/r_{\alpha_m})}{u_{m,\lambda}}
\to C_\lambda\in(0,\infty).
\]
Let
\[
q_*=2\bar\Phi(\sqrt{C_\lambda}),
\qquad
r_m=\lfloor(q_*-\varepsilon)mp_m\rfloor,
\]
for a fixed \(\varepsilon\in(0,q_*)\). Then there exists a fixed
\(\rho_\varepsilon\in(0,1)\) such that
\[
\mathbb P\left(
\widetilde P^{\,1}_{(j)}
\le (1-\rho_\varepsilon)c_j(\alpha_m),\quad
j=1,\ldots,r_m
\right)\to1.
\]
\end{lemma}

\begin{proof}
Fix any sufficiently small \(\rho\in(0,1)\), to be specified below, and put
\[
\alpha_m^{(1)}=(1-\rho)\alpha_m,
\qquad
\alpha_m^{(2)}=(1-\rho)^2\alpha_m.
\]
Because \(\rho\) is fixed and \(\alpha_m\to0\),
\[
\log\frac{r_{\alpha_m^{(\ell)}}}{r_{\alpha_m}}=O(1),
\qquad \ell=1,2.
\]
Since \(u_{m,\lambda}\to\infty\), both modified levels satisfy the same
first-order calibration,
\[
\frac{2\log(f_m/r_{\alpha_m^{(\ell)}})}{u_{m,\lambda}}
\to C_\lambda,
\qquad \ell=1,2,
\]
and \(\alpha_m^{(\ell)}f_m\to\infty\).

We repeat the independent empirical-crossing argument at the tightened
thresholds
\[
d_j=(1-\rho)c_j(\alpha_m).
\]
For the bulk ranks, whenever \(j(1-\alpha_m)>1\),
\[
c_j(\alpha_m)\ge \frac{j\alpha_m}{m},
\]
and hence, with \(q_j=j/(mp_m)\),
\[
d_j\ge q_j\alpha_m^{(1)}p_m.
\]
Applying the factor-known crossing argument proved in
Lemma~\ref{lem:residual-deterministic-signal-crossing} at nominal level \(\alpha_m^{(1)}\), we obtain, for
some \(0<\eta_\varepsilon\le1\) (shrinking the margin if necessary),
\[
G_m(d_j)
\ge G_m(q_j\alpha_m^{(1)}p_m)
\ge (1+\eta_\varepsilon)q_j
\]
uniformly for \(q_j\le q_*-\varepsilon\), once \(m\) is sufficiently
large.  The choice of a fixed \(\rho\) does not change the limiting crossing
point because its contribution to the logarithmic calibration is only
\(O(1)/u_{m,\lambda}=o(1)\).

For the earliest ranks we need a tightened initial crossing.  Since
\[
c_1(\alpha_m^{(2)})
=\frac{(1-\rho)^2\alpha_m}{m+(1-\rho)^2\alpha_m}
<\frac{(1-\rho)\alpha_m}{m+\alpha_m}
=(1-\rho)c_1(\alpha_m)
\]
for all sufficiently large \(m\), monotonicity of \(G_m\) and the
initial-crossing argument in the proof of
Lemma~\ref{lem:residual-initial-crossing} at level \(\alpha_m^{(2)}\) give
\[
mp_mG_m(d_1)
\ge mp_mG_m\{c_1(\alpha_m^{(2)})\}
\to\infty.
\]

Choose integers \(L_m\to\infty\) such that
\[
L_m=o(mp_m),
\qquad
L_m=o\{mp_mG_m(d_1)\}.
\]
Conditional on \(K_m=k\), the number of signal \(p\)-values below \(d_j\)
is \(\operatorname{Bin}\{k,G_m(d_j)\}\).  The preceding tightened initial
crossing handles \(1\le j<L_m\) exactly as in the independent GBS proof.
For \(L_m\le j\le r_m\), on the event
\[
K_m\ge(1-\eta_\varepsilon/4)mp_m,
\]
whose probability tends to one, the bulk bound gives
\[
K_mG_m(d_j)\ge(1+\eta_\varepsilon/2)j
\]
for all sufficiently large \(m\), uniformly over these ranks.  A Chernoff
bound followed by a union bound therefore yields
\[
\mathbb P\left(
N_1(d_j)\ge j,\quad j=1,\ldots,r_m
\right)\to1.
\]
Since \(N_1(d_j)\ge j\) is equivalent to
\(\widetilde P^{\,1}_{(j)}\le d_j\), the result follows.  Taking, for
example, any sufficiently small fixed \(\rho=\rho_\varepsilon\) completes
the proof.
\end{proof}

\begin{lemma}[Residual empirical signal-crossing bound]
	\label{lem:residual-empirical-signal-crossing}
	Assume the conditions of Lemmas~\ref{lem:residual-deterministic-signal-crossing},
	\ref{lem:residual-initial-crossing}, and
	\ref{lem:uniform-residual-perturbation}. Assume also
    \(\alpha_m f_m\to\infty\), as in
    Lemma~\ref{lem:tightened-factor-known-crossing}. Fix
	\(\varepsilon\in(0,q_*)\), and define
	\[
	r_m=\lfloor(q_*-\varepsilon)mp_m\rfloor.
	\]
	Let
	\[
	\widehat P^{\,1}_{(1)}
	\le
	\cdots
	\le
	\widehat P^{\,1}_{(K_m)}
	\]
	denote the ordered residualized \(p\)-values corresponding to the
	non-null hypotheses. Then
	\[
	\mathbb P\left(
	\widehat P^{\,1}_{(j)}\le c_j,\ 
	j=1,\ldots,r_m
	\right)
	\to1.
	\]
	Consequently,
	\[
	\mathbb P\left(
	R_m^{\mathrm{rGBS}}\ge r_m
	\right)
	\to1.
	\]
\end{lemma}

\begin{proof}
	Let \(\rho=\rho_\varepsilon\in(0,1)\) be the constant supplied by Lemma~\ref{lem:tightened-factor-known-crossing}. Let
	\[
	\widetilde P_i=2\bar\Phi(|Y_i(Z)|),
	\qquad
	\widehat P_i=2\bar\Phi(|\widehat Y_i|)
	\]
	denote the factor-known and residualized \(p\)-values, respectively.
	Let
	\[
	\widetilde P^{\,1}_{(1)}\le\cdots\le
	\widetilde P^{\,1}_{(K_m)}
	\quad\text{and}\quad
	\widehat P^{\,1}_{(1)}\le\cdots\le
	\widehat P^{\,1}_{(K_m)}
	\]
	denote the corresponding ordered signal \(p\)-values.
	
	By Lemma~\ref{lem:tightened-factor-known-crossing}, applied to the
	factor-known residual experiment, we have
	\begin{equation}
		\label{eq:factor-known-tight-crossing}
		\mathbb P\left(
		\widetilde P^{\,1}_{(j)}
		\le
		(1-\rho)c_j,
		\ j=1,\ldots,r_m
		\right)
		\to1.
	\end{equation}
	
	We now transfer this event to the residualized \(p\)-values. For
	\(0<t<1\), write
	\[
	z_t=\Phi^{-1}\left(1-\frac t2\right).
	\]
	Since \(\rho\in(0,1)\) is fixed, there exists \(a_\rho>0\) such that,
	for all sufficiently small \(t\),
	\begin{equation}
		\label{eq:quantile-separation-empirical}
		z_{(1-\rho)t}
		\ge
		z_t+\frac{a_\rho}{z_t}.
	\end{equation}
	Indeed, by Mills' ratio,
	\[
	\frac{\bar\Phi(z_t+a/z_t)}{\bar\Phi(z_t)}
	\to e^{-a},
	\qquad t\downarrow0,
	\]
	Choose \(a_\rho\) so that
	\[
	0<a_\rho<-\log(1-\rho),
	\]
	so that \(e^{-a_\rho}>1-\rho\). Hence, for all sufficiently small
	\(t\),
	\[
	\bar\Phi(z_t+a_\rho/z_t)>(1-\rho)\bar\Phi(z_t),
	\]
	which, since \(x\mapsto\bar\Phi(x)\) is decreasing, is equivalent to
	\eqref{eq:quantile-separation-empirical}.
	
	By Lemma~\ref{lem:uniform-residual-perturbation},
	\[
	z_{c_1}\max_{1\le i\le m}|\Delta_{-i}|=o_{\mathbb P}(1).
	\]
	Hence, with probability tending to one,
	\begin{equation}
		\label{eq:uniform-delta-control-empirical}
		z_{c_1}\max_{1\le i\le m}|\Delta_{-i}|
		\le
		\frac{a_\rho}{2}.
	\end{equation}
	Since the GBS critical constants are increasing, \(c_j\ge c_1\), and
	therefore \(z_{c_j}\le z_{c_1}\). Thus, on
	\eqref{eq:uniform-delta-control-empirical},
	\[
	|\Delta_{-i}|
	\le
	\frac{a_\rho}{2z_{c_j}},
	\qquad
	i=1,\ldots,m,\quad j=1,\ldots,r_m.
	\]
	
	Now suppose that, for some non-null index \(i\),
	\[
	\widetilde P_i\le (1-\rho)c_j.
	\]
	Then
	\[
	|Y_i(Z)|\ge z_{(1-\rho)c_j}.
	\]
	By \eqref{eq:quantile-separation-empirical},
	\[
	|Y_i(Z)|
	\ge
	z_{c_j}+\frac{a_\rho}{z_{c_j}}.
	\]
	Therefore, on \eqref{eq:uniform-delta-control-empirical},
	\[
	|\widehat Y_i|
	=
	|Y_i(Z)-\Delta_{-i}|
	\ge
	|Y_i(Z)|-|\Delta_{-i}|
	>
	z_{c_j}.
	\]
	Hence
	\[
	\widehat P_i\le c_j.
	\]
	Consequently, on an event whose probability tends to one,
	\[
	\widetilde P^{\,1}_{(j)}\le (1-\rho)c_j,
	\ j=1,\ldots,r_m
	\quad\Longrightarrow\quad
	\widehat P^{\,1}_{(j)}\le c_j,
	\ j=1,\ldots,r_m.
	\]
	Combining this implication with
	\eqref{eq:factor-known-tight-crossing} yields
	\[
	\mathbb P\left(
	\widehat P^{\,1}_{(j)}\le c_j,
	\ j=1,\ldots,r_m
	\right)
	\to1.
	\]
On this event, the first \(r_m\) ordered residualized signal \(p\)-values
cross the corresponding GBS critical constants. Hence the residualized
GBS step-down procedure rejects at least \(r_m\) hypotheses, and therefore
\(R_m^{\mathrm{rGBS}}\ge r_m\). Therefore
	\[
	\mathbb P\left(
	R_m^{\mathrm{rGBS}}\ge r_m
	\right)\to1.
	\]
	This completes the proof.
\end{proof}

\begin{lemma}[Residual Type-II bound]
	\label{lem:residual-type-two-bound}
	Assume the conditions of Lemma~\ref{lem:residual-empirical-signal-crossing}.
	Then, for every \(\varepsilon\in(0,q_*)\),
	\[
	\mathbb E\left(T_m^{\mathrm{rGBS}}\right)
	\le
	mp_m(1-q_*+\varepsilon)+o(mp_m).
	\]
	Consequently,
	\[
	\mathbb E\left(T_m^{\mathrm{rGBS}}\right)
	\le
	mp_m\{2\Phi(\sqrt{C_\lambda})-1+\varepsilon\}
	+
	o(mp_m).
	\]
\end{lemma}
\begin{proof}
	Fix \(\varepsilon\in(0,q_*)\) and put
	\[
	r_m=\lfloor(q_*-\varepsilon)mp_m\rfloor .
	\]
	Let
	\[
	S_m^{\mathrm{rGBS}}
	\]
	denote the number of true discoveries made by the residualized GBS
	procedure. Then
	\[
	T_m^{\mathrm{rGBS}}
	=
	K_m-S_m^{\mathrm{rGBS}}.
	\]
	
	Define the event
	\[
	A_m
	=
	\left\{
	\widehat P^{\,1}_{(j)}\le c_j,\ 
	j=1,\ldots,r_m
	\right\}.
	\]
	By Lemma~\ref{lem:residual-empirical-signal-crossing},
	\[
	\mathbb P(A_m)\to1.
	\]
	On \(A_m\), for every \(j\le r_m\), the \(j\)th ordered
	residualized \(p\)-value among all hypotheses is no larger than the
	corresponding signal order statistic. Hence
	\[
	\widehat P_{(j)}
	\le
	\widehat P^{\,1}_{(j)}
	\le c_j,
	\qquad j=1,\ldots,r_m,
	\]
	and therefore the step-down rejection count satisfies
	\[
	R_m^{\mathrm{rGBS}}\ge r_m.
	\]
	We next verify that the first \(r_m\) non-null hypotheses themselves are
	rejected.  If \(i\) is one of the indices corresponding to
	\(\widehat P^{\,1}_{(1)},\ldots,\widehat P^{\,1}_{(r_m)}\), then
	\[
	\widehat P_i
	\le
	\widehat P^{\,1}_{(r_m)}
	\le c_{r_m}
	\le c_{R_m^{\mathrm{rGBS}}}.
	\]
	If such an index were not among the \(R_m^{\mathrm{rGBS}}\) rejected
	hypotheses, then, unless \(R_m^{\mathrm{rGBS}}=m\) (in which case there
	is nothing to prove), the \(R_m^{\mathrm{rGBS}}+1\)st ordered residualized
	\(p\)-value would satisfy
	\[
	\widehat P_{(R_m^{\mathrm{rGBS}}+1)}
	\le c_{R_m^{\mathrm{rGBS}}}
	< c_{R_m^{\mathrm{rGBS}}+1},
	\]
	contradicting maximality of the step-down rejection count.  Thus all of
	these \(r_m\) non-null hypotheses are rejected, and hence
	\[
	S_m^{\mathrm{rGBS}}\ge r_m
	\]
	on \(A_m\). Therefore, on \(A_m\),
	\[
	T_m^{\mathrm{rGBS}}
	=
	K_m-S_m^{\mathrm{rGBS}}
	\le
	(K_m-r_m)_+.
	\]
	Always,
	\[
	T_m^{\mathrm{rGBS}}\le K_m.
	\]
	Thus
	\[
	T_m^{\mathrm{rGBS}}
	\le
	(K_m-r_m)_+
	+
	K_m I(A_m^c).
	\]
	Taking expectations,
	\begin{equation}
		\label{eq:type-two-decomposition}
		\mathbb E(T_m^{\mathrm{rGBS}})
		\le
		\mathbb E\{(K_m-r_m)_+\}
		+
		\mathbb E\{K_m I(A_m^c)\}.
	\end{equation}
	
	Since \(K_m\sim\mathrm{Bin}(m,p_m)\) and \(mp_m\to\infty\),
	\[
	\frac{K_m}{mp_m}\to1
	\qquad\text{in }L^1.
	\]
	Therefore,
	\[
	\mathbb E|K_m-mp_m|
	=
	o(mp_m).
	\]
	Using
	\[
	(K_m-r_m)_+
	\le
	(mp_m-r_m)+|K_m-mp_m|,
	\]
	we obtain
	\[
	\mathbb E\{(K_m-r_m)_+\}
	\le
	mp_m-r_m+o(mp_m).
	\]
	Since
	\[
	r_m=\lfloor(q_*-\varepsilon)mp_m\rfloor,
	\]
	we have
	\[
	mp_m-r_m
	=
	mp_m(1-q_*+\varepsilon)+O(1).
	\]
	Because \(mp_m\to\infty\),
	\[
	\mathbb E\{(K_m-r_m)_+\}
	\le
	mp_m(1-q_*+\varepsilon)+o(mp_m).
	\]
	
	Next, by Cauchy--Schwarz,
	\[
	\mathbb E\{K_m I(A_m^c)\}
	\le
	\{\mathbb E(K_m^2)\}^{1/2}
	\mathbb P(A_m^c)^{1/2}.
	\]
	Since \(K_m\sim\mathrm{Bin}(m,p_m)\),
	\[
	\mathbb E(K_m^2)
	=
	mp_m(1-p_m)+(mp_m)^2
	=
	O\{(mp_m)^2\}.
	\]
	Also, \(\mathbb P(A_m^c)\to0\). Hence
	\[
	\mathbb E\{K_m I(A_m^c)\}
	=
	o(mp_m).
	\]
	Combining this with
	\eqref{eq:type-two-decomposition}, we get
	\[
	\mathbb E(T_m^{\mathrm{rGBS}})
	\le
	mp_m(1-q_*+\varepsilon)+o(mp_m).
	\]
	Finally, since
	\[
	q_*=2\bar\Phi(\sqrt{C_\lambda}),
	\]
	we have
	\[
	1-q_*
	=
	2\Phi(\sqrt{C_\lambda})-1.
	\]
	Therefore
	\[
	\mathbb E(T_m^{\mathrm{rGBS}})
	\le
	mp_m\{2\Phi(\sqrt{C_\lambda})-1+\varepsilon\}
	+
	o(mp_m).
	\]
	This completes the proof.
\end{proof}

\begin{lemma}[Shifted uniform order-statistic inequality]
\label{lem:shifted-order-statistic}
Let \(U_{(1)}\le\cdots\le U_{(n)}\) be the order statistics of
\(n\) independent \(\mathrm{Unif}(0,1)\) random variables.  Let \(k\ge0\)
and \(0<\beta<1\), and, for \(\ell=1,\ldots,n\), define
\[
b_{\ell,k}^{(n)}(\beta)
=
\frac{(k+\ell)\beta}
{n+k+1-(k+\ell)(1-\beta)}.
\]
Then
\begin{equation}
\label{eq:shifted-order-statistic-bound}
\sum_{j=1}^n
\mathbb P\left\{
U_{(\ell)}\le b_{\ell,k}^{(n)}(\beta),\
\ell=1,\ldots,j
\right\}
\le
\frac{\beta}{1-\beta}(k+1).
\end{equation}
\end{lemma}

\begin{proof}
For \(j=1,\ldots,n\), put
\[
A_j^{(n,k,\beta)}
=
\left\{
U_{(\ell)}\le b_{\ell,k}^{(n)}(\beta),\
\ell=1,\ldots,j
\right\},
\qquad
S_n(k,\beta)=\sum_{j=1}^n\mathbb P(A_j^{(n,k,\beta)}).
\]
We prove by induction on \(n\) that
\begin{equation}
\label{eq:shifted-induction-claim}
S_n(k,\beta)
\le
\frac{\beta}{1-\beta}(k+1)
\end{equation}
for every \(n\ge1\), \(k\ge0\), and \(0<\beta<1\).

For \(n=1\),
\[
S_1(k,\beta)
=b_{1,k}^{(1)}(\beta)
=
\frac{(k+1)\beta}{1+(k+1)\beta}
\le
\frac{\beta}{1-\beta}(k+1).
\]
Now suppose that \eqref{eq:shifted-induction-claim} holds for \(n-1\),
where \(n\ge2\).  Write
\[
b_1=b_{1,k}^{(n)}(\beta)
=
\frac{(k+1)\beta}{n+(k+1)\beta}.
\]
Condition on \(U_{(1)}=u\), where \(0<u<b_1\).  Given this value, the
transformed remaining order statistics
\[
W_{(\ell)}
=
\frac{U_{(\ell+1)}-u}{1-u},
\qquad \ell=1,\ldots,n-1,
\]
have the joint distribution of the order statistics of \(n-1\)
independent \(\mathrm{Unif}(0,1)\) random variables.  Therefore
\begin{equation}
\label{eq:shifted-recursion}
S_n(k,\beta)
=
\int_0^{b_1}n(1-u)^{n-1}
\{1+\widetilde S_{n-1}(u)\}\,du,
\end{equation}
where
\[
\widetilde S_{n-1}(u)
=
\sum_{j=1}^{n-1}
\mathbb P\left\{
W_{(\ell)}
\le
\frac{b_{\ell+1,k}^{(n)}(\beta)-u}{1-u},\
\ell=1,\ldots,j
\right\}.
\]
Put
\[
x=\frac{u}{1-u},
\qquad
\beta_u=\frac{\beta+x}{1+x},
\qquad
k_u=\frac{(k+1)\beta-nx}{\beta+x}.
\]
Since \(0<u<b_1\), we have
\[
0<x<\frac{(k+1)\beta}{n},
\]
so \(k_u>0\), while \(0<\beta_u<1\).  A direct odds-ratio calculation
shows that, for every \(\ell=1,\ldots,n-1\),
\begin{equation}
\label{eq:shifted-boundary-identity}
\frac{b_{\ell+1,k}^{(n)}(\beta)-u}{1-u}
=
b_{\ell,k_u}^{(n-1)}(\beta_u).
\end{equation}
Indeed, using
\[
\frac{b_{\ell,k}^{(n)}(\beta)}
{1-b_{\ell,k}^{(n)}(\beta)}
=
\frac{(k+\ell)\beta}{n+1-\ell},
\]
the odds corresponding to the left-hand side of
\eqref{eq:shifted-boundary-identity} equal
\[
\frac{(k+\ell+1)\beta-(n-\ell)x}
{(1+x)(n-\ell)},
\]
which is exactly
\[
\frac{(k_u+\ell)\beta_u}{n-\ell},
\]
the odds corresponding to \(b_{\ell,k_u}^{(n-1)}(\beta_u)\).
Thus
\[
\widetilde S_{n-1}(u)=S_{n-1}(k_u,\beta_u).
\]
By the induction hypothesis,
\[
\widetilde S_{n-1}(u)
\le
\frac{\beta_u}{1-\beta_u}(k_u+1)
=
\frac{(k+2)\beta-(n-1)x}{1-\beta}.
\]
Hence
\begin{align}
S_n(k,\beta)
&\le
\frac{1}{1-\beta}
\int_0^{b_1} n(1-u)^{n-1}
\left\{1+(k+1)\beta-(n-1)\frac{u}{1-u}\right\}\,du.
\label{eq:shifted-integral}
\end{align}
Let \(A=(k+1)\beta\), so \(b_1=A/(n+A)\).  The integral in
\eqref{eq:shifted-integral}, before multiplication by \((1-\beta)^{-1}\),
is
\[
I
=
\int_0^{b_1}
\left[n(1+A)(1-u)^{n-1}
-n(n-1)u(1-u)^{n-2}\right]du.
\]
Since
\[
\frac{d}{du}\{nu(1-u)^{n-1}\}
=
n(1-u)^{n-1}-n(n-1)u(1-u)^{n-2},
\]
we obtain
\[
I
=nb_1(1-b_1)^{n-1}
+A\{1-(1-b_1)^n\}.
\]
Using \(b_1=A/(n+A)\) and \(1-b_1=n/(n+A)\), the first term is
\[
nb_1(1-b_1)^{n-1}
=A\left(\frac{n}{n+A}\right)^n,
\]
and therefore \(I=A\).  Consequently
\[
S_n(k,\beta)
\le
\frac{A}{1-\beta}
=
\frac{\beta}{1-\beta}(k+1),
\]
which completes the induction.
\end{proof}

\begin{lemma}[Finite-sample Type-I bound for factor-known GBS]
\label{lem:factor-known-gbs-type-one}
Consider the factor-known residual experiment with independent \(p\)-values
\(\widetilde P_i=2\bar\Phi(|Y_i(Z)|)\), and apply the GBS step-down rule at a
level \(\beta\in(0,1)\).  If \(K_m=\sum_{i=1}^m\nu_i\) is the number of
non-null coordinates and \(\widetilde V_m(\beta)\) is the number of false
discoveries, then
\[
\mathbb E\{\widetilde V_m(\beta)\mid K_m=k\}
\le
\frac{\beta}{1-\beta}(k+1).
\]
Consequently,
\begin{equation}
\label{eq:factor-known-gbs-type-one-unconditional}
\mathbb E\{\widetilde V_m(\beta)\}
\le
\frac{\beta}{1-\beta}(mp_m+1).
\end{equation}
\end{lemma}

\begin{proof}
Fix \(K_m=k\), put \(n=m-k\), and let
\(U_{(1)}\le\cdots\le U_{(n)}\) denote the ordered null factor-known
\(p\)-values.  Conditional on the non-null configuration (and on their
\(p\)-values), these are the order statistics of \(n\) independent
\(\mathrm{Unif}(0,1)\) variables.

If the GBS procedure makes at least \(j\) false discoveries, then the
\(j\) smallest null \(p\)-values are rejected.  For
\(\ell=1,\ldots,j\), let \(R_\ell\) be the overall rank of \(U_{(\ell)}\).
At most \(k\) non-null \(p\)-values can precede \(U_{(\ell)}\), so
\[
R_\ell\le k+\ell.
\]
Because a rejected hypothesis of overall rank \(R_\ell\) satisfies the
corresponding step-down inequality and because the GBS critical constants
are increasing,
\[
U_{(\ell)}
\le c_{R_\ell}(\beta)
\le c_{k+\ell}(\beta),
\qquad \ell=1,\ldots,j.
\]
Thus
\[
\{\widetilde V_m(\beta)\ge j\}
\subseteq
\left\{
U_{(\ell)}\le c_{k+\ell}(\beta),\
\ell=1,\ldots,j
\right\}.
\]
Since \(n=m-k\),
\[
c_{k+\ell}(\beta)
=
\frac{(k+\ell)\beta}
{n+k+1-(k+\ell)(1-\beta)}
=
b_{\ell,k}^{(n)}(\beta).
\]
Using \(\mathbb E(Y)=\sum_{j\ge1}\mathbb P(Y\ge j)\) for a nonnegative
integer-valued random variable and then
Lemma~\ref{lem:shifted-order-statistic},
\begin{align*}
\mathbb E\{\widetilde V_m(\beta)\mid K_m=k\}
&\le
\sum_{j=1}^n
\mathbb P\left\{
U_{(\ell)}\le b_{\ell,k}^{(n)}(\beta),\
\ell=1,\ldots,j
\right\}\\
&\le
\frac{\beta}{1-\beta}(k+1).
\end{align*}
Averaging over \(K_m\sim\mathrm{Bin}(m,p_m)\) yields
\[
\mathbb E\{\widetilde V_m(\beta)\}
\le
\frac{\beta}{1-\beta}\{\mathbb E(K_m)+1\}
=
\frac{\beta}{1-\beta}(mp_m+1).
\]
\end{proof}

\begin{lemma}[Inflated-level coupling up to an integer half-path]
\label{lem:inflated-level-coupling}
Fix $\rho>0$ and put
\[
\alpha_m'=(1+\rho)\alpha_m,
\qquad
r_m^\star=\left\lfloor\frac m2\right\rfloor.
\]
For all sufficiently large $m$, $\alpha_m'<1$.  Let
$\widetilde R_m'$ and $\widetilde{\mathcal R}_m'$ denote the rejection
count and rejection set of factor-known GBS run at level $\alpha_m'$, and
let $\mathcal R_m^{\mathrm{rGBS}}$ denote the rGBS rejection set at level
$\alpha_m$.

There exist constants $\eta_\rho>0$ and $a_\rho>0$ such that, with
\[
D_m=
\left\{
 z_{c_1(\alpha_m)}\max_{1\le i\le m}|\Delta_{-i}|
 \le a_\rho
\right\},
\]
the following assertions hold for all sufficiently large $m$:
\begin{enumerate}[label=(\roman*)]
\item uniformly for $1\le j\le r_m^\star$,
\[
 c_j(\alpha_m')\ge(1+\eta_\rho)c_j(\alpha_m);
\]
\item on $D_m$, uniformly over $i=1,\ldots,m$ and
$1\le j\le r_m^\star$,
\begin{equation}
\label{eq:inflated-level-pvalue-transfer}
 \widehat P_i\le c_j(\alpha_m)
 \quad\Longrightarrow\quad
 \widetilde P_i\le c_j(\alpha_m');
\end{equation}
\item on $D_m\cap\{R_m^{\mathrm{rGBS}}\le r_m^\star\}$,
\begin{equation}
\label{eq:rgbs-factor-known-set-inclusion}
 \mathcal R_m^{\mathrm{rGBS}}
 \subseteq
 \widetilde{\mathcal R}_m';
\end{equation}
while on $D_m\cap\{R_m^{\mathrm{rGBS}}>r_m^\star\}$,
\[
 \widetilde R_m'\ge r_m^\star.
\]
\end{enumerate}
\end{lemma}

\begin{proof}
Write
\[
A_{m,j}=m+1-j.
\]
A direct calculation gives
\[
\frac{c_j(\alpha_m')}{c_j(\alpha_m)}
=
\frac{(1+\rho)\{A_{m,j}+j\alpha_m\}}
{A_{m,j}+(1+\rho)j\alpha_m}.
\]
Hence
\[
\frac{c_j(\alpha_m')}{c_j(\alpha_m)}-1
=
\frac{\rho A_{m,j}}
{A_{m,j}+(1+\rho)j\alpha_m}.
\]
If $j\le r_m^\star$, then $j/A_{m,j}\le1$, so, because
$\alpha_m\to0$, the last display is bounded below by a fixed positive
constant for all sufficiently large $m$.  This proves (i).

Moreover $c_{r_m^\star}(\alpha_m)=O(\alpha_m)\to0$.  By part (i),
there is a fixed $\eta_\rho>0$ such that
\[
c_j(\alpha_m')\ge(1+\eta_\rho)c_j(\alpha_m),
\qquad 1\le j\le r_m^\star,
\]
for all sufficiently large $m$.  Choose any fixed
\[
0<a_\rho<\log(1+\eta_\rho).
\]
For $t\downarrow0$, Mills' ratio gives
\[
\frac{\bar\Phi(z_t-a_\rho/z_t)}{\bar\Phi(z_t)}
\longrightarrow e^{a_\rho}.
\]
Because $e^{a_\rho}<1+\eta_\rho$ and
$c_{r_m^\star}(\alpha_m)\to0$, the preceding convergence applies
uniformly to $t=c_j(\alpha_m)$, $1\le j\le r_m^\star$, for all
sufficiently large $m$.  Hence
\[
2\bar\Phi\left(z_{c_j(\alpha_m)}-
\frac{a_\rho}{z_{c_j(\alpha_m)}}\right)
<
(1+\eta_\rho)c_j(\alpha_m)
\le c_j(\alpha_m').
\]
Since $t\mapsto z_t$ is decreasing, this is equivalent to
\begin{equation}
\label{eq:inflated-level-quantile-gap}
 z_{c_j(\alpha_m)}-
 \frac{a_\rho}{z_{c_j(\alpha_m)}}
 \ge
 z_{c_j(\alpha_m')},
 \qquad 1\le j\le r_m^\star.
\end{equation}
On $D_m$, monotonicity of $j\mapsto c_j(\alpha_m)$ gives
\[
|\Delta_{-i}|
\le
\frac{a_\rho}{z_{c_1(\alpha_m)}}
\le
\frac{a_\rho}{z_{c_j(\alpha_m)}}.
\]
Thus, if $\widehat P_i\le c_j(\alpha_m)$, then
\[
|Y_i(Z)|
\ge
|\widehat Y_i|-|\Delta_{-i}|
\ge
z_{c_j(\alpha_m)}-
\frac{a_\rho}{z_{c_j(\alpha_m)}}
\ge
z_{c_j(\alpha_m')},
\]
which proves \eqref{eq:inflated-level-pvalue-transfer}.

Suppose now that $R_m^{\mathrm{rGBS}}=r\le r_m^\star$.  For every
$j\le r$, the step-down definition gives
$\widehat P_{(j)}\le c_j(\alpha_m)$.  By
\eqref{eq:inflated-level-pvalue-transfer}, at least $j$ factor-known
$p$-values are no larger than $c_j(\alpha_m')$, and hence
\[
\widetilde P_{(j)}\le c_j(\alpha_m'),
\qquad j=1,\ldots,r.
\]
Therefore $\widetilde R_m'\ge r$.  If $i\in\mathcal R_m^{\mathrm{rGBS}}$,
then $\widehat P_i\le c_r(\alpha_m)$ and hence
$\widetilde P_i\le c_r(\alpha_m')\le c_{\widetilde R_m'}(\alpha_m')$.
If $\widetilde R_m'=m$, this already implies that $i$ is rejected.
Otherwise, if $i$ did not belong to the factor-known GBS rejection set,
then the first $\widetilde R_m'$ ordered factor-known $p$-values together
with $\widetilde P_i$ would give at least $\widetilde R_m'+1$ factor-known
$p$-values no larger than $c_{\widetilde R_m'}(\alpha_m')$.  Hence
\[
\widetilde P_{(\widetilde R_m'+1)}
\le c_{\widetilde R_m'}(\alpha_m')
< c_{\widetilde R_m'+1}(\alpha_m'),
\]
which would force the next step-down inequality to hold and contradict
maximality of $\widetilde R_m'$.  Thus $i$ belongs to the factor-known
GBS rejection set, proving \eqref{eq:rgbs-factor-known-set-inclusion}.

Finally, if $R_m^{\mathrm{rGBS}}>r_m^\star$, the first
$r_m^\star$ rGBS step-down inequalities hold.  Applying the same transfer
argument through rank $r_m^\star$ gives
$\widetilde R_m'\ge r_m^\star$.
\end{proof}

\begin{lemma}[Quantitative control of the residual perturbation event]
\label{lem:quantitative-residual-perturbation}
Assume
\[
mp_m\to\infty,
\qquad
\alpha_m\to0,
\qquad
\alpha_m f_m\to\infty,
\]
and
\[
\frac{2\log(f_m/r_{\alpha_m})}{u_{m,\lambda}}
\to C_\lambda\in(0,\infty).
\]
For the event $D_m$ in Lemma~\ref{lem:inflated-level-coupling}, there are
constants $C_1,C_2,c_1,c_2>0$ such that, for all sufficiently large $m$,
\begin{equation}
\label{eq:Dm-quantitative-tail}
\mathbb P(D_m^c)
\le
C_1\exp\left\{-c_1\frac{m}{\log^2m}\right\}
+
C_2m\exp\left\{-c_2\frac{m^2}{\log^2m}\right\}.
\end{equation}
In particular, if also $\alpha_m\delta_m\to0$, then
\begin{equation}
\label{eq:Dm-risk-negligible}
\delta_m m\mathbb P(D_m^c)=o(mp_m).
\end{equation}
\end{lemma}

\begin{proof}
Since $mp_m\to\infty$, eventually $p_m>1/m$.  Since
$\alpha_m f_m\to\infty$ and $f_m\sim p_m^{-1}$, eventually
$\alpha_m>p_m$.  Hence
\[
\log(m/\alpha_m)=O(\log m).
\]
The GBS calibration then implies
\[
u_{m,\lambda}=O(\log m).
\]
Also
\[
z_{c_1(\alpha_m)}^2=O\{\log(m/\alpha_m)\}=O(\log m).
\]

Recall
\[
B_i=\varepsilon_i+\frac{\theta_i}{a_\lambda},
\qquad
\Delta_{-i}=\frac{1}{m-1}\left(\sum_{j=1}^mB_j-B_i\right).
\]
The $B_i$ are independent and centered.  Their moment generating functions
satisfy
\[
\mathbb E e^{tB_i}
=(1-p_m)e^{t^2/2}+p_me^{(1+u_{m,\lambda})t^2/2}
\le
e^{(1+u_{m,\lambda})t^2/2},
\]
so each $B_i$ is sub-Gaussian with variance proxy
$1+u_{m,\lambda}$.  Therefore
\[
\max_i|\Delta_{-i}|
\le
\frac{|\sum_{j=1}^mB_j|}{m-1}
+
\frac{\max_i|B_i|}{m-1},
\]
and standard sub-Gaussian tail bounds, with threshold proportional to
$m/z_{c_1(\alpha_m)}$, give
\[
\mathbb P(D_m^c)
\le
C_1\exp\left\{-c_1
\frac{m}{(1+u_{m,\lambda})z_{c_1(\alpha_m)}^2}\right\}
+
C_2m\exp\left\{-c_2
\frac{m^2}{(1+u_{m,\lambda})z_{c_1(\alpha_m)}^2}\right\}.
\]
Using the preceding logarithmic bounds yields
\eqref{eq:Dm-quantitative-tail}.

Finally, $\alpha_m\delta_m\to0$ implies eventually
$\delta_m<1/\alpha_m$.  Since eventually $\alpha_m>p_m>1/m$,
\[
\frac{\delta_m}{p_m}
<
\frac{1}{\alpha_mp_m}
<
\frac{1}{p_m^2}
<
m^2.
\]
The super-polynomial bound \eqref{eq:Dm-quantitative-tail} therefore gives
$(\delta_m/p_m)\mathbb P(D_m^c)\to0$, which is equivalent to
\eqref{eq:Dm-risk-negligible}.
\end{proof}

\begin{lemma}[Residual Type-I risk bound]
\label{lem:residual-type-one-risk}
Assume the conditions of Lemmas~\ref{lem:inflated-level-coupling} and
\ref{lem:quantitative-residual-perturbation}, and suppose
$\alpha_m\delta_m\to0$.  Then
\begin{equation}
\label{eq:residual-type-one-risk-negligible}
\delta_m\,\mathbb E(V_m^{\mathrm{rGBS}})=o(mp_m).
\end{equation}
Equivalently,
\[
\delta_{0,m}\mathbb E(V_m^{\mathrm{rGBS}})
=o(mp_m\delta_{A,m}).
\]
\end{lemma}

\begin{proof}
Fix $\rho>0$ and use the notation of
Lemma~\ref{lem:inflated-level-coupling}.  Let
$\widetilde V_m'=\widetilde V_m(\alpha_m')$.  The pathwise conclusions of
that lemma imply
\begin{equation}
\label{eq:pathwise-type-one-master}
V_m^{\mathrm{rGBS}}
\le
\widetilde V_m'
+mI(D_m^c)
+mI\{\widetilde R_m'\ge r_m^\star\}.
\end{equation}
By Lemma~\ref{lem:factor-known-gbs-type-one}, since
$\alpha_m'=(1+\rho)\alpha_m$,
\begin{equation}
\label{eq:inflated-factor-known-V-bound}
\mathbb E(\widetilde V_m')=O(\alpha_mmp_m).
\end{equation}

Put
\[
s_m^\star=\left\lfloor\frac m4\right\rfloor.
\]
Since $\widetilde R_m'\le K_m+\widetilde V_m'$, for all sufficiently
large $m$,
\[
\{\widetilde R_m'\ge r_m^\star\}
\subseteq
\{K_m\ge s_m^\star\}
\cup
\{\widetilde V_m'\ge r_m^\star-s_m^\star\}.
\]
Because $p_m\to0$, a Chernoff bound for
$K_m\sim\mathrm{Bin}(m,p_m)$ gives
\[
\mathbb P(K_m\ge s_m^\star)\le e^{-c m\log(1/p_m)}
\]
for some $c>0$ and all sufficiently large $m$.  Also, by Markov's
inequality and $r_m^\star-s_m^\star\asymp m$,
\[
m\mathbb P\{\widetilde V_m'\ge r_m^\star-s_m^\star\}
=O\{\mathbb E(\widetilde V_m')\}
=O(\alpha_mmp_m).
\]
Thus
\begin{equation}
\label{eq:half-path-probability-bound}
m\mathbb P(\widetilde R_m'\ge r_m^\star)
=O(\alpha_mmp_m)+o(mp_m/\delta_m).
\end{equation}
Taking expectations in \eqref{eq:pathwise-type-one-master}, multiplying
by $\delta_m$, and using
\eqref{eq:inflated-factor-known-V-bound},
\eqref{eq:Dm-risk-negligible}, and
\eqref{eq:half-path-probability-bound}, we obtain
\[
\delta_m\mathbb E(V_m^{\mathrm{rGBS}})
=
O(\alpha_m\delta_m\,mp_m)+o(mp_m)
=o(mp_m),
\]
because $\alpha_m\delta_m\to0$.
\end{proof}

The Bayes risk of rGBS, distinct from its rejection count
\(R_m^{\mathrm{rGBS}}\), is
\begin{equation}
\label{eq:rgbs-bayes-risk-definition}
R_{m,\lambda}^{\mathrm{rGBS}}
=\delta_{0,m}\mathbb E(V_m^{\mathrm{rGBS}})
+\delta_{A,m}\mathbb E(T_m^{\mathrm{rGBS}}).
\end{equation}

\subsection{Proof of Theorem~\ref{thm:rgbs-abos}}

\begin{proof}
First note that \eqref{eq:rgbs-abos-calibration} and
\eqref{eq:rgbs-abos-compatibility} imply the Bayes-Oracle calibration.
Indeed,
\[
\log\{u_{m,\lambda}f_m^2\delta_m^2\}
=
2\log(f_m/r_{\alpha_m})
+
\log u_{m,\lambda}
+
2\log(r_{\alpha_m}\delta_m).
\]
After division by $u_{m,\lambda}$, the first term converges to
$C_\lambda$, while the remaining terms vanish.  Moreover, Lemma~\ref{lem:quantitative-residual-perturbation} gives
$u_{m,\lambda}=O(\log m)$ and $\log(m/\alpha_m)=O(\log m)$. Hence
\[
\frac{u_{m,\lambda}\{1+p_m u_{m,\lambda}\}}{m}
=O\left(\frac{\log^2m}{m}\right)\to0,
\]
\[
\frac{\log(m/\alpha_m)\{1+p_m u_{m,\lambda}\}}{m}
=O\left(\frac{\log^2m}{m}\right)\to0,
\]
and
\[
\frac{\sqrt{\log(m/\alpha_m)}
\sqrt{(1+u_{m,\lambda})\log m}}{m}
=O\left(\frac{\log^{3/2}m}{m}\right)\to0.
\]
Thus every factor-estimation and uniform residual-perturbation condition used
above follows from the stated ABOS assumptions. Therefore
Theorem~\ref{thm:one-factor-oracle-risk} yields
\begin{equation}
\label{eq:oracle-risk-for-rgbs-abos}
R_{m,\lambda}^{\mathrm{BO}}
=
mp_m\delta_{A,m}
\{2\Phi(\sqrt{C_\lambda})-1\}
\{1+o(1)\}.
\end{equation}

By Lemma~\ref{lem:residual-type-one-risk},
\[
\delta_{0,m}\mathbb E(V_m^{\mathrm{rGBS}})
=o(mp_m\delta_{A,m}).
\]
On the other hand, Lemma~\ref{lem:residual-type-two-bound} gives, for
every fixed $\varepsilon\in(0,q_*)$,
\[
\mathbb E(T_m^{\mathrm{rGBS}})
\le
mp_m\{2\Phi(\sqrt{C_\lambda})-1+\varepsilon\}
+o(mp_m).
\]
Hence
\[
\limsup_{m\to\infty}
\frac{R_{m,\lambda}^{\mathrm{rGBS}}}
{mp_m\delta_{A,m}}
\le
2\Phi(\sqrt{C_\lambda})-1+\varepsilon.
\]
Letting $\varepsilon\downarrow0$ gives
\begin{equation}
\label{eq:rgbs-risk-limsup}
\limsup_{m\to\infty}
\frac{R_{m,\lambda}^{\mathrm{rGBS}}}
{mp_m\delta_{A,m}}
\le
2\Phi(\sqrt{C_\lambda})-1.
\end{equation}

Finally, the genuine Bayes Oracle minimizes Bayes risk among all
$\boldsymbol{X}$-measurable multiple-testing rules, and rGBS is $\boldsymbol{X}$-measurable.
Therefore
\[
R_{m,\lambda}^{\mathrm{rGBS}}
\ge
R_{m,\lambda}^{\mathrm{BO}}.
\]
Combining this lower bound with \eqref{eq:oracle-risk-for-rgbs-abos} and
\eqref{eq:rgbs-risk-limsup} proves both
\eqref{eq:rgbs-abos-ratio} and \eqref{eq:rgbs-risk-expansion}.
\end{proof}

\section{Heterogeneous model: auxiliary results and proofs}
\label{sec:heterogeneous-one-factor}

The homogeneous theory above extends to a known heterogeneous loading
structure.  This section develops that extension while retaining the same
sparse two-groups model for the means.  The main difference is that the
factor-known signal strength, and therefore the local Bayes-Oracle miss
probability, now vary with the coordinate.

Consider
\begin{equation}
\label{eq:heterogeneous-factor-model}
X_i=\theta_i+\lambda_{i,m}Z+\sigma_{i,m}\varepsilon_i,
\qquad
\sigma_{i,m}^2=1-\lambda_{i,m}^2,
\qquad i=1,\ldots,m,
\end{equation}
where \(Z,\varepsilon_1,\ldots,\varepsilon_m\) are mutually independent
standard normal random variables.  The deterministic loading array is known.
The latent indicators and means satisfy
\[
\nu_i\stackrel{\mathrm{i.i.d.}}{\sim}\operatorname{Bernoulli}(p_m),
\qquad
\theta_i\mid\nu_i=0=0,
\qquad
\theta_i\mid\nu_i=1\sim N(0,\psi_m^2),
\]
independently across coordinates and independently of the Gaussian errors.
As before, put
\[
f_m=\frac{1-p_m}{p_m},
\qquad
\delta_m=\frac{\delta_{0,m}}{\delta_{A,m}},
\qquad
r_{\alpha_m}=\frac{\alpha_m}{1-\alpha_m}.
\]

We impose the loading conditions
\begin{equation}
\label{eq:hetero-loading-upper}
\sup_m\max_{1\le i\le m}|\lambda_{i,m}|
\le \bar\lambda<1
\end{equation}
and
\begin{equation}
\label{eq:hetero-loading-information}
\liminf_{m\to\infty}
\frac1m\sum_{i=1}^m\lambda_{i,m}^2
\ge \kappa>0.
\end{equation}
Thus
\[
\sigma_{i,m}\ge \underline\sigma
:=\sqrt{1-\bar\lambda^2}>0
\]
uniformly in \(i\) and \(m\).

\subsection{Factor-known experiment and heterogeneous Bayes Oracle}

If \(Z\) were observed, define
\begin{equation}
\label{eq:hetero-factor-known-residual}
Y_i(Z)
=
\frac{X_i-\lambda_{i,m}Z}{\sigma_{i,m}}.
\end{equation}
Then the coordinates are mutually independent and
\[
Y_i(Z)\mid\nu_i=0\sim N(0,1),
\qquad
Y_i(Z)\mid\nu_i=1\sim N(0,1+u_{i,m}),
\]
where
\begin{equation}
\label{eq:hetero-effective-strength}
u_{i,m}
=
\frac{\psi_m^2}{\sigma_{i,m}^2}.
\end{equation}
The factor-known two-sided p-values
\[
\widetilde P_i=2\bar\Phi\{|Y_i(Z)|\}
\]
are therefore independent, with iid \(\operatorname{Unif}(0,1)\) null
p-values.

The coordinatewise factor-known Bayes Oracle rejects \(H_{0i}\) when
\[
|Y_i(Z)|>c_{i,m}^{\mathrm{BO}},
\]
where
\begin{equation}
\label{eq:hetero-oracle-threshold}
(c_{i,m}^{\mathrm{BO}})^2
=
\left(1+\frac1{u_{i,m}}\right)
\left[
\log\{u_{i,m}f_m^2\delta_m^2\}
+
\log\left(1+\frac1{u_{i,m}}\right)
\right].
\end{equation}

Define the finite-\(m\) average factor-known detection probability
\begin{equation}
\label{eq:Qm-definition}
Q_m
=
\frac1m\sum_{i=1}^m
2\bar\Phi\{\sqrt{C_0}\,\sigma_{i,m}\}.
\end{equation}
Under \eqref{eq:hetero-loading-upper}, there exist constants
\(0<q_-<q_+<1\) such that
\begin{equation}
\label{eq:Qm-bounds}
q_-\le Q_m\le q_+
\end{equation}
for all sufficiently large \(m\).

\begin{remark}[Computing the observed-data Oracle]
\label{rem:hetero-oracle-quadrature}
The one-factor representation also permits numerical evaluation of the
observed-data Oracle without enumerating all signal configurations, as in
the companion BSD manuscript. For \(d\in\{0,1\}\), write
\[
f_{d,i}(x_i\mid z)
=\phi_{\sigma_{i,m}^2+d\psi_m^2}(x_i-\lambda_{i,m}z),
\qquad
h_i(x_i\mid z)=(1-p_m)f_{0,i}(x_i\mid z)+p_m f_{1,i}(x_i\mid z),
\]
where \(\phi_v\) denotes the \(N(0,v)\) density. Conditional independence
and integration over the latent factor give
\[
\mathbb P(\nu_i=1\mid \boldsymbol{X}=\boldsymbol{x})
=\frac{\displaystyle\int_{\mathbb R}\phi(z)\,p_m f_{1,i}(x_i\mid z)
\prod_{j\ne i}h_j(x_j\mid z)\,dz}
{\displaystyle\int_{\mathbb R}\phi(z)\prod_{j=1}^m h_j(x_j\mid z)\,dz}.
\]
This is an exact representation; quadrature supplies a numerical
approximation. It integrates out \(Z\), rather than treating its realized
value as observed. The asymptotic Oracle-risk results below do not require
numerical evaluation of these integrals.
\end{remark}

\begin{lemma}[Factor-known heterogeneous Oracle risk]
\label{lem:hetero-factor-known-oracle}
Assume \eqref{eq:hetero-loading-upper} and
\[
p_m\to0,\qquad mp_m\to\infty,\qquad
\alpha_m\to0,\qquad \alpha_m f_m\to\infty,
\]
together with
\begin{equation}
\label{eq:hetero-gbs-calibration}
\frac{2\log(f_m/r_{\alpha_m})}{\psi_m^2}
\to C_0\in(0,\infty)
\end{equation}
and
\begin{equation}
\label{eq:hetero-compatibility}
\frac{\log\{1/(r_{\alpha_m}\delta_m)\}}{\psi_m^2}
\to0.
\end{equation}
Then, uniformly over \(1\le i\le m\),
\begin{equation}
\label{eq:hetero-oracle-threshold-limit}
\frac{(c_{i,m}^{\mathrm{BO}})^2}{1+u_{i,m}}
=
C_0\sigma_{i,m}^2+o(1).
\end{equation}
Moreover,
\begin{equation}
\label{eq:hetero-oracle-type-two}
\mathbb P\{|Y_i(Z)|\le c_{i,m}^{\mathrm{BO}}\mid\nu_i=1\}
=
2\Phi\{\sqrt{C_0}\sigma_{i,m}\}-1+o(1)
\end{equation}
uniformly in \(i\), while
\begin{equation}
\label{eq:hetero-oracle-type-one}
\sup_i
f_m\delta_m\,2\bar\Phi(c_{i,m}^{\mathrm{BO}})
\to0.
\end{equation}
Consequently, the factor-known Bayes-Oracle risk satisfies
\begin{equation}
\label{eq:hetero-factor-known-risk}
\mathcal R_m^{\mathrm{BO},Z}
=
mp_m\delta_{A,m}\{1-Q_m+o(1)\}.
\end{equation}
\end{lemma}

\begin{proof}
By \eqref{eq:hetero-oracle-threshold},
\[
\frac{(c_{i,m}^{\mathrm{BO}})^2}{1+u_{i,m}}
=
\frac{
\log\{u_{i,m}f_m^2\delta_m^2\}
+
\log(1+u_{i,m}^{-1})
}{u_{i,m}}.
\]
Using
\[
\log\{u_{i,m}f_m^2\delta_m^2\}
=
2\log(f_m/r_{\alpha_m})
+\log u_{i,m}
+2\log(r_{\alpha_m}\delta_m),
\]
the fact that \(u_{i,m}=\psi_m^2/\sigma_{i,m}^2\), and the uniform bounds
\(\underline\sigma\le\sigma_{i,m}\le1\), conditions
\eqref{eq:hetero-gbs-calibration}--\eqref{eq:hetero-compatibility}
give \eqref{eq:hetero-oracle-threshold-limit}.  In particular,
\(\psi_m^2\to\infty\), and \eqref{eq:hetero-oracle-type-two} follows.

Also
\[
\frac{\log(f_m\delta_m)}{\psi_m^2}
=
\frac{\log(f_m/r_{\alpha_m})}{\psi_m^2}
+
\frac{\log(r_{\alpha_m}\delta_m)}{\psi_m^2}
\to \frac{C_0}{2}.
\]
Thus \(f_m\delta_m\to\infty\).  To verify the Type-I term explicitly,
write \(u=u_{i,m}\) and \(c=c_{i,m}^{\mathrm{BO}}\).  From
\eqref{eq:hetero-oracle-threshold},
\[
\frac{c^2}{2}-\log(f_m\delta_m)
=
\frac12\log u
+
\frac{1}{2u}\log\{u f_m^2\delta_m^2\}
+
\frac12\left(1+\frac1u\right)\log(1+u^{-1}).
\]
The last two terms are uniformly bounded, because
\(\log(f_m\delta_m)=O(\psi_m^2)\), \(u_{i,m}\asymp\psi_m^2\), and
\(\underline\sigma\le\sigma_{i,m}\le1\).  Hence, uniformly in \(i\),
\[
f_m\delta_m\,\phi(c_{i,m}^{\mathrm{BO}})
\le C u_{i,m}^{-1/2}.
\]
Also \(c_{i,m}^{\mathrm{BO}}\asymp\sqrt{u_{i,m}}\) uniformly.  Mills'
upper bound therefore gives
\[
f_m\delta_m\,\bar\Phi(c_{i,m}^{\mathrm{BO}})
\le
\frac{f_m\delta_m\,\phi(c_{i,m}^{\mathrm{BO}})}
     {c_{i,m}^{\mathrm{BO}}}
\le C u_{i,m}^{-1}\to0
\]
uniformly, proving \eqref{eq:hetero-oracle-type-one}.  Summing the
coordinatewise Type-I and Type-II risks proves
\eqref{eq:hetero-factor-known-risk}.
\end{proof}

\subsection{Leave-one-out GLS residualization}

For each \(i\), put
\begin{equation}
\label{eq:hetero-Dminus}
D_{-i}
=
\sum_{j\ne i}\frac{\lambda_{j,m}^2}{\sigma_{j,m}^2}
\end{equation}
and define the leave-one-out GLS factor estimator
\begin{equation}
\label{eq:hetero-Zhat}
\widehat Z_{-i}
=
\frac{
\displaystyle\sum_{j\ne i}
\lambda_{j,m}X_j/\sigma_{j,m}^2
}{D_{-i}}.
\end{equation}
The corresponding feasible residual is
\begin{equation}
\label{eq:hetero-Yhat}
\widehat Y_i
=
\frac{X_i-\lambda_{i,m}\widehat Z_{-i}}{\sigma_{i,m}}
=
Y_i(Z)-\Delta_{i,m},
\end{equation}
where
\begin{equation}
\label{eq:hetero-Delta}
\Delta_{i,m}
=
\frac{\lambda_{i,m}}{\sigma_{i,m}D_{-i}}
\left[
\sum_{j\ne i}\frac{\lambda_{j,m}\theta_j}{\sigma_{j,m}^2}
+
\sum_{j\ne i}\frac{\lambda_{j,m}}{\sigma_{j,m}}\varepsilon_j
\right].
\end{equation}
By construction,
\begin{equation}
\label{eq:hetero-own-independence}
\Delta_{i,m}\perp\{Y_i(Z),\nu_i\}.
\end{equation}

\begin{lemma}[GLS perturbation bound]
\label{lem:hetero-gls-perturbation}
Assume \eqref{eq:hetero-loading-upper}--\eqref{eq:hetero-loading-information},
\[
p_m\to0,\qquad mp_m\to\infty,\qquad
\alpha_m\to0,\qquad \alpha_m f_m\to\infty,
\]
and the calibration condition \eqref{eq:hetero-gbs-calibration}.
Then \(D_{-i}\asymp m\) uniformly in \(i\).  Conditional on
\(\boldsymbol{\nu}_{-i}=(\nu_j:j\ne i)\),
\[
\Delta_{i,m}\sim N\{0,\tau_{i,m}^2(\boldsymbol{\nu})\},
\]
where
\begin{equation}
\label{eq:hetero-tau}
\tau_{i,m}^2(\boldsymbol{\nu})
=
\frac{\lambda_{i,m}^2}{\sigma_{i,m}^2}
\frac{
D_{-i}
+
\psi_m^2
\displaystyle\sum_{j\ne i}
\nu_j\lambda_{j,m}^2/\sigma_{j,m}^4
}{D_{-i}^2}.
\end{equation}
There is a constant \(C<\infty\), depending only on the loading bounds,
such that
\begin{equation}
\label{eq:hetero-tau-uniform}
\sup_{i,\boldsymbol{\nu}}\tau_{i,m}^2(\boldsymbol{\nu})
\le C\frac{1+\psi_m^2}{m}.
\end{equation}
Consequently,
\begin{equation}
\label{eq:hetero-max-delta-rate}
\max_i|\Delta_{i,m}|
=
O_{\mathbb P}\left\{
\sqrt{\frac{(1+\psi_m^2)\log m}{m}}
\right\}.
\end{equation}
\begin{equation}
\label{eq:hetero-scaled-delta}
z_{c_1(\alpha_m)}\max_i|\Delta_{i,m}|=o_{\mathbb P}(1).
\end{equation}
More precisely, for every fixed \(a>0\),
\begin{equation}
\label{eq:hetero-D-tail}
\mathbb P\left\{
z_{c_1(\alpha_m)}\max_i|\Delta_{i,m}|>a
\right\}
\le
2m\exp\left\{-c\frac{m}{\log^2m}\right\}
\end{equation}
for all sufficiently large \(m\).
\end{lemma}

\begin{proof}
By \eqref{eq:hetero-loading-information} and
\(\sigma_{j,m}^2\le1\),
\[
D_{-i}\ge
\sum_{j\ne i}\lambda_{j,m}^2
\ge cm
\]
uniformly for all sufficiently large \(m\); the reverse bound follows from
\eqref{eq:hetero-loading-upper}.  Formula \eqref{eq:hetero-tau} follows
directly from \eqref{eq:hetero-Delta}.  The same loading bounds imply
\[
\sum_{j\ne i}
\nu_j\frac{\lambda_{j,m}^2}{\sigma_{j,m}^4}
\le Cm,
\]
which proves \eqref{eq:hetero-tau-uniform}.

Conditional on \(\boldsymbol{\nu}\), a union bound and the Gaussian tail inequality give
\[
\mathbb P\{\max_i|\Delta_{i,m}|>x\mid\boldsymbol{\nu}\}
\le
2m\exp\left\{-\frac{cmx^2}{1+\psi_m^2}\right\}.
\]
The same bound holds unconditionally.  The calibration
\eqref{eq:hetero-gbs-calibration}, together with
\(p_m\to0\), \(\alpha_m\to0\), and \(\alpha_m f_m\to\infty\), implies
\[
\psi_m^2=O(\log m),
\qquad
\log(m/\alpha_m)=O(\log m),
\qquad
z_{c_1(\alpha_m)}^2=O(\log m).
\]
The asserted conclusions follow.
\end{proof}

\subsection{Proof of Theorem~\ref{thm:hetero-observed-oracle}}

\begin{proof}
The factor-known Oracle has access to the larger information set
\((\boldsymbol{X},Z)\), and therefore
\[
\mathcal R_m^{\mathrm{BO},Z}\le \mathcal R_m^{\mathrm{BO},\boldsymbol{X}}.
\]
For the reverse comparison, consider the \(\boldsymbol{X}\)-measurable plug-in rule
\[
\widehat\phi_i^{\,O}
=
I\{|\widehat Y_i|>c_{i,m}^{\mathrm{BO}}\}.
\]
The crucial point is to use the coordinatewise independence
\eqref{eq:hetero-own-independence}, rather than conditioning on a global
maximum-perturbation event.

Conditional on \(\boldsymbol{\nu}_{-i}\), under \(H_{0i}\),
\[
\widehat Y_i\sim N\{0,1+\tau_{i,m}^2(\boldsymbol{\nu})\}.
\]
Uniformly in \(i\) and \(\boldsymbol{\nu}\),
\[
(c_{i,m}^{\mathrm{BO}})^2\tau_{i,m}^2(\boldsymbol{\nu})
=
O\left(\frac{\log^2m}{m}\right)
=o(1).
\]
Hence Gaussian tail stability gives
\[
\mathbb P_0\{|\widehat Y_i|>c_{i,m}^{\mathrm{BO}}\mid\boldsymbol{\nu}_{-i}\}
=
2\bar\Phi(c_{i,m}^{\mathrm{BO}})\{1+o(1)\}
\]
uniformly.

Under \(H_{1i}\),
\[
\widehat Y_i
\sim
N\{0,1+u_{i,m}+\tau_{i,m}^2(\boldsymbol{\nu})\},
\]
and
\[
\sup_{i,\boldsymbol{\nu}}
\frac{\tau_{i,m}^2(\boldsymbol{\nu})}{1+u_{i,m}}
=O(m^{-1}).
\]
Combining this with
\eqref{eq:hetero-oracle-threshold-limit} yields
\[
\mathbb P_1\{|\widehat Y_i|\le c_{i,m}^{\mathrm{BO}}\mid\boldsymbol{\nu}_{-i}\}
=
2\Phi\{\sqrt{C_0}\sigma_{i,m}\}-1+o(1)
\]
uniformly.  Thus the plug-in rule has risk
\[
\mathcal R_m(\widehat\phi^{\,O})
=
\mathcal R_m^{\mathrm{BO},Z}+o(mp_m\delta_{A,m}).
\]
Since the genuine observed-data Oracle minimizes Bayes risk among all
\(\boldsymbol{X}\)-measurable rules,
\[
\mathcal R_m^{\mathrm{BO},\boldsymbol{X}}
\le \mathcal R_m(\widehat\phi^{\,O}).
\]
The result follows by sandwiching with
Lemma~\ref{lem:hetero-factor-known-oracle}.
\end{proof}

\subsection{Heterogeneous residualized GBS}

Define
\begin{equation}
\label{eq:hrgbs-pvalues}
\widehat P_i=2\bar\Phi(|\widehat Y_i|),
\end{equation}
and let
\(\widehat P_{(1)}\le\cdots\le\widehat P_{(m)}\) denote the ordered
residualized p-values.  The heterogeneous residualized GBS
(\(\mathrm{hrGBS}\)) procedure uses the unchanged GBS critical constants
\[
c_j(\alpha_m)
=
\frac{j\alpha_m}{m+1-j(1-\alpha_m)},
\qquad j=1,\ldots,m,
\]
and rejection count
\begin{equation}
\label{eq:hrgbs-rejection-count}
R_m^{\mathrm{hrGBS}}
=
\max\left\{
r:
\widehat P_{(j)}\le c_j(\alpha_m)
\text{ for every }j=1,\ldots,r
\right\}.
\end{equation}

For the factor-known signal p-value distribution define
\[
G_{i,m}(t)
=
2\bar\Phi\left\{
\frac{z_t}{\sqrt{1+u_{i,m}}}
\right\},
\qquad
\overline G_m(t)
=
\frac1m\sum_{i=1}^mG_{i,m}(t).
\]

\begin{lemma}[Deterministic heterogeneous signal crossing]
\label{lem:hetero-deterministic-crossing}
Assume \eqref{eq:hetero-loading-upper},
\[
p_m\to0,\qquad mp_m\to\infty,\qquad
\alpha_m\to0,\qquad \alpha_m f_m\to\infty,
\]
and \eqref{eq:hetero-gbs-calibration}.  Then, for every fixed
\(0<\varepsilon<q_-\), there exists \(\eta_\varepsilon>0\) such that
\begin{equation}
\label{eq:hetero-deterministic-crossing}
\overline G_m(q\alpha_mp_m)
\ge(1+\eta_\varepsilon)q
\end{equation}
uniformly over
\[
0<q\le Q_m-\varepsilon
\]
for all sufficiently large \(m\).
\end{lemma}

\begin{proof}
Put \(q_m^0=Q_m-\varepsilon\) and
\(t_m=q_m^0\alpha_mp_m\).  Since \(q_m^0\) is bounded away from zero,
\[
\frac{z_{t_m}^2}{\psi_m^2}\to C_0.
\]
Therefore, uniformly in \(i\),
\[
\frac{z_{t_m}}{\sqrt{1+u_{i,m}}}
=
\sqrt{C_0}\sigma_{i,m}+o(1),
\]
and hence
\[
\overline G_m(t_m)=Q_m+o(1).
\]
For each \(i\), the normal-hazard argument in the proof of
Lemma~\ref{lem:residual-deterministic-signal-crossing}, with
\(S=\sqrt{1+u_{i,m}}\), shows that \(G_{i,m}(t)/t\) is non-increasing in \(t\).
Therefore \(\overline G_m(t)/t\) is also non-increasing.  For
\(0<q\le q_m^0\),
\[
\overline G_m(q\alpha_mp_m)
\ge
\frac{q}{q_m^0}\overline G_m(t_m).
\]
Since \(Q_m/(Q_m-\varepsilon)\) is uniformly separated from one,
\eqref{eq:hetero-deterministic-crossing} follows.
\end{proof}

\begin{lemma}[Initial heterogeneous signal crossing]
\label{lem:hetero-initial-crossing}
Under the assumptions of Lemma~\ref{lem:hetero-deterministic-crossing},
\begin{equation}
\label{eq:hetero-initial-crossing}
mp_m\overline G_m\{c_1(\alpha_m)\}\to\infty.
\end{equation}
\end{lemma}

\begin{proof}
Since \(u_{i,m}\ge\psi_m^2\),
\[
\overline G_m(c_1)
\ge
2\bar\Phi\left\{
\frac{z_{c_1}}{\sqrt{1+\psi_m^2}}
\right\}.
\]
Set
\[
A_m=\log(1/\alpha_m),\qquad
B_m=\log(1/p_m),\qquad
D_m=\log(mp_m).
\]
Then
\[
\log(1/c_1)=A_m+B_m+D_m+O(1),
\qquad
\psi_m^2=\frac{2}{C_0}(A_m+B_m)\{1+o(1)\}.
\]
Writing
\[
y_m=\frac{z_{c_1}}{\sqrt{1+\psi_m^2}},
\]
we avoid any second-order Gaussian-quantile expansion.  Since
\[
z_t^2\le 2\log(1/t),\qquad 0<t<1,
\]
we have
\[
\frac{y_m^2}{2D_m}
\le
\frac{\log(1/c_1)}{(1+\psi_m^2)D_m}
\le
C\left\{
\frac1{D_m}
+
\frac1{A_m+B_m}
\right\}
+o(1)
\to0,
\]
because \(D_m=\log(mp_m)\to\infty\) and
\(A_m+B_m\to\infty\).  Hence \(y_m^2=o(D_m)\).
By Mills' lower bound,
\[
\bar\Phi(y_m)
\ge
\frac{c\,e^{-y_m^2/2}}{1+y_m}
\]
for all sufficiently large \(m\).  Therefore
\[
mp_m\bar\Phi(y_m)
=
e^{D_m}\bar\Phi(y_m)
\to\infty,
\]
which proves the claim.
\end{proof}

\begin{lemma}[Tightened factor-known heterogeneous crossing]
\label{lem:hetero-tightened-crossing}
Assume the conditions of Lemmas~\ref{lem:hetero-deterministic-crossing}
and \ref{lem:hetero-initial-crossing}.  Fix \(0<\varepsilon<q_-\) and put
\begin{equation}
\label{eq:hetero-rm}
r_m=\left\lfloor(Q_m-\varepsilon)mp_m\right\rfloor.
\end{equation}
Then there exists a fixed \(\rho_\varepsilon\in(0,1)\) such that
\begin{equation}
\label{eq:hetero-tightened-crossing}
\mathbb P\left\{
\widetilde P^1_{(j)}
\le(1-\rho_\varepsilon)c_j(\alpha_m),
\quad j=1,\ldots,r_m
\right\}\to1,
\end{equation}
where \(\widetilde P^1_{(j)}\) are the ordered factor-known signal
p-values.
\end{lemma}

\begin{proof}
Fix a sufficiently small \(\rho\in(0,1)\) and put
\[
\alpha_m^{(1)}=(1-\rho)\alpha_m,
\qquad
\alpha_m^{(2)}=(1-\rho)^2\alpha_m.
\]
Because \(\rho\) is fixed,
\[
\log\{r_{\alpha_m^{(\ell)}}/r_{\alpha_m}\}=O(1),
\qquad \ell=1,2,
\]
and \(\psi_m^2\to\infty\); hence both modified levels retain the
calibration constant \(C_0\).

Let
\[
d_j=(1-\rho)c_j(\alpha_m),
\qquad
N_{1,m}(t)
=
\sum_{i=1}^m
I\{\nu_i=1,\ \widetilde P_i\le t\}.
\]
The summands in \(N_{1,m}(t)\) are independent Bernoulli random variables,
not necessarily identically distributed, with success probabilities
\[
p_mG_{i,m}(t).
\]
Consequently,
\begin{equation}
\label{eq:hetero-N1-unconditional-mean}
\mathbb E N_{1,m}(t)
=
mp_m\overline G_m(t).
\end{equation}
This unconditional representation is important: no conditioning on the
realized signal set is needed.

For the bulk ranks, choose \(L_m\to\infty\) below and consider
\(L_m<j\le r_m\).  Since \(L_m\to\infty\), eventually \(j\ge2\), and
\[
c_j(\alpha_m)\ge\frac{j\alpha_m}{m}.
\]
Writing \(q_j=j/(mp_m)\), we therefore have
\[
d_j
\ge
q_j\alpha_m^{(1)}p_m.
\]
Because \(j\le r_m\), \(q_j\le Q_m-\varepsilon+o(1)\).  By choosing
\(\rho>0\) sufficiently small and, if necessary, replacing
\(\varepsilon\) by a fixed smaller positive fraction of itself in the
deterministic-crossing lemma, there is a constant
\(\eta_\varepsilon>0\) such that
\[
\overline G_m(d_j)
\ge
(1+\eta_\varepsilon)q_j
\]
uniformly over \(L_m<j\le r_m\).  Hence, by
\eqref{eq:hetero-N1-unconditional-mean},
\[
\mathbb E N_{1,m}(d_j)
\ge
(1+\eta_\varepsilon)j.
\]
The multiplicative Chernoff inequality for a Poisson-binomial sum now gives
\begin{equation}
\label{eq:hetero-bulk-chernoff}
\mathbb P\{N_{1,m}(d_j)<j\}
\le
\exp(-C_\varepsilon j)
\end{equation}
uniformly over the bulk ranks.

For the earliest ranks,
\[
c_1(\alpha_m^{(2)})
<
(1-\rho)c_1(\alpha_m)=d_1
\]
eventually.  Lemma~\ref{lem:hetero-initial-crossing}, applied at
\(\alpha_m^{(2)}\), yields
\[
M_m:=mp_m\overline G_m(d_1)\to\infty.
\]
Choose \(L_m\to\infty\) so that
\[
L_m=o(mp_m),
\qquad
L_m=o(M_m).
\]
Since \(N_{1,m}(d_1)\) is again Poisson-binomial with mean \(M_m\),
a multiplicative Chernoff bound gives
\[
\mathbb P\{N_{1,m}(d_1)<L_m\}
\le
\exp(-cM_m)\to0.
\]
On the complementary event, monotonicity of \(d_j\) implies
\[
N_{1,m}(d_j)\ge N_{1,m}(d_1)\ge L_m\ge j,
\qquad 1\le j\le L_m.
\]
For the bulk ranks, summing \eqref{eq:hetero-bulk-chernoff} gives
\[
\sum_{j=L_m+1}^{r_m}
\mathbb P\{N_{1,m}(d_j)<j\}
\le
\sum_{j=L_m+1}^{\infty}e^{-C_\varepsilon j}
\to0.
\]
Thus, with probability tending to one,
\(N_{1,m}(d_j)\ge j\) for every \(1\le j\le r_m\), which is equivalent
to \eqref{eq:hetero-tightened-crossing}.  Taking the fixed sufficiently
small value of \(\rho\) as \(\rho_\varepsilon\) completes the proof.
\end{proof}

\begin{lemma}[Feasible heterogeneous signal crossing and Type-II bound]
\label{lem:hetero-feasible-crossing}
Assume the conditions of Lemmas~\ref{lem:hetero-tightened-crossing}
and \ref{lem:hetero-gls-perturbation}.  For \(r_m\) in \eqref{eq:hetero-rm},
\begin{equation}
\label{eq:hetero-feasible-crossing}
\mathbb P\left\{
\widehat P^1_{(j)}\le c_j(\alpha_m),
\quad j=1,\ldots,r_m
\right\}\to1.
\end{equation}
Consequently,
\begin{equation}
\label{eq:hetero-type-two-bound}
\mathbb E(T_m^{\mathrm{hrGBS}})
\le
mp_m\{1-Q_m+\varepsilon\}+o(mp_m).
\end{equation}
\end{lemma}

\begin{proof}
Let \(\rho=\rho_\varepsilon\) be supplied by
Lemma~\ref{lem:hetero-tightened-crossing}.  Choose
\[
0<a_\rho<-\log(1-\rho).
\]
Mills' ratio gives, uniformly for sufficiently small \(t\),
\[
z_{(1-\rho)t}\ge z_t+\frac{a_\rho}{z_t}.
\]
By \eqref{eq:hetero-scaled-delta},
\[
z_{c_1}\max_i|\Delta_{i,m}|=o_{\mathbb P}(1).
\]
Hence, with probability tending to one,
\[
\widetilde P_i\le(1-\rho)c_j
\quad\Longrightarrow\quad
\widehat P_i\le c_j
\]
uniformly for \(j\le r_m\).  Combining this implication with
Lemma~\ref{lem:hetero-tightened-crossing} proves
\eqref{eq:hetero-feasible-crossing}.

Let \(A_m\) denote the event in \eqref{eq:hetero-feasible-crossing}.
Exactly as in the homogeneous proof, on \(A_m\) the first \(r_m\)
ordered signal p-values themselves are rejected, so
\(S_m^{\mathrm{hrGBS}}\ge r_m\).  Therefore, pathwise,
\[
T_m^{\mathrm{hrGBS}}
\le
(K_m-r_m)^+ + K_m I(A_m^c).
\]
Since \(K_m\sim\operatorname{Bin}(m,p_m)\) and \(mp_m\to\infty\),
\[
\mathbb E|K_m-mp_m|=o(mp_m).
\]
Therefore
\[
\mathbb E\{(K_m-r_m)^+\}
\le
mp_m-r_m+\mathbb E|K_m-mp_m|
=
mp_m\{1-Q_m+\varepsilon\}+o(mp_m).
\]
Moreover, by Cauchy--Schwarz,
\[
\mathbb E\{K_m I(A_m^c)\}
\le
\{\mathbb E(K_m^2)\}^{1/2}\mathbb P(A_m^c)^{1/2}.
\]
Since
\[
\mathbb E(K_m^2)
=
mp_m(1-p_m)+(mp_m)^2
=
O\{(mp_m)^2\}
\]
and \(\mathbb P(A_m^c)\to0\), the latter term is \(o(mp_m)\).
Combining these bounds proves \eqref{eq:hetero-type-two-bound}.
\end{proof}

\begin{lemma}[Heterogeneous Type-I risk]
\label{lem:hetero-type-one-risk}
Assume the conditions of Lemma~\ref{lem:hetero-gls-perturbation}.  Assume in addition
\begin{equation}
\label{eq:hetero-type-one-loss}
\alpha_m\delta_m\to0.
\end{equation}
Then
\begin{equation}
\label{eq:hetero-type-one-risk}
\delta_m\mathbb E(V_m^{\mathrm{hrGBS}})=o(mp_m).
\end{equation}
\end{lemma}

\begin{proof}
Fix \(\rho>0\), put
\[
\alpha_m'=(1+\rho)\alpha_m,
\qquad
r_m^\star=\lfloor m/2\rfloor,
\]
and let \(\widetilde V_m'\) and \(\widetilde R_m'\) denote the number of
false discoveries and the rejection count of factor-known GBS at level
\(\alpha_m'\).  The deterministic comparison of the GBS critical
constants used in Lemma~\ref{lem:inflated-level-coupling} is unchanged.
Thus there is a fixed \(a_\rho>0\) such that, on
\[
D_m^I=
\left\{
z_{c_1(\alpha_m)}
\max_i|\Delta_{i,m}|\le a_\rho
\right\},
\]
every hrGBS crossing up to rank \(r_m^\star\) implies the corresponding
factor-known GBS crossing at level \(\alpha_m'\).  Consequently,
\[
V_m^{\mathrm{hrGBS}}
\le
\widetilde V_m'
+
mI\{(D_m^I)^c\}
+
mI\{\widetilde R_m'\ge r_m^\star\}.
\]

Conditional on the signal configuration and signal p-values, the
factor-known null p-values are iid uniform.  Hence the shifted
order-statistic inequality, Lemma~\ref{lem:shifted-order-statistic},
gives
\[
\mathbb E(\widetilde V_m'\mid K_m=k)
\le
\frac{\alpha_m'}{1-\alpha_m'}(k+1),
\]
and therefore
\[
\mathbb E(\widetilde V_m')=O(\alpha_mmp_m).
\]
By \eqref{eq:hetero-D-tail},
\[
\mathbb P\{(D_m^I)^c\}
\le
2m\exp\left\{-c\frac{m}{\log^2m}\right\}.
\]
Also
\[
\widetilde R_m'\le K_m+\widetilde V_m'.
\]
With \(s_m^\star=\lfloor m/4\rfloor\),
\[
\{\widetilde R_m'\ge r_m^\star\}
\subseteq
\{K_m\ge s_m^\star\}
\cup
\{\widetilde V_m'\ge r_m^\star-s_m^\star\}.
\]
Since \(p_m\to0\), the binomial Chernoff bound gives, for some \(c>0\),
\[
\mathbb P(K_m\ge s_m^\star)
\le
\exp\{-cm\log(1/p_m)\}.
\]
By Markov's inequality and the preceding bound on
\(\mathbb E(\widetilde V_m')\),
\[
m\mathbb P\{\widetilde V_m'\ge r_m^\star-s_m^\star\}
=
O(\alpha_mmp_m).
\]
Hence
\[
m\mathbb P(\widetilde R_m'\ge r_m^\star)
\le
m e^{-cm\log(1/p_m)}
+
O(\alpha_mmp_m).
\]

Finally, \(mp_m\to\infty\) and \(\alpha_mf_m\to\infty\) imply eventually
\(\alpha_m>p_m>1/m\).  Together with
\(\alpha_m\delta_m\to0\), this yields
\(\delta_m=o(1/\alpha_m)=o(m)\).  Therefore
\[
\delta_m m\mathbb P\{(D_m^I)^c\}=o(mp_m),
\]
by \eqref{eq:hetero-D-tail}, and also
\[
\delta_m m e^{-cm\log(1/p_m)}=o(mp_m).
\]
Finally,
\[
\delta_m\,O(\alpha_mmp_m)=o(mp_m)
\]
because \(\alpha_m\delta_m\to0\).  Multiplying the pathwise inequality
by \(\delta_m\), taking expectations, and combining these three bounds
proves \eqref{eq:hetero-type-one-risk}.
\end{proof}

For clarity, retain \(R_m^{\mathrm{hrGBS}}\) for the rejection count in
\eqref{eq:hrgbs-rejection-count}, and define the Bayes risk of hrGBS by
\begin{equation}
\label{eq:hrgbs-bayes-risk-definition}
\mathcal R_m^{\mathrm{hrGBS}}
=
\delta_{0,m}\mathbb E(V_m^{\mathrm{hrGBS}})
+
\delta_{A,m}\mathbb E(T_m^{\mathrm{hrGBS}}).
\end{equation}

\subsection{Proof of Theorem~\ref{thm:hrgbs-abos}}

\begin{proof}
Fix \(0<\varepsilon<q_-\).  By
Lemma~\ref{lem:hetero-feasible-crossing},
\[
\mathbb E(T_m^{\mathrm{hrGBS}})
\le
mp_m\{1-Q_m+\varepsilon\}+o(mp_m),
\]
whereas Lemma~\ref{lem:hetero-type-one-risk} gives
\[
\delta_{0,m}\mathbb E(V_m^{\mathrm{hrGBS}})
=o(mp_m\delta_{A,m}).
\]
Therefore
\[
\mathcal R_m^{\mathrm{hrGBS}}
\le
mp_m\delta_{A,m}\{1-Q_m+\varepsilon+o(1)\}.
\]
By Theorem~\ref{thm:hetero-observed-oracle},
\[
\mathcal R_m^{\mathrm{BO},\boldsymbol{X}}
=
mp_m\delta_{A,m}\{1-Q_m+o(1)\}.
\]
Since hrGBS is \(\boldsymbol{X}\)-measurable,
\[
\mathcal R_m^{\mathrm{hrGBS}}\ge \mathcal R_m^{\mathrm{BO},\boldsymbol{X}}.
\]
Finally, \eqref{eq:Qm-bounds} implies
\(1-Q_m\ge1-q_+>0\).  Hence
\[
1
\le
\liminf_m\frac{\mathcal R_m^{\mathrm{hrGBS}}}{\mathcal R_m^{\mathrm{BO},\boldsymbol{X}}}
\le
\limsup_m\frac{\mathcal R_m^{\mathrm{hrGBS}}}{\mathcal R_m^{\mathrm{BO},\boldsymbol{X}}}
\le
1+\frac{\varepsilon}{1-q_+}.
\]
Letting \(\varepsilon\downarrow0\) proves
\eqref{eq:hrgbs-abos-ratio}.  Equation
\eqref{eq:hrgbs-risk-expansion} follows from the definition of \(Q_m\).
\end{proof}

\begin{remark}[Recovery of the homogeneous theory]
\label{rem:hetero-recovers-homogeneous}
If \(\lambda_{i,m}\equiv\lambda\), then
\[
Q_m
=
2\bar\Phi\{\sqrt{C_0(1-\lambda^2)}\}.
\]
Writing \(C_\lambda=C_0(1-\lambda^2)\), the risk expansion in
Theorem~\ref{thm:hrgbs-abos} reduces to
\[
\mathcal R_m^{\mathrm{hrGBS}}
=
mp_m\delta_{A,m}
\{2\Phi(\sqrt{C_\lambda})-1\}\{1+o(1)\},
\]
which is exactly the homogeneous rGBS result.
\end{remark}

\section{Proofs for bounded counts and polynomial sparsity}
\label{sec:bounded-polynomial}

The preceding empirical-crossing arguments require a diverging expected
number of signals. We now treat bounded positive expected counts by using
the first GBS critical constant directly. Throughout this subsection the
loss ratio is fixed; these results do not assert coverage of
\(mp_m\to0\).

\subsection{Proof of Theorem~\ref{thm:hrgbs-bounded}}

\begin{proof}
We give the bounds explicitly, without invoking the diverging-count
assumptions of the preceding lemmas. In this proof all uniform error terms
are uniform over the coordinates and, when indicated, over latent
configurations. Write \(u_i=u_{i,m}\), \(\sigma_i=\sigma_{i,m}\), and
\(c_i=c_{i,m}^{\mathrm{BO}}\).

\medskip\noindent\emph{Observed-data Oracle.}
The two-sided bounds on \(mp_m\) imply
\(\log f_m=\log m+O(1)\). Thus
\[
u_i\asymp\log m,\qquad
\frac{\log(u_i f_m^2\delta^2)}{u_i}
=C_0\sigma_i^2+o(1)
\]
uniformly. The exact threshold formula
\eqref{eq:hetero-oracle-threshold} therefore gives
\[
\frac{c_i^2}{1+u_i}=C_0\sigma_i^2+o(1),\qquad
c_i\asymp\sqrt{\log m}.
\]
The same formula and Mills' upper bound give, uniformly,
\[
f_m\delta\,2\bar\Phi(c_i)=O(u_i^{-1})=o(1).
\]
Indeed, writing \(L_i=\log\{(1+u_i)f_m^2\delta^2\}\),
\[
f_m\delta\exp(-c_i^2/2)
=(1+u_i)^{-1/2}\exp\{-L_i/(2u_i)\}
=O(u_i^{-1/2}),
\]
and division by \(c_i\asymp\sqrt{u_i}\) proves the claim.
Consequently the factor-known Oracle risk equals
\[
\mathcal R_m^{\mathrm{BO},Z}
=mp_m\delta_{A,m}\{1-Q_m+o(1)\}.
\]

The loading bounds imply \(D_{-i}\asymp m\). The exact conditional
variance identity \eqref{eq:hetero-tau}, which does not require
\(mp_m\to\infty\), gives
\[
\sup_{i,\boldsymbol{\nu}}\tau_{i,m}^2(\boldsymbol{\nu})
\le C\frac{1+\psi_m^2}{m}
=O\left(\frac{\log m}{m}\right).
\]
Under the null, conditional on \(\boldsymbol{\nu}_{-i}\), the feasible residual is
\(N\{0,1+\tau_{i,m}^2(\boldsymbol{\nu})\}\), independently constructed from its
own latent state. Since
\[
\sup_{i,\boldsymbol{\nu}}c_i^2\tau_{i,m}^2(\boldsymbol{\nu})
=O(\log^2m/m)=o(1),
\]
Mills' bounds yield
\[
\mathbb P\{|\widehat Y_i|>c_i\mid\nu_i=0,\boldsymbol{\nu}_{-i}\}
=2\bar\Phi(c_i)\{1+o(1)\}
\]
uniformly. Under the alternative the conditional variance is
\(1+u_i+\tau_{i,m}^2(\boldsymbol{\nu})\), so
\[
\mathbb P\{|\widehat Y_i|\le c_i\mid\nu_i=1,\boldsymbol{\nu}_{-i}\}
=2\Phi(\sqrt{C_0}\sigma_i)-1+o(1)
\]
uniformly as well. The plug-in rule using thresholds \(c_i\) therefore
has risk \(mp_m\delta_{A,m}\{1-Q_m+o(1)\}\).
Sandwiching the observed-data Oracle risk between the factor-known risk
and this feasible plug-in risk proves \eqref{eq:bounded-hetero-oracle}.

\medskip\noindent\emph{Type-II risk of hrGBS.}
Let \(c_{1,m}=\alpha_m/(m+\alpha_m)\) and \(z_m=z_{c_{1,m}}\).
Every coordinate with \(\widehat P_i\le c_{1,m}\) is rejected: if it has
overall rank \(j\), then for every \(\ell\le j\),
\[
\widehat P_{(\ell)}\le\widehat P_i\le c_{1,m}\le c_\ell(\alpha_m).
\]
Thus, without any independence assumption on the feasible p-values,
\[
\mathbb E(T_m^{\mathrm{hrGBS}})
\le p_m\sum_{i=1}^m
\mathbb P\{|\widehat Y_i|<z_m\mid\nu_i=1\}.
\]
The level assumption gives \(z_m^2=2\log m\{1+o(1)\}\). Hence,
uniformly in \(i,\boldsymbol{\nu}\),
\[
\frac{z_m^2}{1+u_i+\tau_{i,m}^2(\boldsymbol{\nu})}
=C_0\sigma_i^2+o(1).
\]
Using the conditional Gaussian law and averaging over \(\boldsymbol{\nu}_{-i}\),
\begin{equation}
\label{eq:bounded-hetero-type-two}
\mathbb E(T_m^{\mathrm{hrGBS}})
\le mp_m\{1-Q_m+o(1)\}.
\end{equation}

\medskip\noindent\emph{Type-I risk of hrGBS.}
Fix \(\rho>0\), set \(\alpha_m'=(1+\rho)\alpha_m\), and put
\(r_m^\star=\lfloor m/2\rfloor\), \(s_m^\star=\lfloor m/4\rfloor\).
The deterministic coupling proof in
Lemma~\ref{lem:inflated-level-coupling} only uses \(\alpha_m\to0\)
and a bound on the maximum residual perturbation; it does not use
concentration of the signal count. It therefore gives, for a fixed
\(a_\rho>0\),
\[
V_m^{\mathrm{hrGBS}}\le\widetilde V_m'
+mI\{(D_m^I)^c\}+mI\{\widetilde R_m'\ge r_m^\star\},
\qquad
D_m^I=\{z_m\max_i|\Delta_{i,m}|\le a_\rho\}.
\]
The conditional variance bound above and a Gaussian union bound imply
\[
\mathbb P\{(D_m^I)^c\}
\le2m\exp\{-c m/\log^2m\},\qquad
m\mathbb P\{(D_m^I)^c\}=o(1).
\]
For the factor-known procedure, null p-values remain independent uniform
conditional on the signal configuration and its p-values. The shifted
order-statistic bound therefore yields
\[
\mathbb E(\widetilde V_m')
\le\frac{\alpha_m'}{1-\alpha_m'}(mp_m+1)=O(\alpha_m).
\]
Moreover, \(\widetilde R_m'\le K_m+\widetilde V_m'\), so
\[
m\mathbb P(\widetilde R_m'\ge r_m^\star)
\le m\mathbb P(K_m\ge s_m^\star)
+m\mathbb P(\widetilde V_m'\ge r_m^\star-s_m^\star).
\]
Since \(p_m\le b_+/m\), the binomial upper-tail bound gives
\(\mathbb P(K_m\ge s_m^\star)\le e^{-c m\log m}\) eventually.
Markov's inequality bounds the second term by
\(C\mathbb E(\widetilde V_m')=O(\alpha_m)\). Taking expectations in the
pathwise inequality proves
\[
\delta\,\mathbb E(V_m^{\mathrm{hrGBS}})=o(1)=o(mp_m).
\]
Together with \eqref{eq:bounded-hetero-type-two}, this gives an upper
bound \(mp_m\delta_{A,m}\{1-Q_m+o(1)\}\) on hrGBS risk. The
observed-data Oracle is a lower bound. Finally \(1-Q_m\ge1-q_+>0\),
so the additive risk remainders yield the asserted risk ratio.
\end{proof}

\subsection{Proof of Corollary~\ref{cor:rgbs-bounded}}

\begin{proof}
Set \(\lambda_{i,m}\equiv\lambda\). Then the GLS factor estimate is
exactly \(\{\lambda(m-1)\}^{-1}\sum_{j\ne i}X_j\), so hrGBS coincides
with homogeneous rGBS. Since \(\psi_m^2=(1-\lambda^2)u_{m,\lambda}\),
Theorem~\ref{thm:hrgbs-bounded} applies with
\(C_0=C_\lambda/(1-\lambda^2)\) and
\(1-Q_m=2\Phi(\sqrt{C_\lambda})-1\).
\end{proof}

\subsection{Proof of Corollary~\ref{cor:onefactor-polynomial}}

\begin{proof}
For \(0<\beta<1\), \(mp_m\to\infty\) and
\(\log f_m=\beta\log m+O(1)\). The level condition implies
\[
\log(f_m/r_{\alpha_m})=\beta\log m+o(\log m),\qquad
\alpha_m f_m=m^{\beta-o(1)}\to\infty.
\]
Also \(\log\{1/(r_{\alpha_m}\delta)\}=o(\log m)\) and
\(\alpha_m\delta\to0\). Thus the heterogeneous and homogeneous
calibrations and compatibility conditions follow from the respective
signal scalings. Theorems~\ref{thm:hrgbs-abos} and \ref{thm:rgbs-abos}
apply. For \(\beta=1\), \(mp_m\) is bounded above and bounded away
from zero. Theorem~\ref{thm:hrgbs-bounded} and
Corollary~\ref{cor:rgbs-bounded} apply directly, including when
\(mp_m\) does not converge.
\end{proof}

\section{Residual BH: expected-error bounds and proof}
\label{app:rbh}

\begin{lemma}[Expected false counts for independent BH at a small level]
\label{lem:bh-false-count}
Suppose there are $k$ non-null p-values and $n=m-k$ independent uniform
null p-values, independent of the non-null p-values. The non-null p-values
may be arbitrary. For BH at level $0<\beta\le1/(16e)$,
\[
\mathbb E(V\mid K=k)\le40\beta(k+1).
\]
Consequently independent factor-known BH satisfies
$\mathbb E\widetilde V(\beta)\le40\beta(mp_m+1)$.
\end{lemma}
\begin{proof}
If $n=0$, the assertion is immediate. Otherwise let $U_{(\ell)}$ be the
ordered null p-values. If BH makes exactly $v\ge1$ false discoveries,
its total number of rejections is at most $k+v$, and
$U_{(v)}\le\beta(k+v)/m$. Thus, for $j\ge1$, a union bound gives
\[
\Pr(V\ge j\mid K=k)
\le\sum_{\ell=j}^{n}\Pr\{U_{(\ell)}\le\beta(k+\ell)/m\}
\le\sum_{\ell=j}^{n}\{e\beta(1+k/\ell)\}^{\ell}.
\]
The last inequality follows by a union bound over the $\binom n\ell$
subsets of nulls and $\binom n\ell\le(en/\ell)^\ell$;
$\beta(k+\ell)/m\le\beta<1$.
Put $x=\beta(k+1)$. If $x\le1/(4e)$, all bases are at most $ex\le1/4$,
and summing the preceding display over $j$ gives
\[
\mathbb E(V\mid K=k)\le\sum_{\ell\ge1}\ell(ex)^\ell
=\frac{ex}{(1-ex)^2}\le2ex.
\]
If $x>1/(4e)$, set $j_0=\lceil4ex\rceil$. For $\ell\ge j_0$,
$e\beta(1+k/\ell)\le1/16+1/4=5/16$. Therefore
\[
\mathbb E(V\mid K=k)
\le j_0-1+\sum_{\ell\ge j_0}\ell(5/16)^\ell
\le4ex+1\le8ex<40x.
\]
Here $\sum_{\ell\ge1}\ell(5/16)^\ell=80/121<1$.
The bound remains valid if $j_0>n$. It holds conditional on any non-null
configuration and its p-values, so averaging over them and then over
$K_m$ proves the unconditional assertion.
\end{proof}

\begin{lemma}[Fixed-residual perturbation and BH inclusion]
\label{lem:rbh-inclusion}
Under either case of Theorem~\ref{thm:rbh-abos}, there is a fixed $a>0$
such that, on
$\mathcal D_m=\{z_{\alpha_m/m}\max_i|\Delta_{i,m}|\le a\}$,
every rBH rejection belongs to the rejection set of factor-known BH at
level $2\alpha_m$, eventually. Moreover,
\begin{equation}\label{eq:rbh-perturb-tail}
\Pr(\mathcal D_m^c)\le2m\exp\{-cm/(\log m)^2\},\qquad
\frac{1+\delta_m}{p_m}\Pr(\mathcal D_m^c)\to0.
\end{equation}
On the same event one may choose $a$ small enough that, for all
$1\le j\le m$,
\[
\widetilde P_i\le\alpha_mj/(2m)
\quad\Longrightarrow\quad\widehat P_i\le\alpha_mj/m.
\]
\end{lemma}
\begin{proof}
Use the exact conditional variance identity \eqref{eq:hetero-tau} and
$D_{-i}\asymp m$. Under either set of assumptions,
$\psi_m^2=O(\log m)$, $p_m\ge b/m$ and eventually $\alpha_m\ge p_m$.
Thus $z_{\alpha_m/m}^2=O(\log m)$, and the conditional Gaussian union
bound gives \eqref{eq:rbh-perturb-tail}'s probability estimate. This
argument does not require diverging signal counts. In case (i),
$\delta_m\le1/\alpha_m$ eventually; in case (ii) it is fixed.
In both cases $(1+\delta_m)/p_m=O(m^2)$, which proves the weighted limit.

The Mills-ratio transfers used in Lemma~\ref{mrd:lem:transfers} apply
to $t_j=j\alpha_m/m$ uniformly through $j=m$: the largest such threshold
is $\alpha_m\to0$, and $z_{\alpha_m/m}\ge z_{t_j}$.
Choose $a>0$ sufficiently small. The two deterministic implications are
\[
\widehat P_i\le t_j\ \Longrightarrow\ \widetilde P_i\le2t_j,
\qquad \widetilde P_i\le t_j/2\ \Longrightarrow\ \widehat P_i\le t_j.
\]
These comparisons use no independence of the feasible p-values.
If rBH makes $r>0$ rejections, each selected index has
$\widehat P_i\le r\alpha_m/m$, so at least $r$ factor-known p-values
are at most $2r\alpha_m/m$. The factor-known BH count $\widetilde R'$
is therefore at least $r$. Every selected index also satisfies
$\widetilde P_i\le2\alpha_m\widetilde R'/m$, which is precisely the
factor-known BH rejection threshold. (BH rejects all p-values below its
final threshold.) This proves set inclusion, including ties. The case
$r=0$ is trivial.
\end{proof}

\subsection{Proof of Theorem~\ref{thm:rbh-abos}}
\begin{proof}
The two preceding lemmas yield the pathwise and expectation bounds
\[
V_m^B\le\widetilde V(2\alpha_m)+mI(\mathcal D_m^c),
\]
\begin{align*}
\frac{(1+\delta_m)\mathbb EV_m^B}{mp_m}
&\le80\alpha_m(1+\delta_m)\left(1+\frac1{mp_m}\right)\\
&\quad+\frac{1+\delta_m}{p_m}\Pr(\mathcal D_m^c)\to0.
\end{align*}
The small-level restriction in Lemma~\ref{lem:bh-false-count} holds
because $2\alpha_m\to0$. This proves both assertions in
\eqref{eq:rbh-typeI} without a deletion cap or an extra sparsity rate.

\medskip\noindent\emph{Diverging expected counts.}
Fix $0<\epsilon<q_-$ and put $r_m=\lfloor(Q_m-\epsilon)mp_m\rfloor$.
At $t_m=\alpha_m r_m/(2m)$ the calibration gives, uniformly in $i$,
\[
G_{i,m}(t_m)=
2\bar\Phi\left\{\frac{z_{t_m}}{\sqrt{1+\psi_m^2/\sigma_i^2}}\right\}
=2\bar\Phi(\sqrt{C_0}\sigma_i)+o(1).
\]
Indeed $\log(1/t_m)=\log(f_m/r_{\alpha_m})+O(1)$ and
$\psi_m^2\to\infty$. The number $N_m(t_m)$ of signal p-values below
$t_m$ is a sum of independent, possibly nonidentical Bernoulli variables
with mean $mp_m\{Q_m+o(1)\}$. Its mean exceeds $r_m$ by a fixed fraction
of $mp_m$, so a Chernoff bound gives $\Pr\{N_m(t_m)<r_m\}\to0$.
On the complementary event intersected with $\mathcal D_m$, at least
$r_m$ feasible p-values are at most $\alpha_mr_m/m$.
The BH step-up condition therefore gives $R_m^B\ge r_m$.
Consequently
\[
T_m^B=K_m-R_m^B+V_m^B
\le K_m-r_m+V_m^B+r_mI\{R_m^B<r_m\},
\]
and taking expectations proves
$\mathbb ET_m^B\le mp_m(1-Q_m+\epsilon)+o(mp_m)$.
This proof only needs a single deterministic crossing rank, unlike the
all-prefix crossing required for a step-down rule.

\medskip\noindent\emph{Bounded positive expected counts.}
Call a signal strong when $\widetilde P_i\le\alpha_m/(2m)$.
On $\mathcal D_m$, every strong signal has $\widehat P_i\le\alpha_m/m$
and is rejected by BH. If $N_m^{\rm strong}$ counts these signals,
\[
T_m^B\le K_m-N_m^{\rm strong}+K_mI(\mathcal D_m^c).
\]
Since $\mathbb EK_m^2=O((mp_m)^2)$, Cauchy--Schwarz and the vanishing
exceptional probability make the final expectation $o(mp_m)$.
Further $z_{\alpha_m/(2m)}^2\sim2\log m$ and
$\psi_m^2/\log m\to2/C_0$, giving
$\mathbb EN_m^{\rm strong}=mp_m\{Q_m+o(1)\}$ uniformly.
Thus $\mathbb ET_m^B\le mp_m\{1-Q_m+o(1)\}$.

\medskip\noindent\emph{Risk normalization.}
The factor-known Oracle expansion in Lemma~\ref{mrd:lem:oracle} uses
only calibration and loading bounds;
its displayed proof therefore applies in both present cases. Alternatively
use Theorems~\ref{thm:hetero-observed-oracle} and
\ref{thm:hrgbs-bounded} for the observed-data benchmark.
Hence
\[
mp_m\delta_{A,m}\{1-Q_m+o(1)\}
=\mathcal R_m^{\mathrm{BO},Z}
\le\mathcal R_m^{\mathrm{BO},\boldsymbol{X}}\le\mathcal R_m^B.
\]
Combine this lower bound with negligible weighted Type-I error and the
Type-II upper bounds; let $\epsilon\downarrow0$ in the divergent case.
Since $1-Q_m\ge1-q_+>0$, the risk ratio tends to one even if $Q_m$
does not converge. Subtracting the negligible weighted Type-I term also
gives equality to first order in \eqref{eq:rbh-typeII}.
The polynomial specialization follows by the same calibration check as
Corollary~\ref{cor:onefactor-polynomial}.
\end{proof}

\section{Known one-factor covariance and general vanishing sparsity}
\label{app:covariance-sparsity}

\subsection{Proof of Lemma~\ref{lem:known-sigma-loadings}}
\begin{proof}
Suppress $m$. Let $a,b,c$ be distinct indices with nonzero loadings.
The off-diagonal entries satisfy $\Sigma_{ij}=\lambda_i\lambda_j$.
Consequently
\begin{equation}\label{eq:loading-reconstruction}
\lambda_a^2=\frac{\Sigma_{ab}\Sigma_{ac}}{\Sigma_{bc}}>0,
\qquad
\lambda_j=\frac{\Sigma_{aj}}{\lambda_a},\quad j\ne a.
\end{equation}
Choosing either square root for $\lambda_a$ determines every coordinate;
zero loadings are recovered as zero. Any other loading representation
of the same normalized covariance must satisfy these identities and
therefore differs only by the common sign. Thus
$\sigma_j^2=1-\lambda_j^2$ is uniquely determined.
A valid triple can be selected from the known covariance by requiring
its three off-diagonal entries to be nonzero, using a fixed ordering
if more than one triple exists.

The information condition gives $\sum_i\lambda_{i,m}^2\ge\kappa m/2$
eventually. Since each squared loading is at most $\bar\lambda^2$,
the number of nonzero loadings is at least
$\kappa m/(2\bar\lambda^2)$, which exceeds three eventually.
The case $\bar\lambda=0$ is incompatible with $\kappa>0$.
Hence the reconstruction applies for all sufficiently large dimensions.

The covariance matrices and all MRD residuals are unchanged by
$\lambda\mapsto-\lambda$. For the fixed GLS construction,
$D_{-i}=\sum_{j\ne i}\lambda_j^2/\sigma_j^2$ and $\sigma_i$ are
unchanged, whereas $\widehat Z_{-i}$ changes sign. The product
$\lambda_i\widehat Z_{-i}$, and thus every $\widehat Y_i$, is unchanged.
Equivalently, $\lambda_i\lambda_j=\Sigma_{ij}$ yields
\eqref{eq:gls-known-sigma} directly. Therefore both fixed residualized
GBS and residual BH have identical p-values and decisions for either
orientation. The observed Gaussian experiment and its Bayes Oracle
are determined by $\boldsymbol{\Sigma}$, the prior and losses; the factor realization
has not been supplied. This proves the equivalence used in the theorems.
\end{proof}

\subsection{Proof of Corollary~\ref{cor:general-vanishing-known-sigma}}
\begin{proof}
Put $B_m=\log(1/p_m)$ and $A_m=\log(1/\alpha_m)$.
The assumptions give $B_m\to\infty$, $A_m=o(B_m)$,
$\psi_m^2\sim2B_m/C_0$, $\log f_m=B_m+o(1)$ and
$\log(1/r_{\alpha_m})=A_m+o(1)$. Therefore
\[
\alpha_mf_m=\exp\{B_m-A_m+o(1)\}\to\infty,
\qquad
\frac{2\log(f_m/r_{\alpha_m})}{\psi_m^2}\to C_0,
\]
and, for fixed positive $\delta$,
\[
\frac{\log\{1/(r_{\alpha_m}\delta)\}}{\psi_m^2}\to0,
\qquad \alpha_m\delta\to0.
\]
If $mp_m\to\infty$ along a subsequence, all hypotheses of
Theorems~\ref{mrd:thm:diverging}, \ref{thm:hrgbs-abos} and
\ref{thm:rbh-abos}(i) hold along that subsequence, proving the desired
risk expansions and ratios there.

If $mp_m$ is bounded along a subsequence, the positive lower bound on
its liminf gives $b_-\le mp_m\le b_+$ along its tail, for suitable
constants. Thus $B_m=\log m+O(1)$, so
$\psi_m^2/\log m\to2/C_0$ and
$\log(1/\alpha_m)=o(\log m)$. Theorems~\ref{mrd:thm:bounded},
\ref{thm:hrgbs-bounded} and \ref{thm:rbh-abos}(ii) give the same
conclusions along that subsequence. The observed-data Oracle expansion
follows either from those theorems or from the procedure expansion and
its ratio to the Oracle.

To conclude for the full sequence, take an arbitrary subsequence. If
it has infinitely many expected counts below some finite constant,
extract a bounded further subsequence. Otherwise its expected counts
tend to infinity. Each case therefore has a further subsequence along
which the normalized risk error tends to zero for every procedure and
for the Oracle. A subsequence bounded away from zero in absolute error
would contradict this property. Hence all normalized expansions hold
on the full sequence. Finally $1-Q_m\ge1-q_+>0$ converts the expansions
to the asserted risk ratios without requiring a limit of $Q_m$.
\end{proof}

\begingroup
\footnotesize
\setlength{\bibsep}{6pt}
\bibliography{ABOS_One_Factor_Dependence_References}
\endgroup
\end{document}